\RequirePackage{fix-cm}
\documentclass[smallextended]{svjour3}       % onecolumn (second format)
\smartqed  % flush right qed marks, e.g. at end of proof
\usepackage{graphicx}
\usepackage{mathptmx}      % use Times fonts if available on your TeX system
\usepackage{hyperref} % for autoref
\usepackage{amsmath} % for eqref
\usepackage{subcaption} % for subfigure
\usepackage{amsfonts} % for mathbb
\usepackage{bm} % for bm
\usepackage{bbm} % for indicator function
\usepackage{stmaryrd} % for \llbracket and \rrbracket 
\usepackage{makecell} % for multiple lines within a cell of table
\usepackage{multirow} % for multirow within a table
\usepackage[misc]{ifsym} % for symbol Letter
\usepackage{appendix} % for appendix
\usepackage{color} % for color
\usepackage{fullpage}

\usepackage{mathtools}
\DeclarePairedDelimiter\ceil{\lceil}{\rceil}

\usepackage{algorithm}
\usepackage{algorithmic}
\begin{document}

%%%%%%%%%%%%%%%%%%%%%%%%%%%%%%%%%%%%%%%%%%%%%%%%%%%%%
%%%%%%%%%%%%%%%%%%%%%%%%%%%%%%%%%%%%%%%%%%%%%%%%%%%%%
\title{Lift-and-Embed Learning Method for Solving Scalar Hyperbolic Equations with Discontinuous Solutions}

\titlerunning{Lift-and-Embed Learning Method for Solving Scalar Hyperbolic Equations} % if too long for running head

\author{Qi Sun \textsuperscript{{1,2}}
        \and Zhenjiang Liu \textsuperscript{1}        
        \and Lili Ju \textsuperscript{{3}}    
        \and Xuejun Xu \textsuperscript{{1,2}}        
        }

\authorrunning{Q. Sun, Z. Liu, X. Xu, and L. Ju} % if too long for running head

\institute{ 
 \begin{itemize}
  \item[\textsuperscript{\Letter}] {Qi Sun} \\ qsun\_irl@tongji.edu.cn\\
  \item[] Zhenjiang Liu \\ 2333705@tongji.edu.cn \\
  \item[] Lili Ju \\ ju@math.sc.edu \\  
  \item[] Xuejun Xu \\ xuxj@tongji.edu.cn 
  \at
  \item[\textsuperscript{1}] School of Mathematical Sciences, Tongji University, Shanghai 200092, China
  \item[\textsuperscript{2}] Key Laboratory of Intelligent Computing and Applications (Ministry of Education),  Tongji University, Shanghai 200092, China
  \item[\textsuperscript{3}] Department of Mathematics, University of South Carolina, Columbia, SC 29208, USA
 \end{itemize}
}

\date{Received: date / Accepted: date}

\maketitle
%%%%%%%%%%%%%%%%%%%%%%%%%%%%%%%%%%%%%%%%%%%%%%%%%%%%%
%%%%%%%%%%%%%%%%%%%%%%%%%%%%%%%%%%%%%%%%%%%%%%%%%%%%%

%%%%%%%%%%%%%%%%%%%%%%%%%%%%%%%%%%%%%%%%%%%%%%%%%%%%%
%%%%%%%%%%%%%%%%%%%%%%%%%%%%%%%%%%%%%%%%%%%%%%%%%%%%%
\begin{abstract}

Deep learning methods, which exploit auto-differentiation to compute derivatives without dispersion or dissipation errors, have recently emerged as a compelling alternative to classical mesh-based numerical schemes for solving hyperbolic conservation laws. However, solutions to hyperbolic problems are often piecewise smooth, posing challenges for training of neural networks to capture solution discontinuities and jumps across  interfaces. In this paper, we propose a novel lift-and-embed learning method to  effectively resolve these challenges. The proposed method comprises three innovative components: (i) embedding the Rankine-Hugoniot condition within a one-order higher-dimensional space by including an augmented variable; (ii) utilizing neural networks to handle the increased dimensionality and address both linear and nonlinear problems within a unified mesh-free learning framework; and (iii) projecting the trained model back onto the original physical domain to obtain the approximate solution. Notably, the location of discontinuities also can be treated as trainable parameters in our method and inferred concurrently with the training of neural network solutions. With collocation points sampled only on piecewise surfaces rather than fulfilling the whole lifted space, we demonstrate through extensive numerical experiments that our method can efficiently and accurately solve scalar hyperbolic equations with discontinuous solutions without spurious smearing or oscillations.

\end{abstract}
%%%%%%%%%%%%%%%%%%%%%%%%%%%%%%%%%%%%%%%%%%%%%%%%%%%%%
%%%%%%%%%%%%%%%%%%%%%%%%%%%%%%%%%%%%%%%%%%%%%%%%%%%%%

%%%%%%%%%%%%%%%%%%%%%%%%%%%%%%%%%%%%%%%%%%%%%%%%%%%%%
%%%%%%%%%%%%%%%%%%%%%%%%%%%%%%%%%%%%%%%%%%%%%%%%%%%%%
\section{Introduction}

The modeling of various problems in natural sciences and engineering is rooted in the concept of conservation laws, e.g., the Euler equations in fluid dynamics \cite{leveque2002finite}, the Lighthill-Whitham-Richards equations for traffic flow \cite{di2023physics}, and Maxwell’s equations of electromagnetics \cite{hesthaven2017numerical}, all of which are mathematically described as hyperbolic partial differential equations \cite{evans2022partial}. However, numerical simulation of these systems remains challenging, primarily due to discontinuities induced by non-smooth initial data or intersecting characteristics. Traditional numerical methods using uniform meshes, such as the upwind and Lax-Wendroff schemes \cite{gustafsson2013time}, suffer from dispersion or dissipation issues near discontinuities. A variety of refinements, including but not limited to ENO/WENO-based finite volume schemes \cite{shu1988efficient,harten1997uniformly,leveque2002finite,hesthaven2017numerical}, discontinuous/adaptive finite element methods \cite{cockburn2012discontinuous,godlewski2013numerical,demkowicz2010class,dahmen2012adaptive}, have achieved broad success in capturing sharp solution transitions. However, anomalous solutions, such as the carbuncle phenomenon \cite{peery1988blunt,kitamura2012carbuncle}, are observed when using flux limiters within high-resolution methods, necessitating specialized adjustments through ad hoc stabilization techniques. In addition to the escalating algorithmic complexity, accurately resolving discontinuities often requires high-quality meshes and anisotropic refinements, the generation of which could be time-consuming for complex interfaces and high-dimensional problems.

With the rapid advancement of hardware and software resources, artificial neural network-based deep learning methods have emerged as a prominent approach in contemporary scientific computing, particularly valued for their universal approximation capability \cite{hornik1989multilayer}, flexible meshless implementation \cite{raissi2019physics}, and effectiveness in addressing the curse of dimensionality \cite{barron1993universal}. One representative work is the physics-informed neural networks (PINNs) \cite{lagaris1998artificial,raissi2019physics,karniadakis2021physics}, which trains network models to minimize loss functions derived from the residual of differential equations at selected collocation points. Unlike numerical differentiation employed in solving hyperbolic equations \cite{godlewski2013numerical}, auto-differentiation within neural networks \cite{paszke2017automatic} enables precise derivative evaluation and has attracted significant attention for its ability to mitigate stability issues \cite{lorin2024non,arora2023deep,patel2022thermodynamically,sun2020convolution}. Unfortunately, standard network architectures often struggle to capture solution discontinuities and jumps arising from shock waves or contact discontinuities \cite{patel2022thermodynamically,wang2024understanding,fuks2020limitations,huang2023limitations}, spurring ongoing efforts to enhance their ability for resolving discontinuities.

For instance, discontinuous solutions can be artificially smoothed by incorporating an artificial viscosity term \cite{patel2022thermodynamically} into the hyperbolic system, but it would incur substantial costs due to the necessity of computing second-order derivatives through automatic differentiation. To maintain the simplicity of the equations, adaptive activation functions \cite{jagtap2020adaptive,della2023discontinuous} have been proposed to facilitate the use of non-differentiable activation functions for shock capturing, while another line of methods resort to the adaptive sampling of collocation points near regions of discontinuity \cite{mao2020physics,gao2023failure}. Recent studies suggest that the primary challenge in addressing discontinuous solutions lies in their rigorous mathematical definition as weak solutions \cite{lorin2024non}, which motivates the design of loss functions based on control volume formulations \cite{patel2022thermodynamically}, variational formulations \cite{de2024wpinns}, and some of their variants \cite{liu2024discontinuity}. To further enhance computational efficiency, domain decomposition strategies are integrated into the training of local network solutions as well as the interfaces that separate them \cite{lorin2024non}. Meanwhile, considerable efforts have been devoted to combining traditional numerical methods with deep learning techniques \cite{bar2019learning,morand2024deep,chen2021deep,feng2024hybrid}, aiming to harness the strengths of both paradigms. 

In contrast to prior methods, in this work we introduce a novel lift-and-embed learning framework for solving scalar hyperbolic equations with discontinuous solutions. First, the Rankine-Hugoniot relation \cite{leveque2002finite}, derived from the local weak form of conservation laws across discontinuities, is embedded in a one-order higher-dimensional space by augmenting the solution representation with an auxiliary variable. Notably, the augmented variable endows our solution ansatz with a piecewise smoothness property, thereby eliminating the need for multi-valued function approximations and enabling neural networks to easily reconstruct sharp discontinuities. Besides, our reformulated equations are defined on piecewise surfaces rather than being distributed over the entire lifted space, ensuring the generation process of collocation points remains unaffected despite the increased dimensionality. Second, neural networks are utilized to parametrize the unknown solution, which alleviate the curse of dimensionality and serve as universal function approximators for addressing both linear and quasi-linear problems. Upon the completion of training process \cite{karniadakis2021physics}, the approximate solution is obtained by projecting the trained network solution back onto the original physical domain. In particular, building on deep learning advances in inverse problems \cite{raissi2019physics}, we are also able to incorporate unknown discontinuity locations as additional trainable parameters, which then can be learned jointly with the neural network solution. Finally, experimental studies on a wide range of benchmark scalar hyperbolic  problems are reported, wherein all discontinuous solutions are well catched without the manifestation of spurious numerical smearing or oscillations.
	
The remainder of this work is organized as follows.  A brief introduction to hyperbolic equations is first provided in Subsection \ref{SubSection-Preliminaries}.
Section \ref{Section-RelatedWork} presents relevant literature on existing traditional and deep learning-based numerical methods for hyperbolic equations, followed by an exploration of our motivations from the perspective of function approximation. Section \ref{Section-Method} is devoted to detailing our lifted-and-embed learning method and  corresponding algorithms with identified or unknown discontinuity locations, supplemented by two benchmark problems to enhance clarity. Next, numerical experiments conducted on various problems are reported in Section \ref{Section-Experiment}, demonstrating the effectiveness and efficiency of our proposed methods in handling discontinuous solutions. Finally, a concluding remark and future work are given in Section \ref{Section-Conclusion}.

%%%%%%%%%%%%%%%%%%%%%%%%%%%%%%%%%%%%%%%%%%%%%%%%%%%%%
\subsection{Hyperbolic Conservation Laws}\label{SubSection-Preliminaries}

\begin{figure}[t!]
\centering
\begin{subfigure}[b]{0.24\textwidth}
\centering
\includegraphics[width=\textwidth]{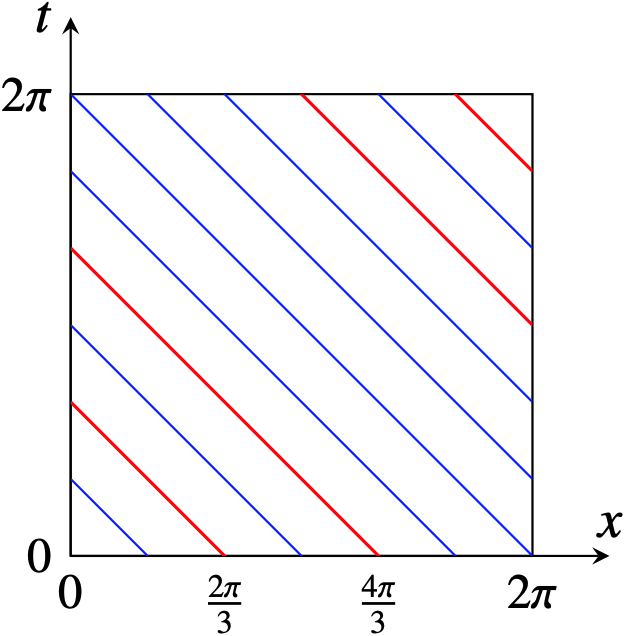}
\caption{Convection equation \eqref{Exp1-Eqns-u-exact}}
\label{fig-characteristics-linear-convection}
\end{subfigure}
\hspace{0.2cm}
\begin{subfigure}[b]{0.225\textwidth}
\centering
\includegraphics[width=\textwidth]{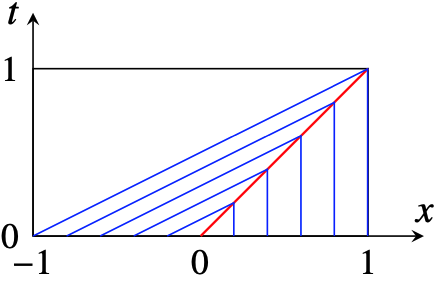}
\caption{Burgers' equation \eqref{Exp2-Eqns-u-exact}}
\label{fig-characteristics-inviscid-burgers-1}
\end{subfigure}
\hspace{0.2cm}
\begin{subfigure}[b]{0.2225\textwidth}
\centering
\includegraphics[width=\textwidth]{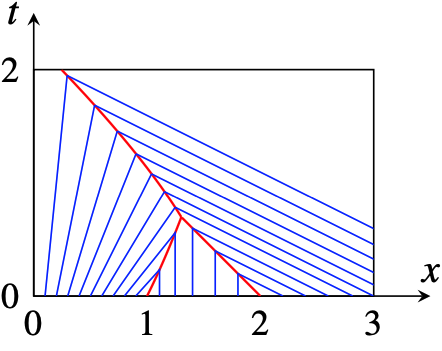}
\caption{Burgers' equation \eqref{Exp6-Eqns-u-exact}}
\label{fig-characteristics-inviscid-burgers-2}
\end{subfigure}
\hspace{0.2cm}
\begin{subfigure}[b]{0.24\textwidth}
\centering
\includegraphics[width=\textwidth]{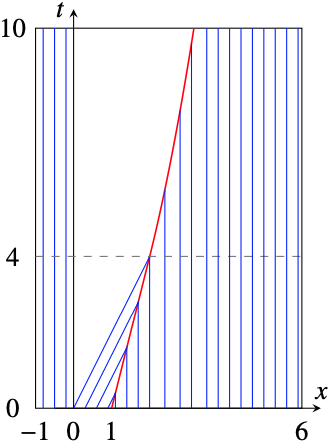}
\caption{Burgers' equation \eqref{ExpSM1-Eqns-u-exact}}
\label{fig-characteristics-inviscid-burgers-3}
\end{subfigure}
\vspace{-0.1cm}
\caption{Characteristic and discontinuity curves (marked in blue and red colors, respectively) for linear and quasi-linear hyperbolic equations.}
\label{fig-characteristics}
\vspace{-0.3cm}
\end{figure}

Scalar conservation laws arise from various physical and engineering disciplines, which are typically formulated as hyperbolic equations of the following form \cite{evans2022partial}
\begin{equation}
\left\{
\begin{array}{cl}
	\displaystyle \partial_t u(x,t) + \nabla_x \cdot f(u(x,t)) = 0, \ & \ \ \ \text{in}\ \ \Omega\times (0,T],\\
	u(x,0) = u_0(x), \ & \ \ \ \text{on}\ \ \Omega, \\
	u(x,t) = g(x,t), \ & \ \ \ \text{on}\ \ \Sigma \times (0, T],
\end{array}\right.
\label{GeneralProb-StrongForm}
\end{equation}
where $\Omega\subset \mathbb{R}^d$ $(d\in \mathbb{N}_+)$ is an open bounded domain with Lipschitz boundary $\partial\Omega$, $[0,T]$ the temporal interval, $u(x,t)\in\mathbb{R}^1$ a conserved quantity, $f(u) = (f_1(u),\cdots,f_d(u))$ the outward spatial flux vector, $\nabla_x\cdot$ a divergence operator with respect to the $x$-variable, $u_0(x)$ the initial data, and $g(x,t)$ the inflow data defined on parts of the boundary $\Sigma\subset\partial\Omega$ where characteristics enter the computational domain. Notably, classical solutions to the problem \eqref{GeneralProb-StrongForm} may not exist due to discontinuous initial data or intersecting characteristics (see \autoref{fig-characteristics}), thereby necessitating the introduction of weak solutions \cite{evans2022partial,godlewski2013numerical}. 
		
The concept of weak derivatives enables the analysis of solutions to the problem \eqref{GeneralProb-StrongForm} that exhibit only piecewise smoothness. Given an orientable surface of discontinuity $\Gamma$, we denote by $\bm{n} = (n_t, n_{x_1}, \cdots, n_{x_d})^T\neq (0,0,\cdots,0)^T$ a normal vector to it and by $u^+(x,t)$ and $u^-(x,t)$ the limits of $u(x,t)$ on each side of $\Gamma$, namely, for any $\varepsilon>0$,
\begin{equation*}
	u^{\pm}(x,t) = \lim_{\varepsilon\to0} u( (x,t)\pm\varepsilon\bm{n} ).
\end{equation*}
Note that if $(n_{x_1}, n_{x_2}, \cdots, n_{x_d}) \neq (0,0,\cdots,0)$, the normal vector can be rewritten as $\bm{n} = (-s, \bm{\nu})^T$ where $s\in\mathbb{R}$ and $\bm{\nu}=(\nu_1,\nu_2,\cdots,\nu_d)$ is a unit vector in $\mathbb{R}^d$ \cite{godlewski2013numerical}. By denoting $\llbracket \cdot \rrbracket$ the difference of quantity across the discontinuity surface, i.e., $\llbracket u \rrbracket = u^+ - u^-$ and $\llbracket f(u) \rrbracket = f(u^+) - f(u^-)$, the Rankine-Hugoniot relation is introduced below to ensure the consistency of conservation laws across discontinuities \cite{godlewski2013numerical}.
	
\begin{lemma} \textnormal{\cite{godlewski2013numerical}}
\label{Thm-Speed-Discontinuities}
Let $u\in L^\infty(\Omega\times [0,T] )$  be a piecewise $C^1$ function, then $u(x,t)$ is a weak solution of \eqref{GeneralProb-StrongForm} if and only if the two following conditions are satisfied:
\begin{enumerate}
\item[$(1)$] $u(x,t)$ is a classical solution of \eqref{GeneralProb-StrongForm} in the domains where $u(x,t)$ is $C^1$;
\item[$(2)$] $u(x,t)$ satisfies the jump condition along the surfaces of discontinuity, i.e.,
\begin{equation}    	
s\llbracket u \rrbracket = \sum_{j=1}^d \nu_j \llbracket f_j(u) \rrbracket
\label{RH-jump-condition}
\end{equation}    
where $s$ can be interpreted as the speed of propagation of the discontinuity.
\end{enumerate}
\end{lemma}

\begin{remark}
For a smooth discontinuity curve $\Gamma = (t, \gamma(t))$ in one-dimension, we have $\bm{\nu} = 1$ and $s(t) = \frac{d \gamma(t)}{dt}$.
\end{remark}

%--------------------------------%
%%%--------------------------------%

%%%%%%%%%%%%%%%%%%%%%%%%%%%%%%%%%%%%%%%%%%%%%%%%%%%%%
%%%%%%%%%%%%%%%%%%%%%%%%%%%%%%%%%%%%%%%%%%%%%%%%%%%%%

%%%%%%%%%%%%%%%%%%%%%%%%%%%%%%%%%%%%%%%%%%%%%%%%%%%%%
%%%%%%%%%%%%%%%%%%%%%%%%%%%%%%%%%%%%%%%%%%%%%%%%%%%%%
\section{Related Work and Motivations}\label{Section-RelatedWork}

%%%%%%%%%%%%%%%%%%%%%%%%%%%%%%%%%%%%%%%%%%%%%%%%%%%%%
\subsection{Related Work}

In contrast to the traditional mesh-based discretization of differential operators \cite{godlewski2013numerical}, deep learning methods leverage automatic differentiation to compute derivatives \cite{lagaris1998artificial,raissi2019physics}, providing a mesh-free alternative that can circumvent dispersion and dissipation issues encountered in the numerical simulation of hyperbolic problems. However, solutions to hyperbolic equations usually develop discontinuities over time and neural network methods often struggle to capture such solution jumps, necessitating specialized techniques to augment existing learning approaches \cite{karniadakis2021physics} for effective discontinuity resolution.

Instead of employing standard network architectures for solution approximation, adaptive activation functions \cite{jagtap2020adaptive,della2023discontinuous} and attention-based neural networks \cite{rodriguez2022physics} have been proposed to enhance the representation of discontinuous solutions, paired with adaptive sampling strategies \cite{mao2020physics,gao2023failure} that prioritize training data allocation near regions of discontinuity. On the other hand, recent studies suggest that the primary challenge in addressing discontinuous solutions arises from their mathematical characterization as weak solutions, hence independently of the regularity of the activation function being implemented \cite{lorin2024non}. Nevertheless, standard neural networks remain effective for computing rarefaction waves \cite{lorin2024non} that constitute continuous solutions.
Motivated by the method of vanishing viscosity \cite{evans2022partial}, an alternative approach involves introducing a suitable diffusion term to mitigate the problem's hyperbolicity \cite{coutinho2023physics,patel2022thermodynamically}. While this technique regularizes the solution near discontinuities and enables seamless auto-differentiation across the entire computational domain, it converts the original problem into a parabolic type, which can often result in a degradation in the solution accuracy. Additionally, computing second-order derivatives can significantly slow down the training process, as each spatial point requires two rounds of backpropagation instead of a single backward pass for constructing the residual loss function \cite{diab2021pinns}. 

Recall that a key feature of hyperbolic problems lies in their tendency to develop discontinuous solutions within a finite time, which indicates that the underlying equations should be considered in a weak sense rather than in a pointwise (strong) manner. Accordingly, the weak formulation of hyperbolic equations is recast as a min-max optimization problem \cite{de2024wpinns}, where the network solution is maximized with respect to test functions and minimized with respect to trial functions. This approach simplifies to a minimization problem when adopting a mesh-based representation for test functions \cite{kharazmi2019variational}, though it may compromise the meshless advantage offered by neural networks \cite{de2024wpinns}. Alternatively, the weak formulation can be characterized locally by the Rankine-Hugoniot condition \cite{leveque2002finite}, which allows for the design of adaptive weights that can autonomously switch between strong and weak formulation-based loss functions \cite{liu2024discontinuity}. Moreover, the computation of model solutions could be parallelized using domain decomposition strategies \cite{lorin2024non}, with discontinuity interfaces co-evolved during training.

By partitioning the space-time domain into disjoint subregions and applying the Gauss divergence theorem \cite{evans2022partial} to each, hyperbolic equations can also be expressed in an integral form, referred to as the control volume formulation \cite{leveque2002finite}, which admits weak solutions involving discontinuities and therefore facilitates the development of various learning algorithms \cite{patel2022thermodynamically,cen2024deep,cai2023least,cai2022least,cai2023evolving}. Concurrently, much effort has been made to incorporating deep learning techniques into existing numerical methods, such as employing neural networks to learn, either partially or entirely, the finite volume \cite{bar2019learning,bois2023optimal,morand2024deep} or discontinuous Galerkin schemes \cite{chen2021deep} with specialized numerical flux functions. Encoding the distinctive structure of hyperbolic equations \cite{mojgani2023kolmogorov,zhou2024capturing}, such as characteristic curves \cite{braga2022characteristics} or implicit forms \cite{zhang2022implicit}, within the network architecture offers another way of integrating prior domain knowledge.

%%%%%%%%%%%%%%%%%%%%%%%%%%%%%%%%%%%%%%%%%%%%%%%%%%%%%
\subsection{Motivations from A Function Approximation Perspective}

%--------------------------------%
\begin{figure}[!htbp]
\centering
\begin{subfigure}[b]{0.29\textwidth}
\centering
\includegraphics[width=\textwidth]{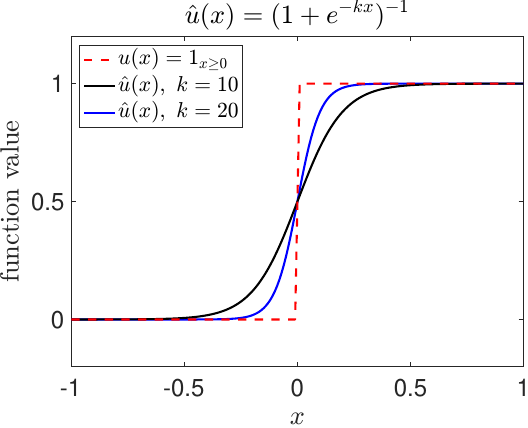}
\caption{Smooth approximation}
\label{fig-heaviside-approx-logistic}
\end{subfigure}
\hfill
\begin{subfigure}[b]{0.3\textwidth}
\centering
\includegraphics[width=\textwidth]{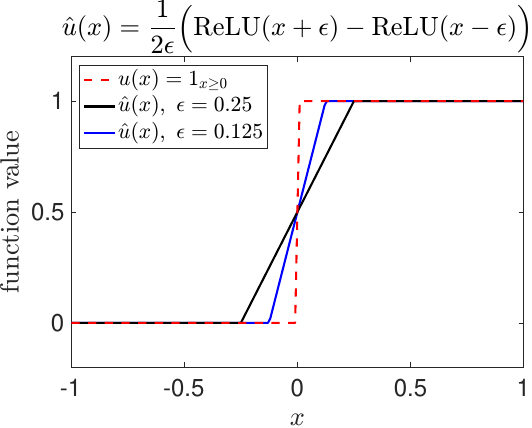}
\caption{Piecewise approximation}
\label{fig-heaviside-approx-pwlinear}
\end{subfigure}
\hfill
\begin{subfigure}[b]{0.33\textwidth}
\centering
\includegraphics[width=\textwidth]{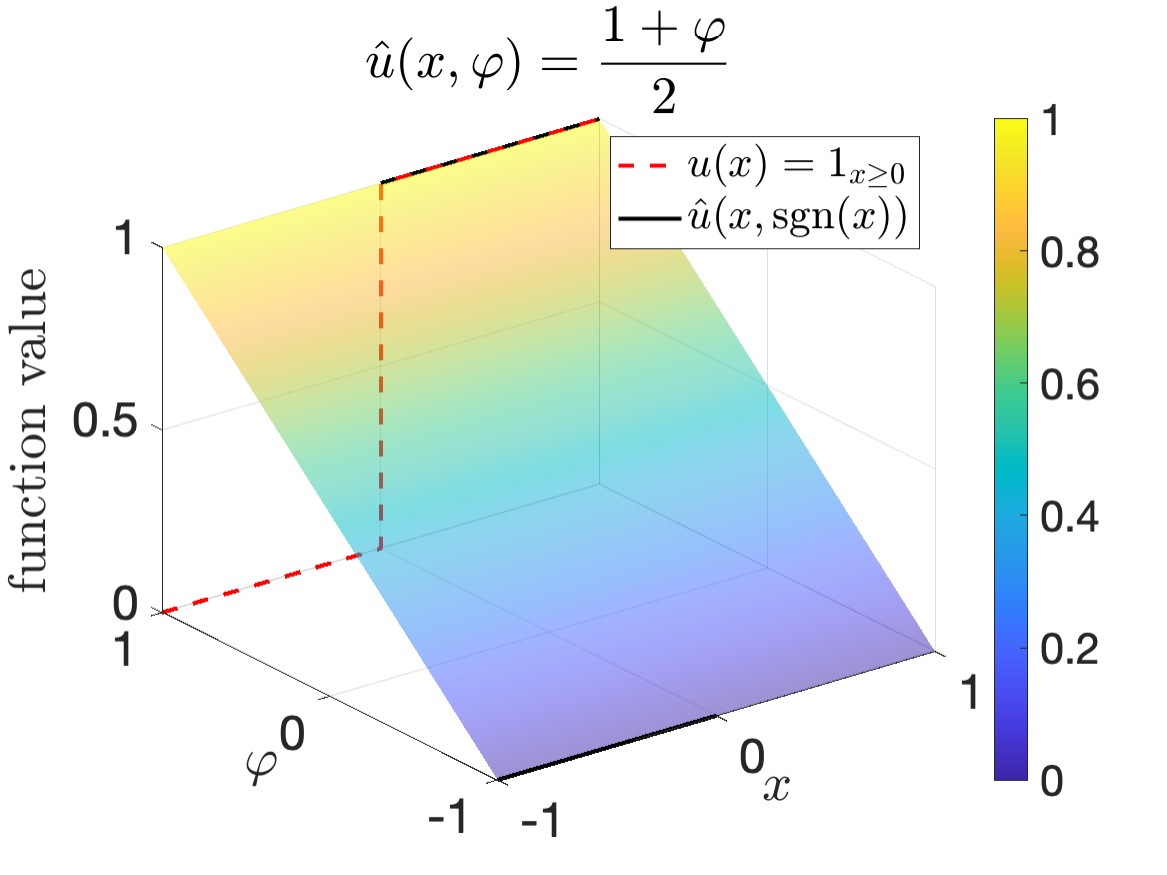}
\caption{Lift-and-embed reconstruction}
\label{fig-heaviside-approx-LP}
\end{subfigure}
\vspace{-0.1cm}
\caption{Approaches for approximating the Heaviside function: (a) a kernel smoothing, (b) a finite element interpolation, (c) a lift-and-embed technique. }
\label{fig-heaviside-approx}
\vspace{-0.3cm}
\end{figure}
%--------------------------------%

To elucidate our motivations in handling discontinuities, we present an illustrative example involving the analytic approximation of Heaviside step function $u(x) = H(x) := \mathbbm{1}_{x\geq 0}$ using kernel smoothing, finite element interpolation, and lift-and-embed technique, as depicted in \autoref{fig-heaviside-approx}. More precisely, a common approach is to apply the convolution operation with suitable kernel functions, e.g.,  
\begin{equation*}
	u(x) \approx  \hat{u}(x)  := ( u * G )(x) = \int_{-\infty}^{+\infty} u(\tau)G(x-\tau) \, d\tau = \frac{1}{1+e^{-kx}},
\end{equation*}
where $G(y) = ke^{-ky}(1+e^{-ky})^{-2}$ is the Sigmoid kernel function. As is well known, a larger value of parameter $k>0$ corresponds to a sharper transition at $x=0$ (see \autoref{fig-heaviside-approx-logistic}), in which the discontinuity is smoothed out in a manner analogous to the method of vanishing viscosity applied to hyperbolic equations \cite{patel2022thermodynamically}.

Alternatively, consider the partition $-1 = x_0 < x_1 < x_2 < x_3 = 1$ of interval $[-1,1]$ with $x_2 = - x_1 = \epsilon$, another widely used approach is to approximate the target function using piecewise linear elements (see \autoref{fig-heaviside-approx-pwlinear}), i.e.,
\begin{equation*}
	u(x) \approx  \hat{u}(x) := \sum_{j=0}^3 u(x_j)\phi_j(x) = \frac{1}{2\epsilon} \Big( \textnormal{ReLU}(x+\epsilon) - \textnormal{ReLU}(x-\epsilon) \Big),
\end{equation*}
where basis functions $\{\phi_j(x)\}_{j=0}^3$ are rewritten in terms of the rectified linear unit function $\textnormal{ReLU}(x)=\max(0,x)$ \cite{xu2020finite}. The latter expression is also known as a free-knot spline with two knots or as a single-hidden-layer feedforward neural network \cite{cai2023evolving}.

Beyond these two representative methods, discontinuous functions can be embedded in a one-order higher-dimensional space to achieve continuity \cite{hu2022discontinuity,tseng2023cusp}. The key idea involves including an augmented variable to allocate the jump in $u(x) = \mathbbm{1}_{x\geq 0}$ to separate branches, and a concrete example is displayed in \autoref{fig-heaviside-approx-LP}, namely,
\begin{equation*}
	u(x) =\hat{u}( x, \textnormal{sgn}(x)) \ \ \ \textnormal{with} \ \ \ \hat{u}(x,\varphi) := \frac{1-\varphi}{2}u^- +\frac{1+\varphi}{2}u^+=\frac{1+\varphi}{2},
\end{equation*}
where $u^- =0$ and $u^+ =1$. Clearly, function $\hat{u}(x,\varphi)$ is no longer discontinuous with respect to the $x$-variable and admits exact recovery by choosing $\varphi = \textnormal{sgn}(x) = 2 (\mathbbm{1}_{x\geq 0} - 0.5)$. While the lifting operation often yields a high-dimensional function that may incur the curse of dimensionality, neural networks are highly effective at addressing this issue \cite{barron1993universal}. Moreover, discontinuity locations can be regarded as unknown parameters to be inferred during training \cite{lorin2024non}, which deep learning-based methods are proven to be well-suited to resolve \cite{raissi2019physics}. 	

%%%%%%%%%%%%%%%%%%%%%%%%%%%%%%%%%%%%%%%%%%%%%%%%%%%%%
%%%%%%%%%%%%%%%%%%%%%%%%%%%%%%%%%%%%%%%%%%%%%%%%%%%%%

%%%%%%%%%%%%%%%%%%%%%%%%%%%%%%%%%%%%%%%%%%%%%%%%%%%%%
%%%%%%%%%%%%%%%%%%%%%%%%%%%%%%%%%%%%%%%%%%%%%%%%%%%%%
\section{The Proposed Lift-and-Embed Learning Method for Scalar Hyperbolic Equations}\label{Section-Method}

In this section, we present a lift-and-embed learning method for solving scalar hyperbolic equations with discontinuous solutions, addressing scenarios with and without a-priori knowledge of discontinuity locations. We begin by embedding the target  equations into a one-order higher-dimensional space, illustrated with two benchmark examples to show the effectiveness in resolving discontinuities while maintaining overall continuity. Next, detailed algorithms are presented, as well as the parametrization of discontinuity interfaces as trainable parameters, which allows us to harness the flexibility of the proposed deep learning method in solving both forward and inverse problems.

%%%%%%%%%%%%%%%%%%%%%%%%%%%%%%%%%%%%%%%%%%%%%%%%%%%%%
\subsection{Hyperbolic Equations Embedded in  A One-order Higher-Dimensional Space}

Here, we focus on a single discontinuity interface that divides the space-time domain into two subregions, which can be generalized to scenarios containing multiple interfaces. Based on the parametrization of interface 
\begin{equation*}
	\Gamma = \Gamma(t, x_1, x_2, \cdots, x_d)\ \ \ \textnormal{or}\ \ \ x_d = \gamma(t, x_1, \cdots, x_{d-1}),
\end{equation*}
we construct the solution ansatz to \eqref{GeneralProb-StrongForm} by incorporating an augmented variable 
\begin{equation*}
	u(x,t) = \hat{u}(x,t,\varphi(x,t))\ \ \ \textnormal{with}\ \ \ \varphi(x,t) = H(x_d - \gamma(t, x_1, \cdots, x_{d-1})),
\end{equation*}
where $H(x) = \mathbbm{1}_{x\geq 0}$ is the Heaviside step function. Though the inclusion of an extra augmented variable appears to increase the dimensionality, its dependence on intrinsic variables $(x,t)$ remains unchanged.
	
Let $\delta(*)$ denote the Dirac delta function, a direct calculation implies that
\begingroup
\renewcommand*{\arraystretch}{1.2}
\begin{equation*}
	\begin{array}{c}
		\displaystyle \partial_t u = \partial_t \hat{u} - \delta(x_d - \gamma) \, \partial_\varphi \hat{u} \, \partial_t \gamma,\\
		\displaystyle \nabla_x \cdot f(u) = f_d'(\hat{u}) (\partial_{x_d} \hat{u} + \delta(x_d-\gamma) \, \partial_\varphi\hat{u})  +  \sum_{j=1}^{d-1} f_j'(\hat{u}) (\partial_{x_j}\hat{u} - \delta(x_d -\gamma) \, \partial_\varphi\hat{u} \, \partial_{x_j} \gamma ),\\
	\end{array}
\end{equation*}
\endgroup
while the Rankine-Hugoniot jump condition \eqref{RH-jump-condition} can be reformulated as
\begin{equation*}    	
   s\big(\hat{u}(x,t,\varphi^+(x,t))-\hat{u}(x,t,\varphi^-(x,t)) \big) = \sum_{j=1}^d \nu_j \big(f_j \big(\hat{u}(x,t,\varphi^+(x,t)) \big) - f_j \big(\hat{u}(x,t,\varphi^-(x,t)) \big) \big) 
\end{equation*} 
along the surface of discontinuity. Then, the problem \eqref{GeneralProb-StrongForm} is embedded into a one-order higher-dimensional space as follows:
\begin{equation}
\left\{
\begin{array}{ll}
\displaystyle \partial_t \hat{u}(x,t,\varphi(x,t)) + \nabla_x \cdot f(\hat{u}(x,t,\varphi(x,t)))  = 0, \ & \ \ \ \textnormal{for}\ \ (x, t) \in \Omega\times (0, T] \setminus \Gamma,\\
\displaystyle s\llbracket \hat{u}(x,t,\varphi(x,t)) \rrbracket = \sum_{j=1}^d \nu_j \llbracket f_j(\hat{u}(x,t,\varphi(x,t))) \rrbracket, \ & \ \ \ \textnormal{for}\ \ (x, t) \in \Gamma, \\
\hat{u}(x, 0, \varphi(x,0)) = u_0(x), \ & \ \ \ \textnormal{for}\ \ (x,t) \in \Omega\times\{0\},\\
\hat{u}(x, t, \varphi(x,t)) = g(x, t), \ & \ \ \ \textnormal{for}\ \ (x, t) \in \Sigma \times (0, T],\\
\end{array}\right.
\label{GeneralProb-Embeded-Form}
\end{equation}
in which the solution jump is imposed on separate branches that greatly differs from its lower-dimensional counterpart \eqref{RH-jump-condition}. Before detailing our learning algorithm, two benchmark examples with analytically defined interfaces of discontinuity are provided to illustrate the effectiveness of equations \eqref{GeneralProb-Embeded-Form} in handling discontinuous solutions.
	
%%%%%%%%%%%%%%%%%%%%%%%%%%%%%%%%%%%%%%%%%%%%%%%%%%%%%
\subsubsection{Example 1: A Linear Convection Equation in One-Dimension}\label{sec-Exp1-1d-Linear-disCsolu-a1}

%---------------------------------------%
\begin{figure}[!htbp]
\centering
\begin{subfigure}{0.32\textwidth}
\centering
\includegraphics[width=\textwidth]{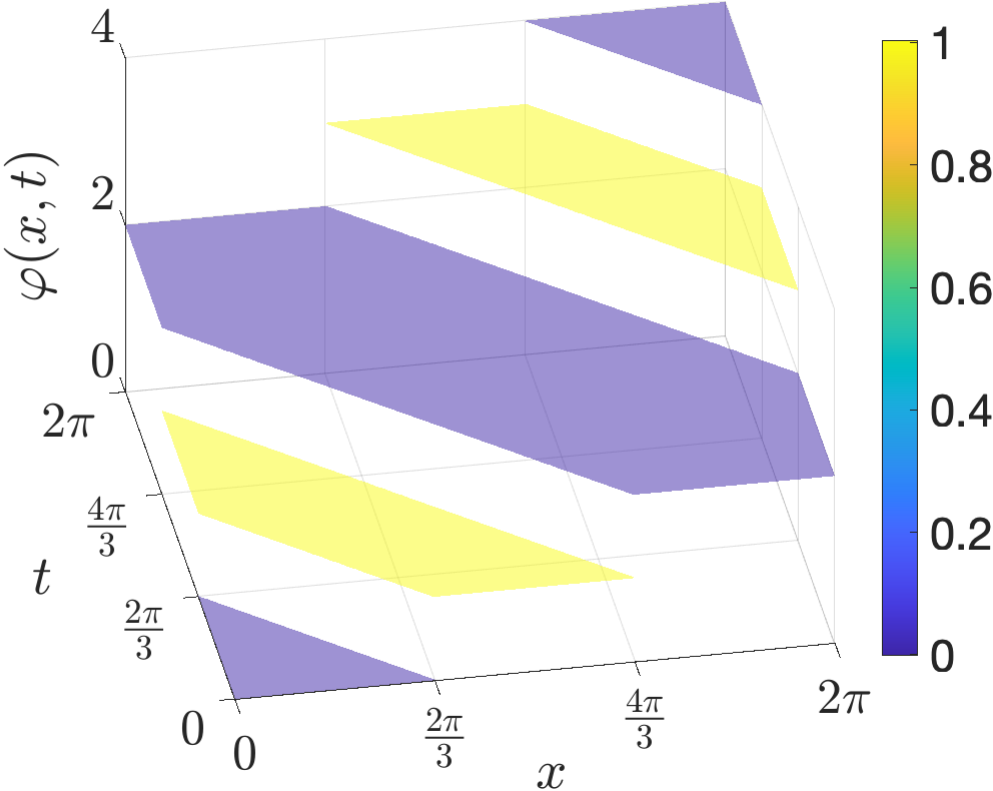}
\caption{$\hat{u}(x,t,\varphi(x,t))$}
\label{Exp1-fig-a}
\end{subfigure}
\hspace{0.15cm}
\begin{subfigure}{0.32\textwidth}
\centering
\includegraphics[width=0.92\textwidth]{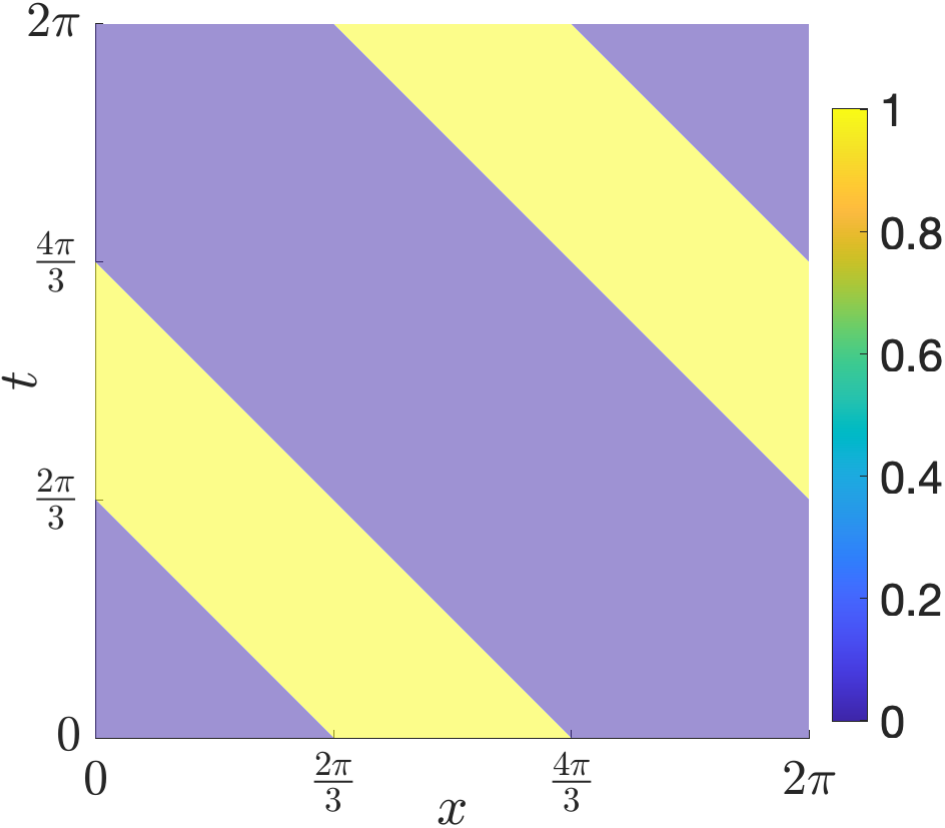}
\caption{$\check{u}(x,t)$ }
\end{subfigure}
\hspace{0.15cm}
\begin{subfigure}{0.32\textwidth}
\centering
\includegraphics[width=0.97\textwidth]{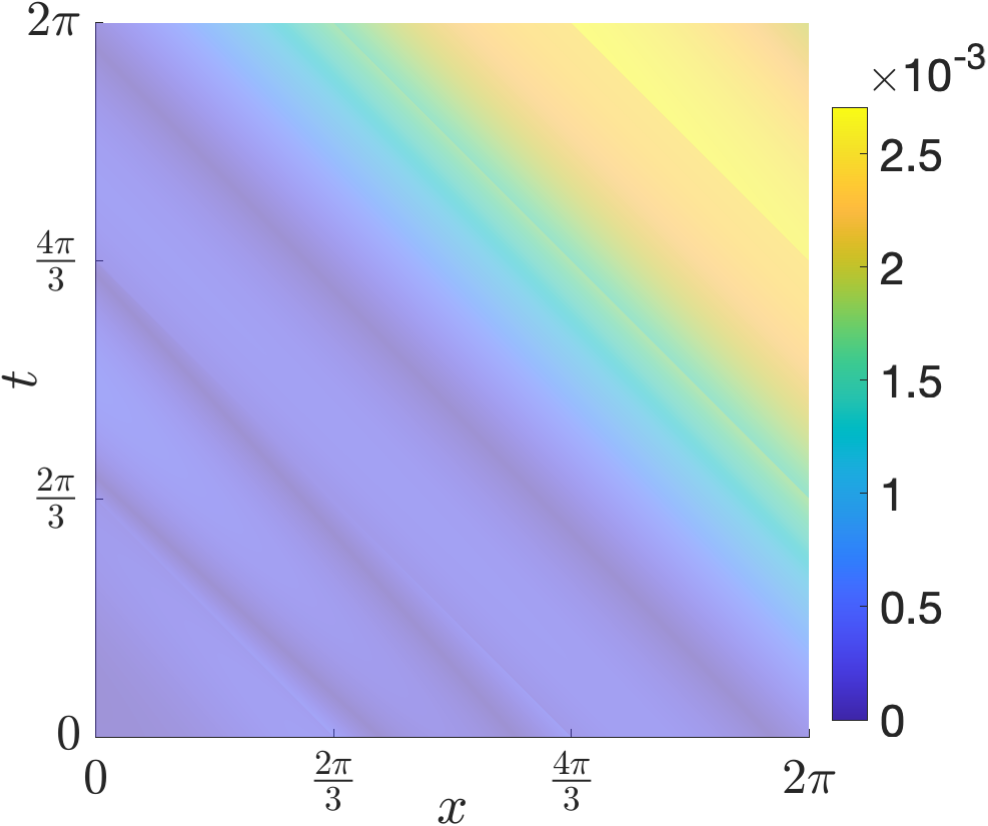}
\caption{$|u(x,t) - \check{u}(x,t)|$}
\end{subfigure}
\vspace{0.3cm}

\centering
\begin{subfigure}{0.32\textwidth}
\centering
\includegraphics[width=\textwidth]{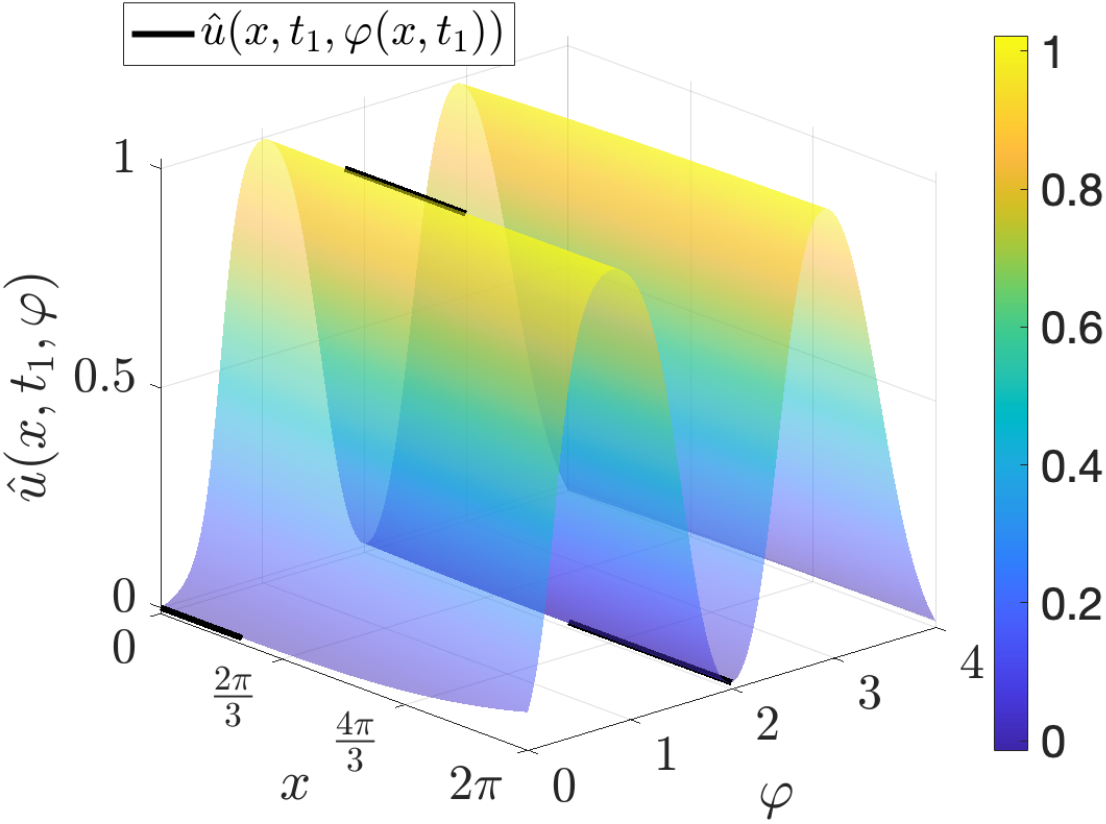}
\caption{$\hat{u}(x,t_1,\varphi)$}
\end{subfigure}
\hspace{0.15cm}
\begin{subfigure}{0.32\textwidth}
\centering
\includegraphics[width=0.9\textwidth]{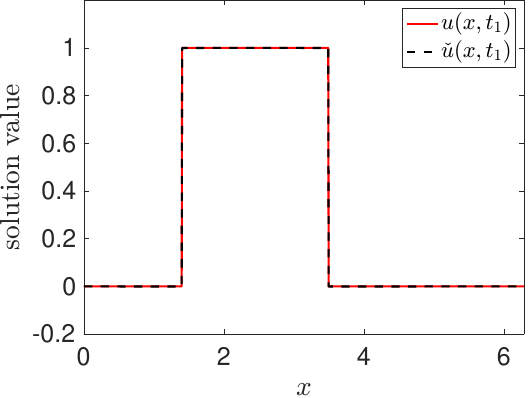}
\caption{$\check{u}(x,t_1)$ and $u(x,t_1)$ }
\end{subfigure}
\hspace{0.15cm}
\begin{subfigure}{0.32\textwidth}
\centering
\includegraphics[width=0.9\textwidth]{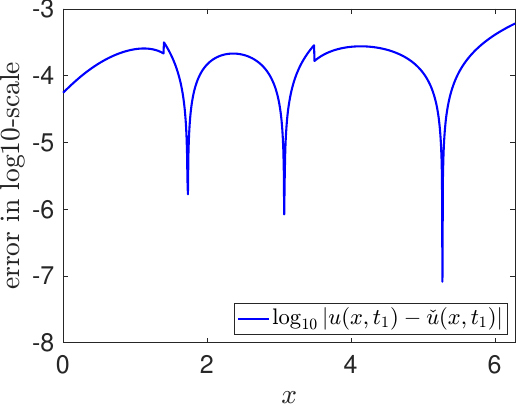}
\caption{$\log_{10}|u(x,t_1)-\check{u}(x,t_1)|$}
\end{subfigure}
\vspace{0.3cm}

\centering
\begin{subfigure}{0.32\textwidth}
\centering
\includegraphics[width=\textwidth]{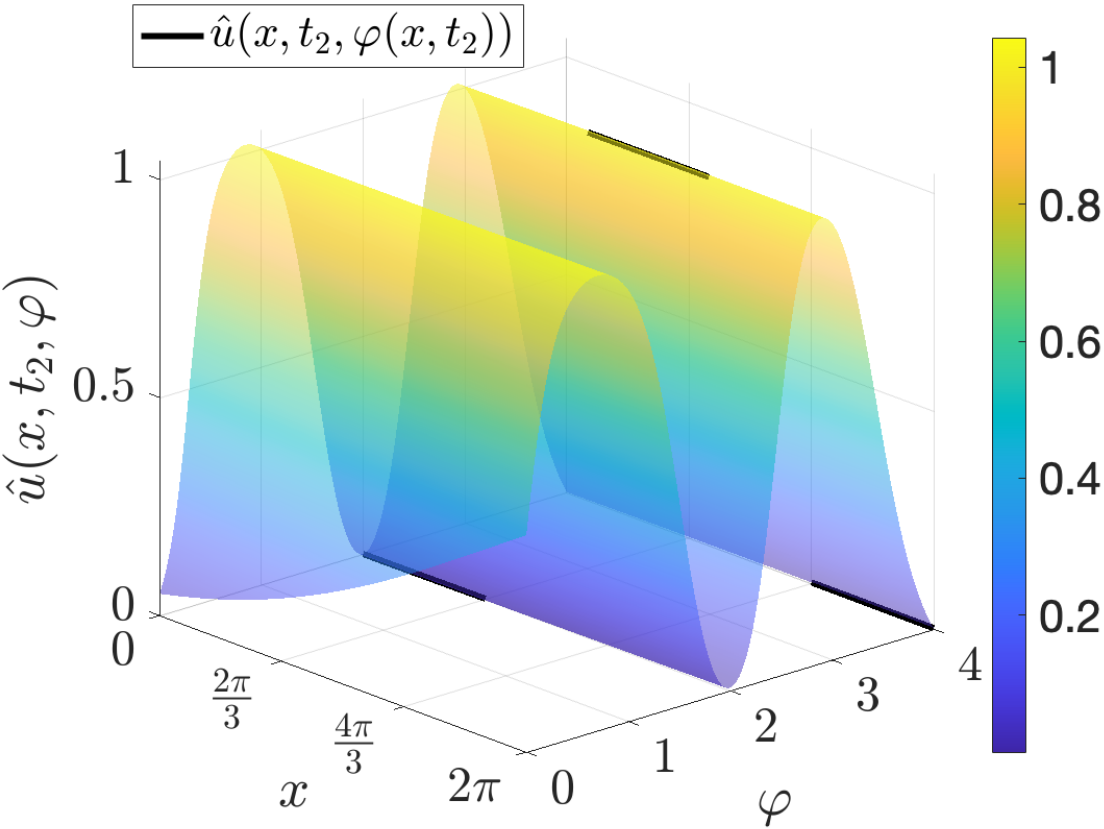}
\caption{$\hat{u}(x,t_2,\varphi)$}
\end{subfigure}
\hspace{0.15cm}
\begin{subfigure}{0.32\textwidth}
\centering
\includegraphics[width=0.9\textwidth]{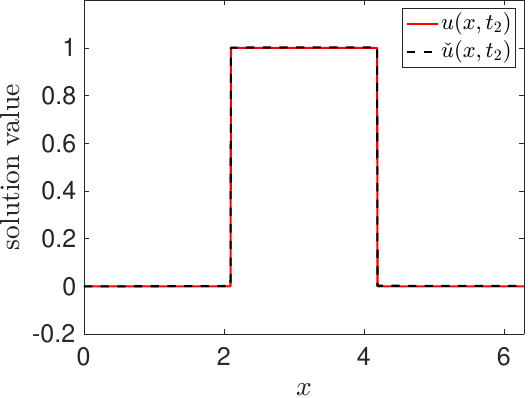}
\caption{$\check{u}(x,t_2)$ and $u(x,t_2)$ }
\end{subfigure}
\hspace{0.15cm}
\begin{subfigure}{0.32\textwidth}
\centering
\includegraphics[width=0.9\textwidth]{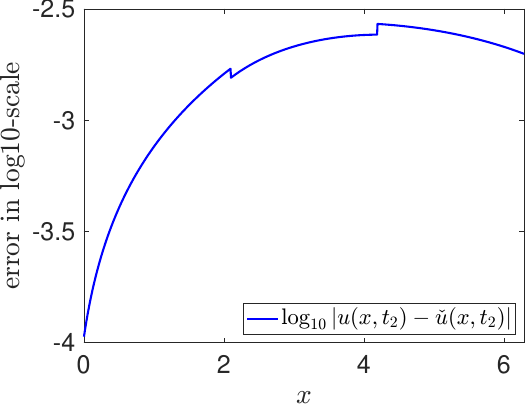}
\caption{$\log_{10}|u(x,t_2) - \check{u}(x,t_2)|$}
\end{subfigure}
\caption{Numerical results for the linear convection equation \eqref{Exp1-Eqns-u-exact} ($t_1 = \frac{2\pi}{9}$, $ t_2 = 2\pi$).}
\label{Exp1-fig-1d-Linear-disCsolu-a1}
\vspace{-0.3cm}
\end{figure}
%---------------------------------------%

Let us consider the following  one-dimensional convection equation \cite{gustafsson2013time}:
\begingroup
\renewcommand*{\arraystretch}{1.2}
\begin{equation}
\left\{
\begin{array}{ll}
\partial_t u(x,t) - \partial_x u(x,t) = 0, \ & \ \ \ \textnormal{for}\ \ (x,t)\in \Omega = (0,2\pi)\times (0,2\pi],\\
u_0(x) = H(x - \frac{2\pi}{3}) - H(x - \frac{4\pi}{3}), \ & \ \ \ \textnormal{for}\ \ x\in (0,2\pi), \\
u(0,t) = u(2\pi, t),\ & \ \ \ \textnormal{for}\ \ t\in (0,2\pi],
\end{array}\right.
\label{Exp1-Eqns-u-exact}
\end{equation}
\endgroup
where the exact solution $u(x,t) = u_0(x+t)$. Notably, the above initial data is given by a square wave, namely,
\begingroup
\renewcommand*{\arraystretch}{1.2}
\begin{equation}
u_0(x)=\left\{
\begin{array}{ll}
0, \ &\ \ \ \textnormal{for}\ \ 0\leq x < \frac{2\pi}{3},\\
1, \ &\ \ \ \textnormal{for}\ \ \frac{2}{3}\pi \leq x < \frac{4\pi}{3}, \\
0, \ &\ \ \ \textnormal{for}\ \ \frac{4}{3}\pi \leq x \leq 2\pi,
\end{array}\right.
\label{Exp1-initial-condition}
\end{equation}
\endgroup
where traditional mesh-based numerical schemes \cite{gustafsson2013time} often encounter issues with dissipation or dispersion near discontinuities  as the solution evolves over time.	

Note that, by Lemma \ref{Thm-Speed-Discontinuities}, the propagation speed of discontinuities for the problem \eqref{Exp1-Eqns-u-exact}
\begin{equation*}
	s = \frac{ \llbracket f(u) \rrbracket }{ \llbracket u \rrbracket }= -1\ \ \ \textnormal{where}\ \ f(u) = - u \ \ \textnormal{and}\ \ \bm{\nu}=1,
\end{equation*}
is the same as the slope of characteristic lines (see \autoref{fig-characteristics-linear-convection}). The augmented variable can then be constructed to embed the problem \eqref{Exp1-Eqns-u-exact} into a one-order higher-dimensional space  (see \autoref{Exp1-fig-a}) as follows:
\begin{equation}
	u(x,t) = \hat{u}(x,t,\varphi(x,t))\ \ \ \textnormal{with} \ \ \ \varphi(x,t) = \sum_{i=1}^2 \sum_{k=0}^{n_i}  H(x - st - x_i - 2k\pi),
	\label{Exp1-solution-ansatz}
\end{equation}
where $x_1=\frac{2\pi}{3}$, $x_2=\frac{4\pi}{3}$, $T=2\pi$, and $n_i = \lceil{ - \frac{sT+x_i}{2\pi} }\rceil$ for $s<0$ and $i=1,2$. Moreover, the solution jump arising from the initial data \eqref{Exp1-initial-condition} implies that
\begin{equation}
	\llbracket u_0(x_i) \rrbracket = \llbracket \hat{u}(x,t,\varphi(x,t)) \rrbracket = \hat{u}(x,t,\varphi^+(x,t)) - \hat{u}(x,t,\varphi^-(x,t)) \ \ \ \textnormal{on}\ \ \Gamma_{i} = \{ ( x(t), t ) \; |\; x(t) = s t + x_i + 2k\pi \}
\label{Exp1-jump-condition}
\end{equation}
where $0\leq k\leq n_i$ and $1\leq i\leq 2$. 

It is noteworthy that the right and left limits of our augmented variable $\varphi(x,t)$ at $\Gamma = \Gamma_1\cup\Gamma_2$ can take distinct values in \eqref{Exp1-jump-condition}, allowing the jump of solution $u(x,t)$ across the discontinuity interface to be exactly captured by the difference of ansatz $\hat{u}(x,t,\varphi(x,t))$ at separate points. On the other hand, a direct computation implies
\begin{equation*}
\partial_t u = \partial_t \hat{u} - s\partial_\varphi \hat{u} \sum_{i=1}^2 \sum_{k=0}^{n_i} \delta( x-st - x_i - 2k\pi) \ \ \textnormal{and} \ \ \partial_x u = \partial_x \hat{u}+ \partial_\varphi \hat{u} \sum_{i=1}^2 \sum_{k=0}^{n_i} \delta( x-st - x_i - 2k\pi).
\end{equation*}
	
To sum up, the hyperbolic system satisfied by our solution ansatz $\hat{u}(x,t,\varphi(x,t))$ takes on the form
\begingroup
\renewcommand*{\arraystretch}{1.2}
\begin{equation}
	\left\{
	\begin{array}{ll}
		\partial_t \hat{u}(x,t,\varphi(x,t)) - \partial_x \hat{u}(x,t,\varphi(x,t)) = 0, \ &\ \ \ \textnormal{for}\ \ (x,t)\in \Omega \setminus \Gamma_1\cup\Gamma_2,\\
		\llbracket \hat{u}(x,t,\varphi(x,t)) \rrbracket  = \llbracket u_0(x_i) \rrbracket, \ &\ \ \ \textnormal{for}\ \ (x,t)\in \Gamma_i \ (i=1,2),\\
		\hat{u}(x,0,\varphi(x,0)) = u_0(x), \ &\ \ \ \textnormal{for}\ \ x\in (0,2\pi), \\
		\hat{u}(0,t,\varphi(x,0)) = \hat{u}(2\pi, t, \varphi(2\pi, t)), \ &\ \ \ \textnormal{for}\ \ t\in (0,2\pi],
	\end{array}\right.
	\label{Eqns-Exp1-u-NN}
\end{equation}
\endgroup
which enables a smooth representation of the discontinuous solution problem (\ref{Exp1-Eqns-u-exact}) by bridging the jump conditions \eqref{Exp1-jump-condition} within a one-order higher-dimensional space (see (d) and (g) in \autoref{Exp1-fig-1d-Linear-disCsolu-a1}). The solution ansatz is then parametrized using a fully-connected neural network \cite{goodfellow2016deep}, followed by a vanilla training process under the pointwise residual-minimization framework \cite{raissi2019physics,karniadakis2021physics}. The trained neural network solution (still denoted as $\hat{u}(x,t,\varphi(x,t))$) and its projection back to the lower-dimensional space (denoted as $\check{u}(x,t)$) are shown in  (b), (e), and (h) of \autoref{Exp1-fig-1d-Linear-disCsolu-a1}, together with their error profiles without notable dissipation or dispersion errors (see (c), (f), and (i) in \autoref{Exp1-fig-1d-Linear-disCsolu-a1}).

%%%%%%%%%%%%%%%%%%%%%%%%%%%%%%%%%%%%%%%%%%%%%%%%%%%%%
\subsubsection{Example 2: An Inviscid Burgers’ Equation in One-Dimension}\label{sec-Exp2-1d-Burgers}

%---------------------------------------%
\begin{figure}[b!]
\centering
\begin{subfigure}{0.32\textwidth}
\centering
\includegraphics[width=\textwidth]{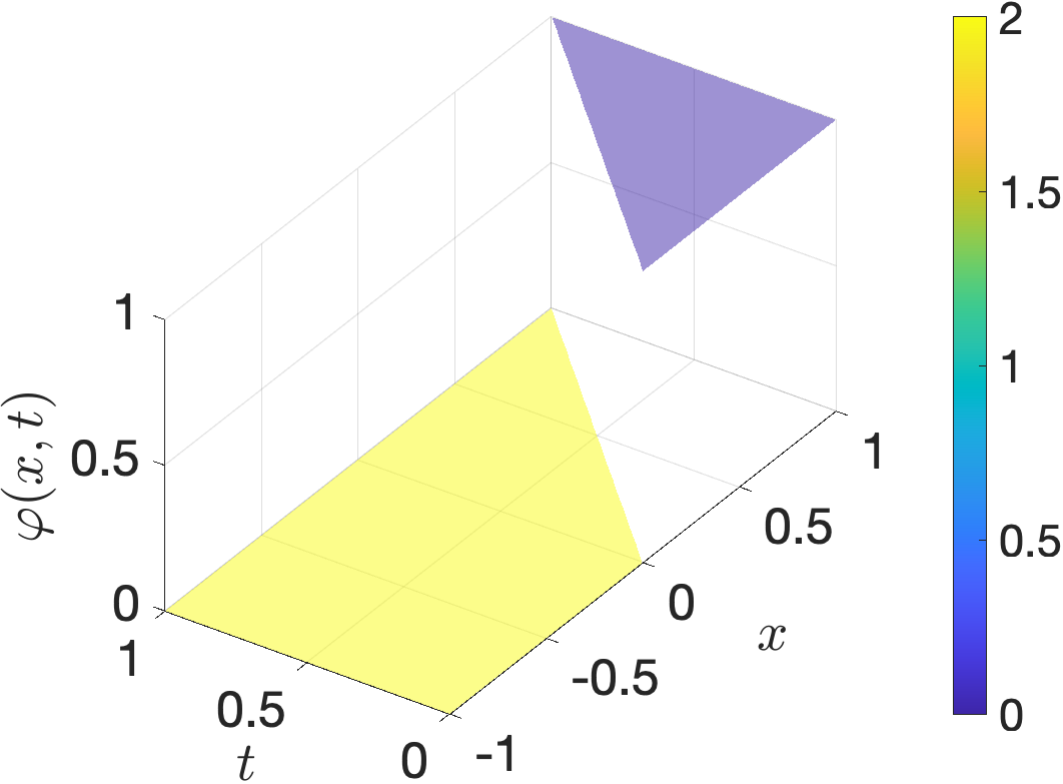}
\caption{$\hat{u}(x,t,\varphi(x,t))$}
\label{Exp2-fig-a}
\end{subfigure}
\hspace{0.15cm}
\begin{subfigure}{0.32\textwidth}
\centering
\includegraphics[width=0.92\textwidth]{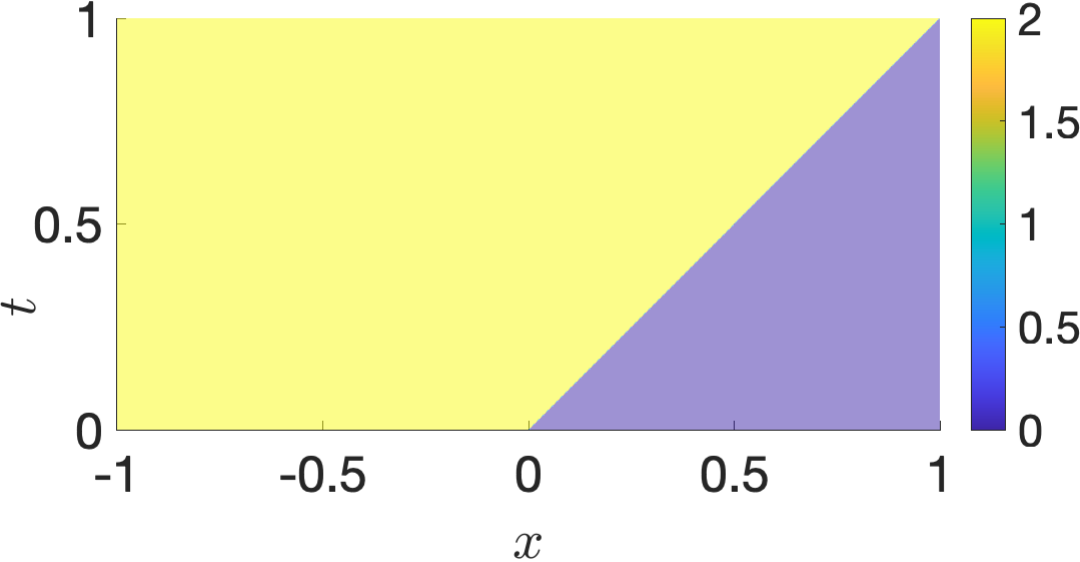}
\caption{$\check{u}(x,t)$ }
\end{subfigure}
\hspace{0.15cm}
\begin{subfigure}{0.32\textwidth}
\centering
\includegraphics[width=\textwidth]{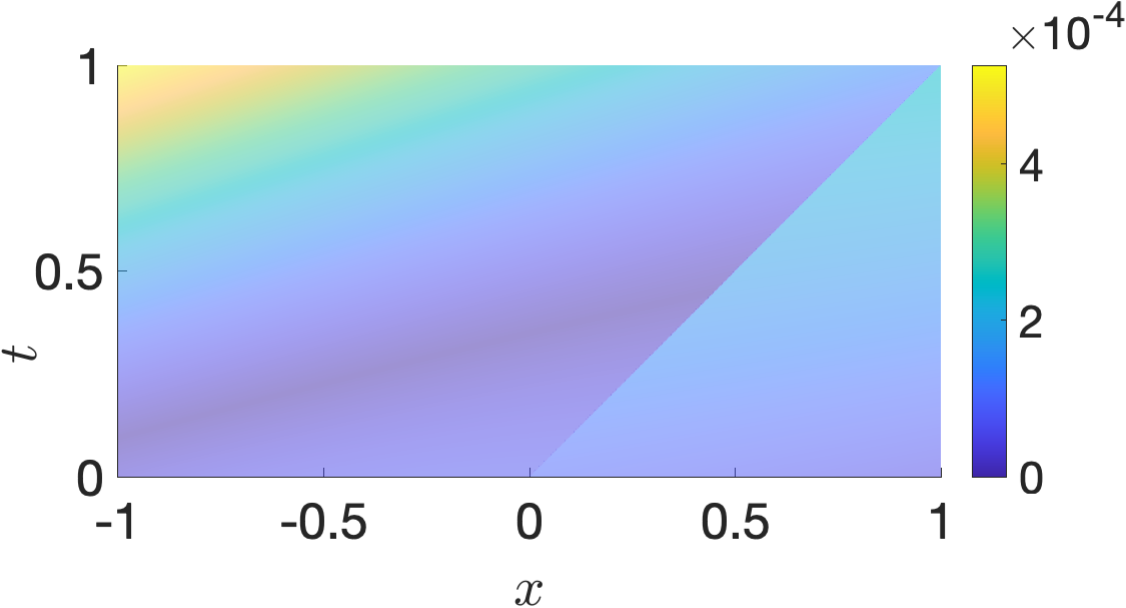}
\caption{$| u(x,t) - \check{u}(x,t) |$}
\end{subfigure}
\vspace{0.3cm}

\centering
\begin{subfigure}{0.32\textwidth}
\centering
\includegraphics[width=\textwidth]{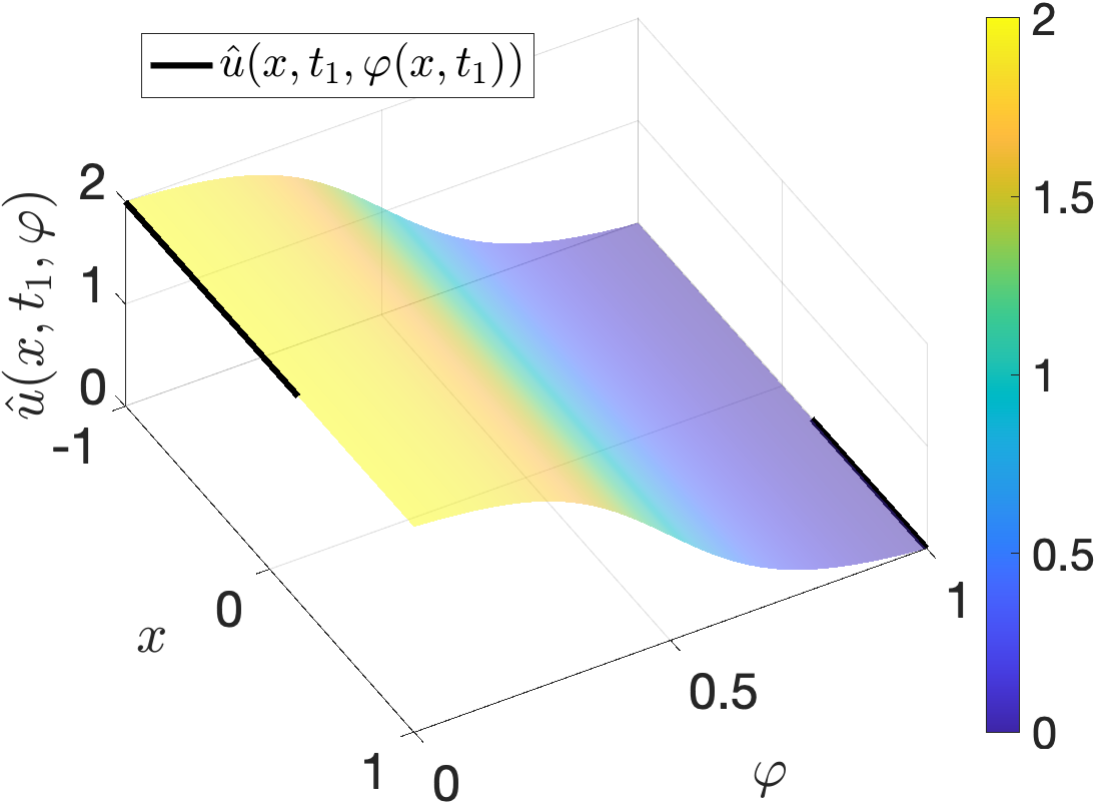}
\caption{$\hat{u}(x,t_1,\varphi)$}
\end{subfigure}
\hspace{0.15cm}
\begin{subfigure}{0.32\textwidth}
\centering
\includegraphics[width=0.9\textwidth]{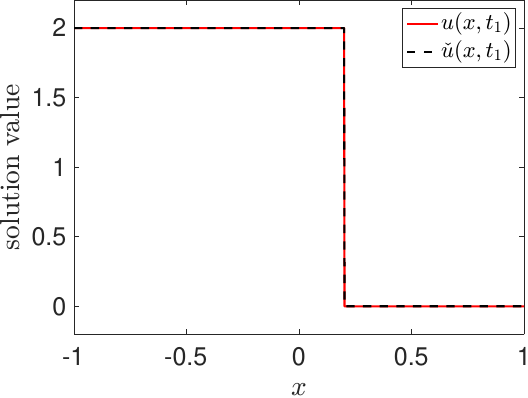}
\caption{$\check{u}(x,t_1)$ and $u(x,t_1)$ }
\end{subfigure}
\hspace{0.15cm}
\begin{subfigure}{0.32\textwidth}
\centering
\includegraphics[width=0.9\textwidth]{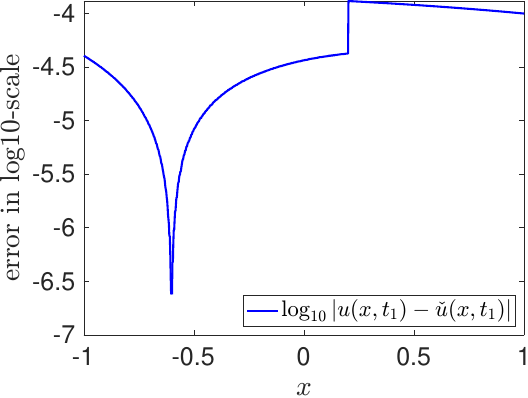}
\caption{$\log_{10}|u(x,t_1) - \check{u}(x,t_1) |$}
\end{subfigure}
\vspace{0.3cm}

\centering
\begin{subfigure}{0.32\textwidth}
\centering
\includegraphics[width=\textwidth]{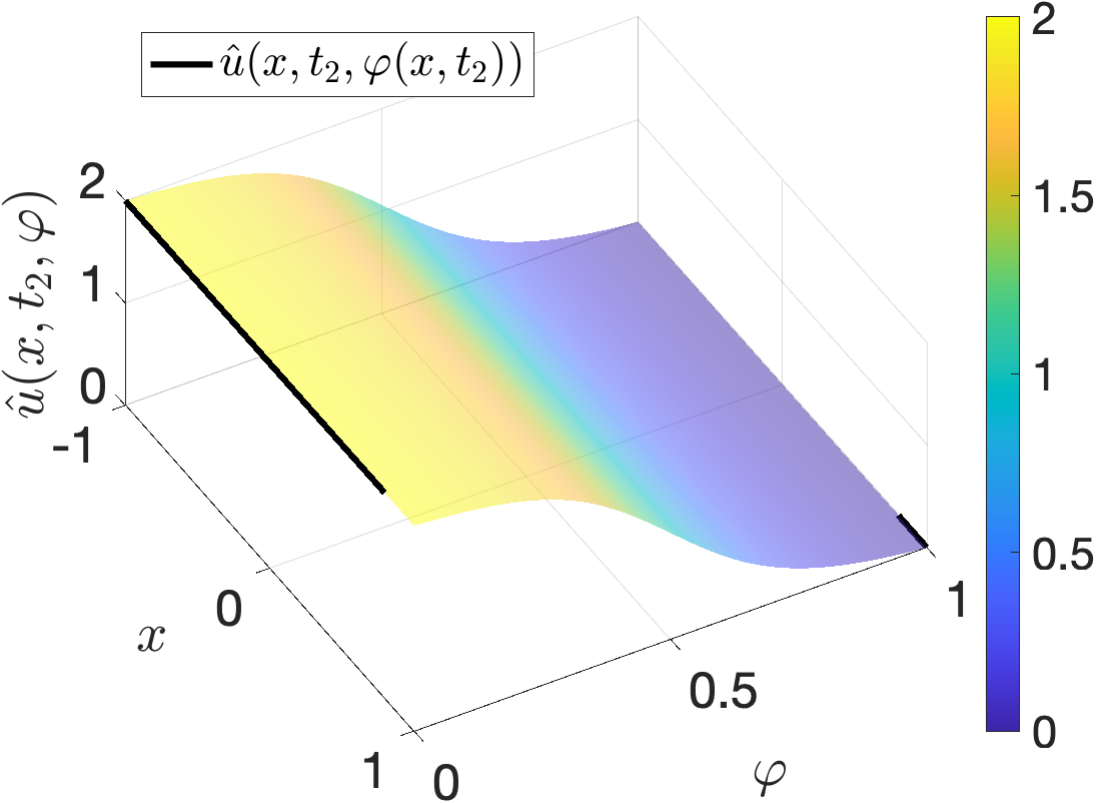}
\caption{$\hat{u}(x,t_2,\varphi)$}
\end{subfigure}
\hspace{0.15cm}
\begin{subfigure}{0.32\textwidth}
\centering
\includegraphics[width=0.9\textwidth]{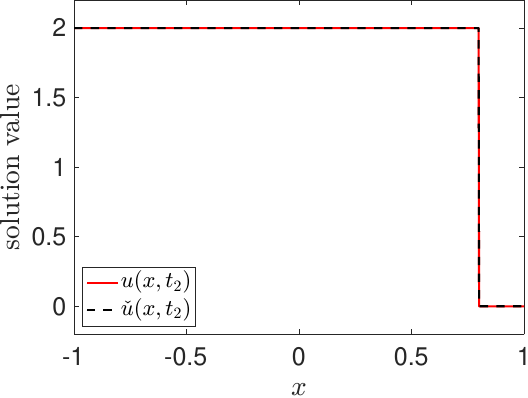}
\caption{$\check{u}(x,t_2)$ and $u(x,t_2)$ }
\end{subfigure}
\hspace{0.15cm}
\begin{subfigure}{0.32\textwidth}
\centering
\includegraphics[width=0.9\textwidth]{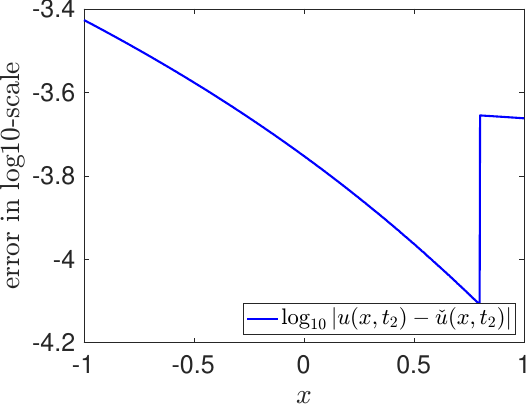}
\caption{$\log_{10}|u(x,t_2) - \check{u}(x,t_2) |$}
\end{subfigure}
\caption{Numerical results for the inviscid Burgers’ equation \eqref{Exp2-Eqns-u-exact} ($t_1=0.2$, $ t_2=0.8$).}
\label{Exp2-fig-1d-Burger-disCsolu}
\vspace{-0.3cm}
\end{figure}
%---------------------------------------%

Next, we consider the following one-dimensional inviscid  Burgers' equation \cite{gustafsson2013time,leveque2002finite}:
\begingroup
\renewcommand*{\arraystretch}{1.2}
\begin{equation}
\left\{
\begin{array}{ll}
\displaystyle \partial_t u(x,t) + u(x,t) \partial_x u(x,t) = 0, \ &\ \ \ \textnormal{for}\ \ (x,t)\in \Omega=(-1, 1)\times (0, 1],\\
\displaystyle u_0(x) = 2 H(-x), \ &\ \ \ \textnormal{for}\ \ x\in (-1, 1), \\
\displaystyle u(-1,t) = 2,\  u(1,t) = 0,\ &\ \ \ \textnormal{for}\ \ t\in (0, 1].
\end{array}\right.
\label{Exp2-Eqns-u-exact}
\end{equation}
\endgroup
In contrast to the linear equations \eqref{Exp1-Eqns-u-exact}, the characteristics of problem \eqref{Exp2-Eqns-u-exact} now cross each other as shown in \autoref{fig-characteristics-inviscid-burgers-1}, leading to the formation of shock discontinuities at intersection points.

Based on the Rankine-Hugoniot relation, the shock speed is a constant, i.e.,
\begin{equation*}
    s = \frac{ \llbracket f(u) \rrbracket }{ \llbracket u \rrbracket } = \frac12 \left( u^+(x,t) + u^-(x,t) \right) = 1\ \ \ \textnormal{where}\ \ f(u) = \frac12 u^2 \ \ \textnormal{and}\ \ \bm{\nu}=1,
\end{equation*}
and therefore the shock curve is a straight line departing from the origin. Then, our solution ansatz is defined as
\begin{equation}
	u(x,t) = \hat{u}(x,t,\varphi(x,t))\ \ \ \textnormal{with} \ \ \ \varphi(x,t) = H(x - st),
	\label{Exp2-uNN-ansatz}
\end{equation}
which also allows for the reformulation of Rankine-Hugoniot jump condition as
\begin{equation*}
	s = \frac{ \llbracket f(\hat{u}) \rrbracket }{ \llbracket \hat{u} \rrbracket } = \frac12 \left( \hat{u}(x,t,\varphi^+(x,t)) + \hat{u}(x,t,\varphi^-(x,t)) \right) \ \ \ \textnormal{on}\ \ \Gamma = \left\{ \left( x(t), t \right) \, \big|\, x(t) = s t \right\}.
\end{equation*}
In other words, the difference in values of $\hat{u}(x,t,\varphi(x,t))$ on opposing sides of the discontinuity interface $\Gamma$ corresponds to the jump of original solution $u(x,t)$ along the shock curve (see \autoref{Exp2-fig-a} for its numerical realization).

On the other hand, a direct computation implies that
\begin{equation*}
\partial_t u = \partial_t \hat{u} - s\partial_\varphi \hat{u} \delta( x-st ) \ \ \ \textnormal{and}\ \ \ \partial_x u = \partial_x \hat{u}+ \partial_\varphi \hat{u} \delta( x-st ),
\end{equation*}
and therefore the hyperbolic system satisfied by our solution ansatz $\hat{u}(x,t,\varphi(x,t))$ is given by 
\begingroup
\renewcommand*{\arraystretch}{1.2}
\begin{equation}
\left\{
\begin{array}{ll}
\partial_t \hat{u}(x,t,\varphi(x,t)) - \hat{u}(x,t,\varphi(x,t)) \partial_x \hat{u}(x,t,\varphi(x,t)) = 0, \ &\ \ \ \textnormal{for}\ \ (x,t)\in \Omega\setminus \Gamma,\\
\hat{u}(x,t,\varphi^+(x,t)) + \hat{u}(x,t,\varphi^-(x,t)) = 2,  \ &\ \ \ \textnormal{for}\ \ (x,t)\in \Gamma,\\
\hat{u}(x,0,\varphi(x,0)) = u_0(x), \ &\ \ \ \textnormal{for}\ \ x\in (-1,1), \\
\hat{u}(-1,t,\varphi(-1,t)) = 2,\ \hat{u}(1, t, \varphi(1, t)) = 0,\ &\ \ \ \textnormal{for}\ \ t\in (0,1],
\end{array}\right.
\label{Eqns-Exp2-u-NN}
\end{equation}
\endgroup
which offers a smooth representation of the shock wave problem \eqref{Exp2-Eqns-u-exact} within a one-order higher-dimensional space. To address the computational challenge posed by the increased dimensionality, a fully-connected neural network \cite{goodfellow2016deep} is deployed to parametrize our solution ansatz, followed by a standard training procedure \cite{raissi2019physics,karniadakis2021physics} as before.

The trained neural network solution and its projection back to the $(x,t)$-plane, denoted as $\hat{u}(x,t,\varphi(x,t))$ and $\check{u}(x,t)$ respectively, are shown in \autoref{Exp2-fig-1d-Burger-disCsolu}, as well as their error profiles. Remarkably, the approximation accuracy for inviscid Burgers’ equation \eqref{Eqns-Exp2-u-NN} is comparable to that of the convection problem \eqref{Eqns-Exp1-u-NN}, which showcases the flexibility of our learning methods in addressing both linear and nonlinear problems. Moreover, the inclusion of augmented variable facilitates the representation of discontinuous solutions on a smooth surface (see (d) and (g) in \autoref{Exp2-fig-1d-Burger-disCsolu}), therefore allowing for a satisfactory reconstruction of discontinuities without laborious adjustment of sampling strategies or penalty viscous terms.
%%%%%%%%%%%%%%%%%%%%%%%%%%%%%%%%%%%%%%%%%%%%%%%%%%%%%
%%%%%%%%%%%%%%%%%%%%%%%%%%%%%%%%%%%%%%%%%%%%%%%%%%%%%

%%%%%%%%%%%%%%%%%%%%%%%%%%%%%%%%%%%%%%%%%%%%%%%%%%%%%
%%%%%%%%%%%%%%%%%%%%%%%%%%%%%%%%%%%%%%%%%%%%%%%%%%%%%
\subsection{The Lift-and-Embed Learning Method}

Thanks to the lift-and-embed operation, a fully connected neural network \cite{goodfellow2016deep} with smooth activation functions, e.g., $\sigma(x)=\tanh(x)$, can be naturally applied to parametrize our solution ansatz (not relabelled), i.e.,
\begin{equation*}
	\hat{u}(x,t,\varphi;\theta) = (L_{D+1} \circ \sigma \circ L_D \circ \cdots \circ \sigma \circ L_1 \circ L_0)(x,t,\varphi)
\end{equation*}
where $D\in\mathbb{N}_+$ denotes the number of hidden layers and $\theta$ indicates the collection of model parameters. It is also noteworthy that the value of input variable $\varphi(x,t)$ can be promptly determined for each sample point $(x,t)$.

%%%%%%%%%%%%%%%%%%%%%%%%%%%%%%%%%%%%%%%%%%%%%%%%%%%%%
\subsubsection{With Identified Discontinuity Locations}

We now introduce our deep learning algorithm under the residual-minimization framework \cite{raissi2019physics}, starting with the generation of collocation points 
\begin{equation*}
	X_{\textnormal{Intrr}} = \big\{(x_n, t_n)\big\}_{n=1}^{N_{\textnormal{Intrr}}},   \quad X_{\textnormal{Shock}} = \big\{(x_n, t_n) \big\}_{n=1}^{N_{\textnormal{Shock}}},  \quad  X_{\textnormal{Bndry}} = \big\{(x_n, t_n)\big\}_{n=1}^{N_{\textnormal{Bndry}}},  \quad  \text{and}  \quad  X_{\textnormal{Initl}} = \big\{(x_n, 0)\big\}_{n=1}^{N_{\textnormal{Initl}}}, 
\end{equation*}
that are sampled uniformly at random from the interior region $\Omega\times(0,T] \setminus \Gamma$, surfaces of discontinuity $\Gamma$, boundaries $\partial\Omega\times(0,T)$ and $\Omega\times\{t=0\}$ of the problem \eqref{GeneralProb-StrongForm}. Here, $N_{\textnormal{Intrr}}$, $N_{\textnormal{Shock}}$, $N_{\textnormal{Bndry}}$, and $N_{\textnormal{Initl}}$ denote the batch sizes of training datasets $X_{\textnormal{Intrr}}$, $X_{\textnormal{Shock}}$, $X_{\textnormal{Bndry}}$, and $X_{\textnormal{Initl}}$, respectively. 

\begin{remark}
Critically, our solution ansatz is defined on piecewise surfaces within a one-order higher-dimensional space (see \autoref{Exp1-fig-a} or \autoref{Exp2-fig-a}), rather than fulfilling the entire $(x,t,\varphi)$-space. Hence, the number of collocation points remains unchanged regardless of the increased dimensionality, which is highly desirable for sample generation.
\end{remark}

Next, by defining empirical loss functions according to the residual of equations \eqref{GeneralProb-Embeded-Form}
\begingroup
\renewcommand*{\arraystretch}{2.5}
\vspace*{-0.2cm}
\begin{equation*}
\begin{array}{c}
\displaystyle L_{\textnormal{Intrr}} (\theta) = \frac{1}{N_{\textnormal{Intrr}}} \sum_{n=1}^{N_{\textnormal{Intrr}}} \big| \partial_t\hat{u}(x_n,t_n, \varphi(x_n, t_n); \theta) - \nabla_x \cdot f(\hat{u}(x_n,t_n, \varphi(x_n, t_n); \theta)) \big|^2,\\
\displaystyle L_{\textnormal{Shock}} (\theta) = \frac{1}{N_{\textnormal{Shock}}} \sum_{n=1}^{N_{\textnormal{Shock}}} \big| s(x_n, t_n) \llbracket \hat{u}(x_n,t_n, \varphi(x_n, t_n); \theta) \rrbracket - \sum_{j=1}^d \nu_j \llbracket f (\hat{u}(x_n,t_n, \varphi(x_n, t_n); \theta)) \rrbracket \big|^2, \\
\displaystyle L_{\textnormal{Bndry}} (\theta) = \frac{1}{N_{\textnormal{Bndry}}} \sum_{n=1}^{N_{\textnormal{Bndry}}} \big|  \hat{u}(x_n,t_n, \varphi(x_n, t_n); \theta) - g(x_n,t_n) \big|^2, \ \ 
\displaystyle L_{\textnormal{Initl}} (\theta) = \frac{1}{N_{\textnormal{Initl}}} \sum_{n=1}^{N_{\textnormal{Initl}}} \big|  \hat{u}(x_n, 0, \varphi(x_n, 0); \theta) - u_0(x_n) \big|^2,
\end{array}
\end{equation*}
\endgroup
the minimization task associated with our deep learning algorithm takes on the form
\begin{equation}
	\theta^* = \operatorname*{arg\, min}_\theta\ L_{\textnormal{Intrr}} (\theta) + \beta_{\textnormal{S}} L_{\textnormal{Shock}} (\theta) + \beta_{\textnormal{B}} L_{\textnormal{Bndry}} (\theta) + \beta_{\textnormal{I}} L_{\textnormal{Initl}} (\theta),
	\label{LE-Loss-Forward}
\end{equation}
where $\beta_{\textnormal{S}}$, $\beta_{\textnormal{B}}$, and $\beta_{\textnormal{I}}>0$ are user-defined penalty parameters. Readers are directed to Algorithm \ref{Algorithm-LPLM} for a detailed illustration of the corresponding lift-and-embed learning algorithm.

%% algorithm %
%%--------------------------------------%
\begin{figure}[t!]
\vspace*{-0.8cm}
\begin{algorithm}[H]
\caption{Lift-and-Embed Learning Algorithm with Identified Discontinuity Locations}
\fontsize{10}{12}\selectfont
\begin{algorithmic}
\STATE{\% \textit{Preparation} }
\STATE{-- construct the augmented variable $\varphi(x, t)$ with identified discontinuity locations;}
%\STATE{\% \textit{Initialization} }
\STATE{-- generate training datasets $X_{\textnormal{Intrr}}$, $X_{\textnormal{Shock}}$, $X_{\textnormal{Bndry}}$, and $X_{\textnormal{Initl}}$, then compute their augmented variables;;}
%\STATE{-- compute the value of augmented variable for each sample;}

\STATE{\% \textit{Training Process}}
\STATE{-- construct and initialize the network model $\hat{u}(x,t,\varphi; \theta)$;}
%\STATE{-- weight initialization for trainable parameters;}
\WHILE{maximum number of epochs is not reached}
\STATE{-- draw mini-batch data uniformly at random from training datasets;}
\STATE{-- network training on the shuffled dataset with a suitable learning rate, i.e., }
\STATE{
\vspace{-.28cm}
\begin{equation*}
	\theta^* = \operatorname*{arg\, min}_\theta \ L_{\textnormal{Intrr}} (\theta) + \beta_{\textnormal{S}} L_{\textnormal{Shock}} (\theta) + \beta_{\textnormal{B}} L_{\textnormal{Bndry}} (\theta) + \beta_{\textnormal{I}} L_{\textnormal{Initl}} (\theta);
\end{equation*}
\vspace{-0.4cm}
}
\ENDWHILE

\STATE{\% \textit{Testing Process} }
\STATE{-- forward pass of the trained model on the testing dataset, i.e.,}
\STATE{
\vspace{-.32cm}
\begin{equation*}
	\check{u}(x,t) = \hat{u}(x,t,\varphi(x, t); \theta^*).
\end{equation*}
\vspace{-0.4cm}
}
\end{algorithmic}
\label{Algorithm-LPLM}
\end{algorithm}
\vspace{-0.6cm}
\end{figure}
%%--------------------------------------%

%%%%%%%%%%%%%%%%%%%%%%%%%%%%%%%%%%%%%%%%%%%%%%%%%%%%%
\subsubsection{With Unknown Discontinuity Locations}

Differing from the previous study, we now explore scenarios in which the location of discontinuity is not known $a$-$priori$. Here, we will focus on the one-dimensional inviscid Burgers' equations, as the discontinuity curves of linear equations could typically be ascertained due to the parallelism of characteristic lines.

Note that in the one-dimensional case, the shock curve has a parametrization of the form $\Gamma = (t, \gamma(t) )$ where $\gamma(0) = x_0$ and $s(t) = \frac{d \gamma(t)}{dt}$ \cite{godlewski2013numerical}. Then, given an approximated shock speed $\hat{s}^{[k]}(t)$ at the $k$-th training epoch, the corresponding shock curve $\hat{\gamma}^{[k]}(t)$ can be determined by solving ordinary differential equations
\begin{equation}
\left\{
\begin{array}{l}
\displaystyle \frac{d \hat{\gamma}^{[k]}(t)}{dt} = \hat{s}^{[k]}(t), \ \ \ \textnormal{for}\ \ t\in (0,T],\\[0.2cm]
\hat{\gamma}^{[k]}(0) = x_0.
\end{array}\right.
\label{ODE-shock-speed-inverse}
\end{equation}
Then, the augmented variable reads
\begin{equation*}
\varphi^{[k]}(x, t) = H(x - \hat{\gamma}^{[k]}(t)),
\end{equation*}
followed by executing Algorithm \ref{Algorithm-LPLM} to obtain the trained neural network solution $\hat{u}^{[k]}(x,t,\varphi^{[k]}(x,t); \theta^*)$. This, in turn, facilitates the employment of Rankine-Hugoniot jump condition \eqref{RH-jump-condition} to adjust the value of shock speed for the next iteration, namely,
\begin{equation*}
	\hat{s}^{[k+1]}(t) =  \operatorname*{arg\, min}_{\hat{s}(t)} \int_\Gamma \big| \hat{s}(t) \llbracket \hat{u}^{[k]}(x,t,\varphi^{[k]}(x,t); \theta^*) \rrbracket - \llbracket f( \hat{u}^{[k]} (x,t,\varphi^{[k]}(x,t); \theta^*)  ) \rrbracket \big|^2 dS.
\end{equation*}
In other words, both neural network model and shock speed are evaluated within an inverse problem framework, in which physics-informed machine learning \cite{raissi2019physics,karniadakis2021physics} has proven its flexibility and effectiveness. 

To numerically realize the inference of shock curve alongside the training of neural network, the problem \eqref{ODE-shock-speed-inverse} is solved with instantaneous values $\{\hat{s}^{[k]}_i = \hat{s}^{[k]}(t_i)\}_{i=1}^n$ on discretized grids $0 = t_1 < \cdots < t_n = T$ to get $\{\hat{\gamma}^{[k]}(t_i)\}_{i=1}^n$. Then, the numerical shock speed $\hat{s}^{[k]}(t)$ and shock curve $\hat{\gamma}^{[k]}(t)$ can be determined by constructing their interpolating polynomials from datasets $\{\hat{s}^{[k]}(t_i)\}_{i=1}^n$ and $\{\hat{\gamma}^{[k]}(t_i)\}_{i=1}^n$ \cite{suli2003introduction}, respectively. Furthermore, the parametrized shock speed $\{\hat{s}^{[k+1]}_i = \hat{s}^{[k+1]}(t_i)\}_{i=1}^n$ can be integrated as leaf nodes within the existing computational graph, enabling simultaneous update with our neural network solution through a single backpropagation process.

More specifically, by including $\{ \hat{s}_i \}_{i=1}^n$ as trainable parameters and defining the loss function (the superscript $[k]$ is omitted for notational simplicity)
\begin{equation}
L_{\textnormal{Shock}}^{\textnormal{inv}} (\theta, \{ \hat{s}_i \}_{i=1}^n) = L_{\textnormal{Shock}} (\theta ) + \frac{1}{n} \sum_{i=1}^{n} \big| \hat{s}_i \llbracket \hat{u}(\gamma(t_i), t_i, \varphi(\gamma(t_i), t_i); \theta) \rrbracket - \llbracket f (\hat{u} (\gamma(t_i), t_i, \varphi(\gamma(t_i), t_i); \theta) ) \rrbracket \big|^2,
\label{Loss-Shock-Update}
\end{equation}
the shock speed is co-updated with our solution ansatz through the backpropagation of neural network, as illustrated in Algorithm \ref{Algorithm-LPLM-inverse}. Notably, the loss function \eqref{Loss-Shock-Update} preserves a similar formulation as that in \eqref{LE-Loss-Forward}, despite the inclusion of extra collocation points for numerical integration.

%% algorithm %
%%--------------------------------------%
\begin{figure}[t!]
\vspace*{-0.8cm}
\begin{algorithm}[H]
\caption{Lift-and-Embed Learning Algorithm with Unknown Discontinuity Locations}
\fontsize{10}{12}\selectfont
\begin{algorithmic}
\STATE{\% \textit{Preparation} }
\STATE{-- generate training datasets $X_{\textnormal{Intrr}}$, $X_{\textnormal{Bndry}}$, $X_{\textnormal{Initl}}$, and the discretization grids $\{ t_i \}_{i=1}^n$;}
\STATE{-- initialize $\{ \hat{s}_i \}_{i=1}^n$ and add them into the optimizer as extra trainable parameters;}

\STATE{\% \textit{Training Process}}
\STATE{-- construct and initialize the network model $\hat{u}(x,t,\varphi; \theta)$;}
\WHILE{maximum number of epochs is not reached}
\STATE{\% \textit{Identification of Shock Speed and Shock Curve} }
\STATE{-- solve \eqref{ODE-shock-speed-inverse} through numerical methods, e.g., explicit Runge-Kutta schemes;}
\STATE{-- reconstruct $\hat{s}(t)$ and $\hat{\gamma}(t)$ from $\{\hat{s}(t_i)\}_{i=1}^n$ and $\{\hat{\gamma}(t_i)\}_{i=1}^n$, e.g. interpolations;}
\STATE{-- update the augmented variable $\varphi(x, t)$ and generate the training dataset $X_{\textnormal{Shock}}$;}
\STATE{\% \textit{Update of Network Solution and Shock Speed Concurrently} }
\STATE{-- draw mini-batch data uniformly at random from training datasets;}
\STATE{-- network training on the shuffled dataset with a suitable learning rate, i.e., }
\STATE{
\vspace{-0.3cm}
\begin{equation*}
\theta^*, \{ \hat{s}_i^* \}_{i=1}^n = \operatorname*{arg\, min}_{\theta,\ \hat{s}_i } L_{\textnormal{Intrr}} (\theta) + \beta_{\textnormal{S}} L_{\textnormal{Shock}}^{\textnormal{inv}} (\theta, \{ \hat{s}_i \}_{i=1}^n ) + \beta_{\textnormal{B}} L_{\textnormal{Bndry}} (\theta) + \beta_{\textnormal{I}} L_{\textnormal{Initl}} (\theta);
\end{equation*}
\vspace{-0.2cm}
}
\ENDWHILE
				
\STATE{\% \textit{Testing Process} }
\STATE{-- forward pass of the trained model on the testing dataset, i.e.,}
\STATE{
\vspace{-.32cm}
\begin{equation*}
\check{u}(x,t) = \hat{u}(x,t,\varphi(x, t); \theta^*).
\end{equation*}
\vspace{-0.4cm}
}
\end{algorithmic}
\label{Algorithm-LPLM-inverse}
\end{algorithm}
\vspace{-0.6cm}
\end{figure}
%%--------------------------------------%

%%%%%%%%%%%%%%%%%%%%%%%%%%%%%%%%%%%%%%%%%%%%%%%%%%%%%
%%%%%%%%%%%%%%%%%%%%%%%%%%%%%%%%%%%%%%%%%%%%%%%%%%%%%

%%%%%%%%%%%%%%%%%%%%%%%%%%%%%%%%%%%%%%%%%%%%%%%%%%%%%
%%%%%%%%%%%%%%%%%%%%%%%%%%%%%%%%%%%%%%%%%%%%%%%%%%%%%
\section{Numerical Experiments}\label{Section-Experiment}

%--------------------------------------%
\begin{table}[t!]
	\renewcommand{\arraystretch}{1.4}
	\fontsize{9.5}{12}\selectfont
	\centering
	\begin{tabular}{ccccccc}
		\hline 
		\multicolumn{2}{c}{} & \makecell{Neural Network\\ (Depth, Width)} & \makecell{Penalty Coeff.\\ ($\beta_{\textnormal{S}}$, $\beta_{\textnormal{B}}$, $\beta_{\textnormal{I}}$)} & \makecell{Train Data Size\\ ($N_{\textnormal{Intrr}}$, $N_{\textnormal{Shock}}$, $N_{\textnormal{Bndry}}$, $N_{\textnormal{Initl}}$)} &  \makecell{Test Data Size \\ ($N_x$, $N_t$)} \\ \hline
		\multirow{3}*{\parbox{1.8cm}{\centering Linear Convection Equations}} & Sec. \ref{sec-Exp1-1d-Linear-disCsolu-a1} & (2, 40) & (400, 1, 400) & (10$k$, 4$k$, 2$k$, 1$k$) & (1$k$, 1$k$) \\ \cline{2-6}
		& Sec. \ref{sec-Exp3-1d-Linear-disCsolu-a1} & (2, 40) & (400, 1, 400) & (10$k$, 4$k$, 2$k$, 1$k$) & (1$k$, 1$k$) \\ \cline{2-6}
		& Sec. \ref{sec-Exp5-1d-Linear-disCsolu-a50} & (6, 40) & (400, 10, 400) & (80$k$, 60$k$, 5$k$, 5$k$) & (1$k$, 1$k$)\\ \hline
		\multirow{6}*{\parbox{1.8cm}{\centering Inviscid Burgers' Equations}} & Sec. \ref{sec-Exp2-1d-Burgers} & (3, 40) & (400, 1, 400) & (10$k$, 1$k$, 1$k$, 1$k$) & (1$k$, 1$k$) \\ \cline{2-6}
		& Sec. \ref{sec-Exp6-1d-Burgers-Merging-Shocks} & (6, 40) & (400, 1, 400) & (80$k$, 15$k$, 5$k$, 5$k$) & (1$k$, 1$k$) \\ \cline{2-6} 
		& Sec. \ref{sec-ExpSM1-1d-Burgers-Rarefaction} & (6, 40) & (400, 1, 400) & (80$k$, 10$k$, 10$k$, 5$k$) & (1$k$, 1$k$) \\ \cline{2-6}
		& Sec. \ref{sec-Exp7-2d-Burgers} & (4, 80) & (50, 1, 400) & (80$k$, 20$k$, 30$k$, 10$k$) & ($1k \times 1k$, 17) \\ \cline{2-6}      
		& Sec. \ref{sec-Exp8-1d-iBurgers} & (4, 40) & (400, 1, 400) & (10$k$, 1$k$, 1$k$, 1$k$) & (1$k$, 1$k$) \\ \cline{2-6}
		& Sec. \ref{sec-ExpSM2-1d-iBurgers-Curved} & (6, 40) & (50, 1, 400) & (80$k$, 5$k$, 5$k$, 5$k$) & (1$k$, 1$k$) \\  \hline    
	\end{tabular}	
	\vspace*{-0.15cm}
	\caption{List of key hyper-parameter configurations for numerical experiments.}
	\label{table-hyper-parameters}
%\end{table}
%%--------------------------------------%
\vspace{0.5cm}
\renewcommand{\arraystretch}{1.8}
\begin{tabular}{ccccccc}
	\hline 
	\multirow{2}*{\parbox{1.8cm}{\centering Linear Convection Equations}} & Sec. \ref{sec-Exp1-1d-Linear-disCsolu-a1} & Sec. \ref{sec-Exp3-1d-Linear-disCsolu-a1} & Sec. \ref{sec-Exp5-1d-Linear-disCsolu-a50} & & & \\ \cline{2-4}
	& $2.01\times 10^{-3}$ & $2.49\times 10^{-5}$ & $7.03\times 10^{-4}$ & & & \\ \hline
	\multirow{2}*{\parbox{1.8cm}{\centering Inviscid Burgers' Equations}} & Sec. \ref{sec-Exp2-1d-Burgers} & Sec. \ref{sec-Exp6-1d-Burgers-Merging-Shocks} & Sec. \ref{sec-ExpSM1-1d-Burgers-Rarefaction} & Sec. \ref{sec-Exp7-2d-Burgers} & Sec. \ref{sec-Exp8-1d-iBurgers} & Sec. \ref{sec-ExpSM2-1d-iBurgers-Curved} \\ \cline{2-7}
	& $9.47\times 10^{-5}$ & $1.85\times10^{-4}$ & $4.85 \times 10^{-3}$ & $1.38 \times 10^{-5}$  & $1.82 \times 10^{-2}$ & $3.64 \times 10^{-2}$ \\ \hline		       
\end{tabular}
\caption{Relative $L^2$ errors for numerical experiments conducted in this study.}
\vspace*{-0.15cm}  
\label{table-relative-errors}
\vspace{-0.3cm}
\end{table}
%--------------------------------------%

In this section, we conduct numerical studies on a series of benchmark problems to demonstrate the effectiveness and flexibility of our methods. More precisely, the inclusion of an augmented variable allows for a representation of jump conditions as differences between distinct collocation points, while the parametrization of solution ansatz with neural networks can mitigate the curse of dimensionality and handle shock waves that may propagate along lines, curves, or a combination of both. Our approach is also capable of addressing both linear and quasi-linear equations, including cases with unknown discontinuity interfaces, within a unified learning framework. Upon the completion of the training process, the approximate solution can be obtained by projecting the trained neural network solution back onto the original physical domain.
 
Note that our method offers the merit of using a smooth neural network to learn discontinuous solutions, hence the fully-connected neural network is deployed with the use of hyperbolic tangent activation functions \cite{goodfellow2016deep}. For the following numerical experiments, different hyper-parameter configurations are implemented and summarized in \autoref{table-hyper-parameters}. When the maximum number of training epochs is reached, the model with the smallest training loss, indicated by $\hat{u}(x,t,\varphi(x,t))$ as before, is adopted for comparison with the exact solution $u(x,t)$, while the projected solution onto the original $(x,t)$-plane is indicated as $\check{u}(x,t)$. All network models are trained using the ADAMW optimizer \cite{kingma2014adam,loshchilov2017decoupled}, with an initial learning rate set to $0.01$ and a step decay schedule that divides it by 10 after specific milestones.

To assess the performance of our trained neural network solutions, we generate testing points $\big\{ (x_i, t_j) \big\}_{i=1,\, j=1}^{N_x,\, N_t}$ through a uniform partition of the computational domain, and then report the relative $L^2$ error in \autoref{table-relative-errors}, i.e.,
\begin{equation*}
\frac{ \lVert u(x,t) - \check{u}(x,t) \rVert_2 }{ \lVert u(x,t) \rVert_2 } 
\ \ \ \textnormal{where}\ \ \ 
\lVert u \rVert_2 =  \sqrt{ \sum\nolimits_{m=1}^M \left| u(x_m,t_m) \right|^2}.
\end{equation*}
All experiments are conducted using PyTorch on Nvidia GeForce RTX 4090 GPU cards \cite{paszke2017automatic}. 
%All source codes are publicly accessible on GitHub (\url{https://github.com/q1sun/LELM_HyperbolicPDE}).

%%%%%%%%%%%%%%%%%%%%%%%%%%%%%%%%%%%%%%%%%%%%%%%%%%%%%
\subsection{Linear Convection Equations}

In this section, we report numerical studies on linear convection equations with a discontinuous initial condition, as well as a large convection coefficient \cite{krishnapriyan2021characterizing,daw2022mitigating}, to validate the effectiveness of our methods. Algorithm \ref{Algorithm-LPLM} is used to train the neural network solution.

%%%%%%%%%%%%%%%%%%%%%%%%%%%%%%%%%%%%%%%%%%%%%%%%%%%%%
\subsubsection{The Benchmark Linear Convection Problem Revisited}\label{sec-Exp3-1d-Linear-disCsolu-a1}

\begin{figure}[!htbp]
\centering
\begin{subfigure}{0.32\textwidth}
\centering
\includegraphics[width=\textwidth]{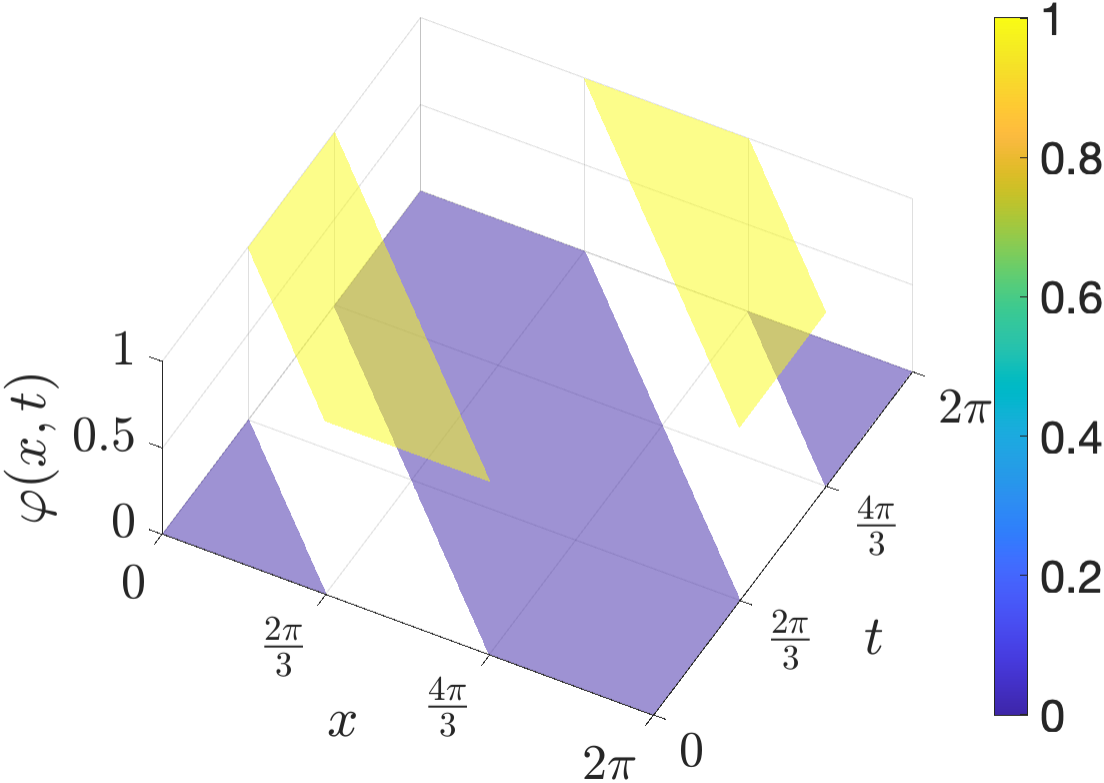}
\caption{$\hat{u}(x,t,\varphi(x,t))$}
\label{Exp3-fig-a}
\end{subfigure}
\hspace{0.15cm}
\begin{subfigure}{0.32\textwidth}
\centering
\includegraphics[width=0.9\textwidth]{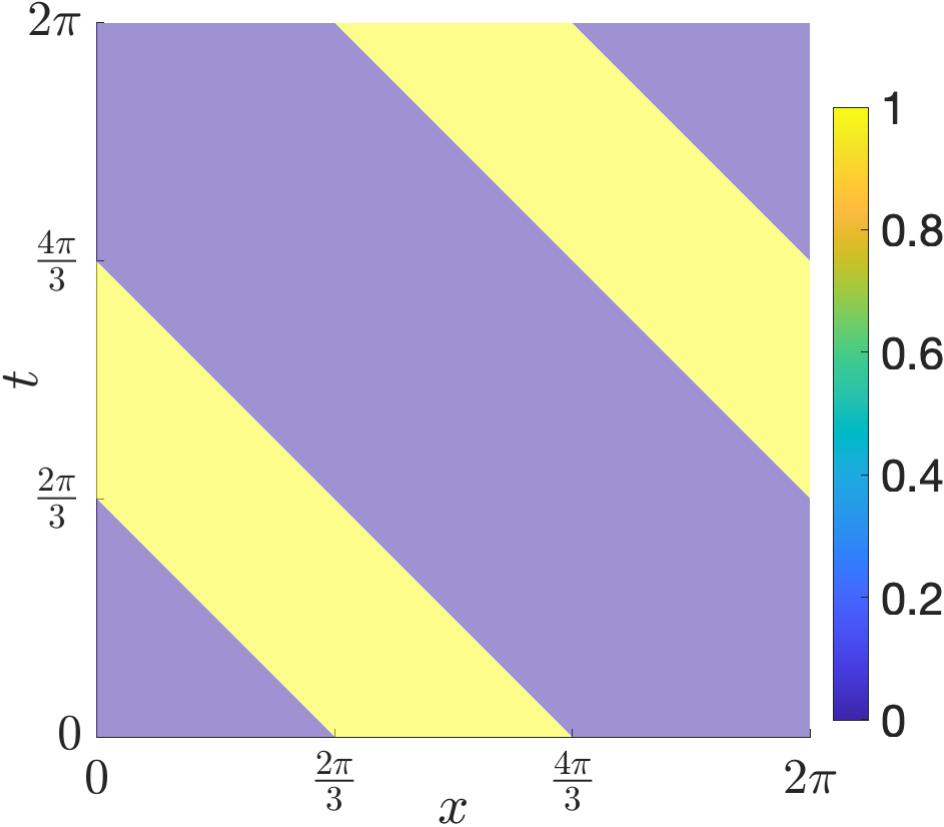}
\caption{$\check{u}(x,t)$}
\end{subfigure}
\hspace{0.15cm}
\begin{subfigure}{0.32\textwidth}
\centering
\includegraphics[width=0.95\textwidth]{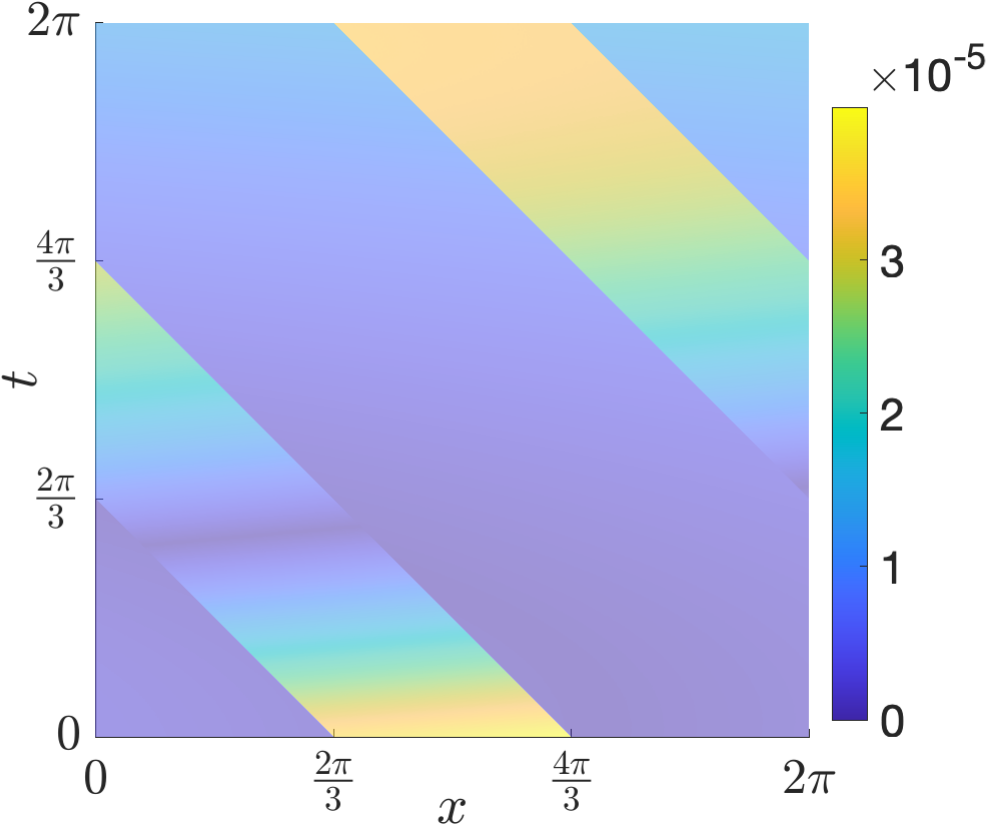}
\caption{$|u(x,t) - \check{u}(x,t) |$}
\end{subfigure}
\vspace{0.3cm}

\centering
\begin{subfigure}{0.32\textwidth}
\centering
\includegraphics[width=\textwidth]{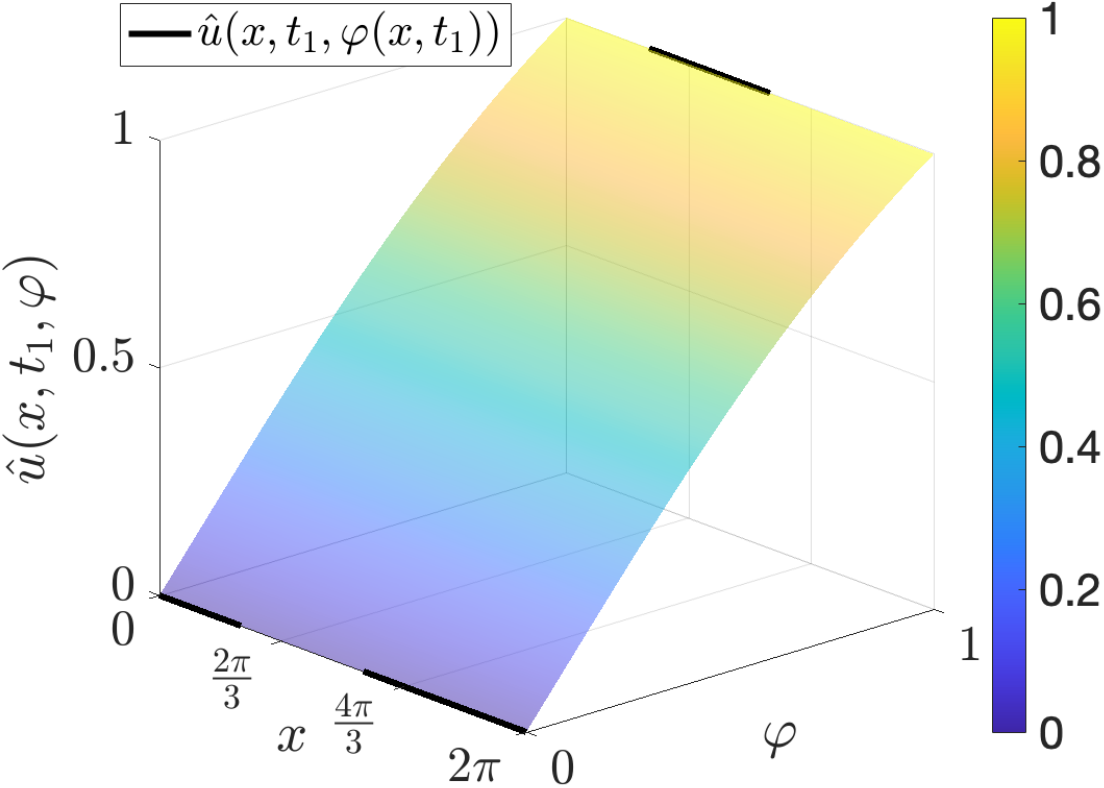}
\caption{$\hat{u}(x,t_1)$}
\end{subfigure}
\hspace{0.15cm}
\begin{subfigure}{0.32\textwidth}
\centering
\includegraphics[width=0.9\textwidth]{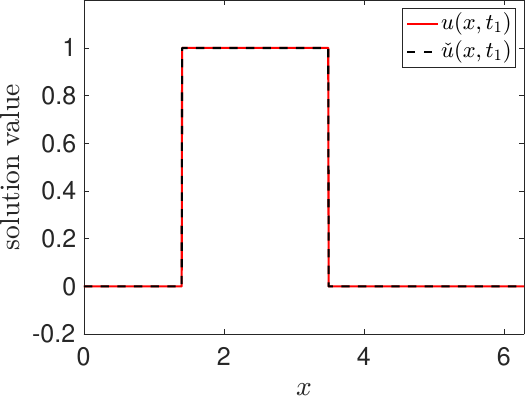}
\caption{$\check{u}(x,t_1)$ and $u(x, t_1)$ }
\end{subfigure}
\hspace{0.15cm}
\begin{subfigure}{0.32\textwidth}
\centering
\includegraphics[width=0.9\textwidth]{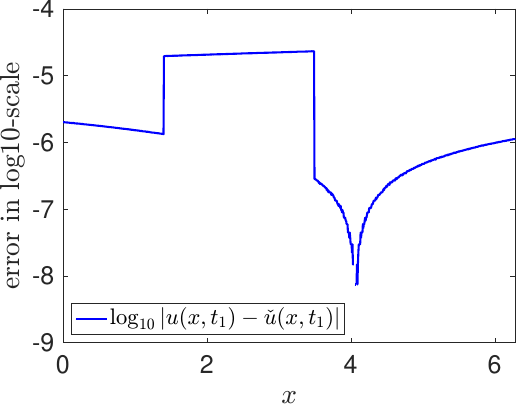}
\caption{$\log_{10}|u(x,t_1) - \check{u}(x,t_1) |$}
\end{subfigure}
\vspace{0.3cm}

\centering
\begin{subfigure}{0.32\textwidth}
\centering
\includegraphics[width=\textwidth]{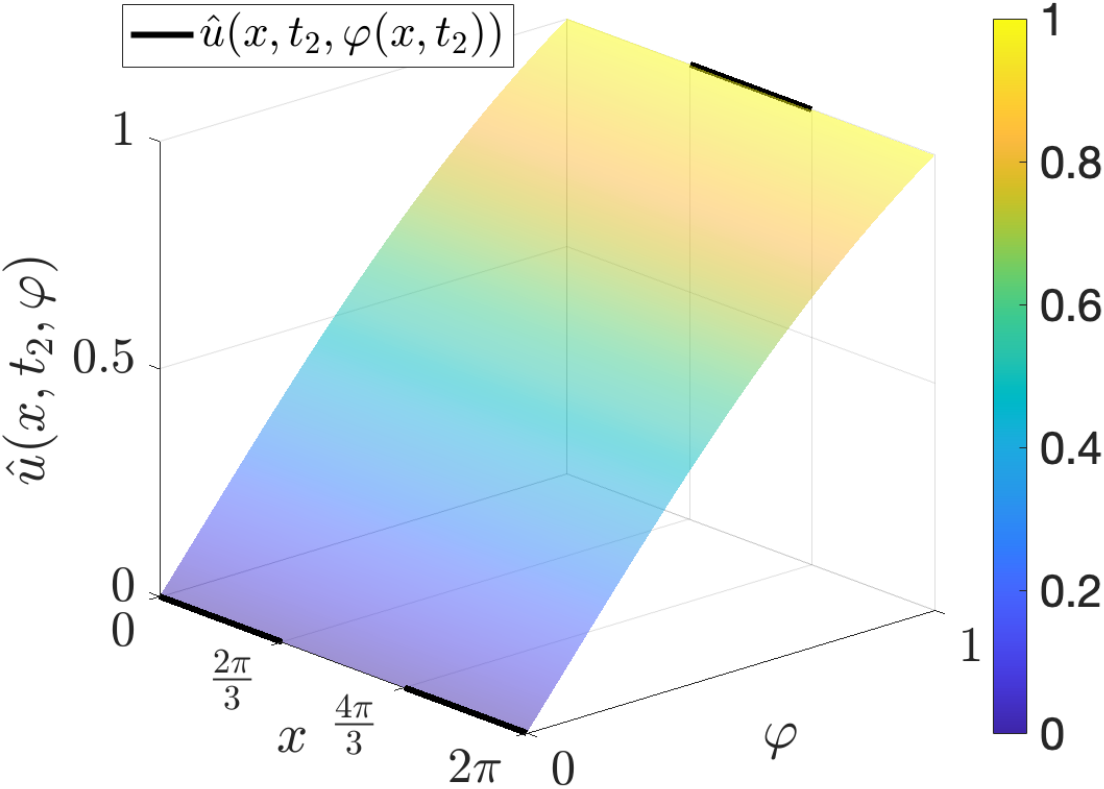}
\caption{$\check{u}(x,t_2)$}
\end{subfigure}
\hspace{0.15cm}
\begin{subfigure}{0.32\textwidth}
\centering
\includegraphics[width=0.9\textwidth]{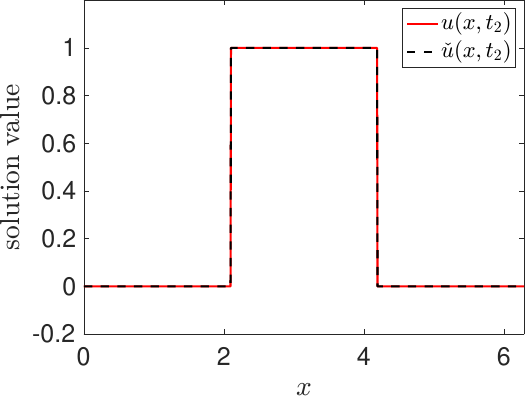}
\caption{$\check{u}(x,t_2)$ and $u(x, t_2)$}
\end{subfigure}
\hspace{0.15cm}
\begin{subfigure}{0.32\textwidth}
\centering
\includegraphics[width=0.9\textwidth]{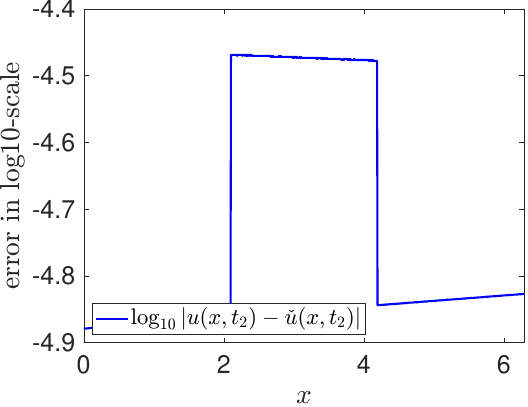}
\caption{$\log_{10}|u(x,t_2) - \check{u}(x,t_2)|$}
\end{subfigure}
\caption{Numerical results for the convection equation \eqref{Exp1-Eqns-u-exact} using another augmented variable ($t_1 = \frac{2\pi}{9}$, $t_2= 2\pi$).}
\label{Exp3-fig-1d-Linear-disCsolu-a1}
\hspace{-0.3cm}
\end{figure}

By revisiting the benchmark linear convection problem \eqref{Exp1-Eqns-u-exact}, we showcase another way of building the augmented variable, as the additional degree of freedom offers non-unique mappings that yield the identical behaviour upon collapsing to the original lower-dimensional domain. To be precise, we consider
\vspace{-0.15cm}
\begin{equation*}
	u(x,t) = \hat{u}(x,t,\varphi(x,t))\ \ \ \textnormal{with} \ \ \ \varphi(x,t) =  \sum_{k=0}^{n} u_0(x-st - 2k\pi),
	\vspace{-0.15cm}
\end{equation*}
where $s=-1$, $n = \ceil*{ - \frac{sT+x_1}{2\pi} }$ with $x_1=\frac{2\pi}{3}$ and $T=2\pi$. Here, the augmented variable $\varphi(x,t)$ is constructed using the initial condition \eqref{Exp1-initial-condition}, which also enables the representation of jump discontinuities as differences across distinct collocation points (see \autoref{Exp3-fig-a}). The lifted hyperbolic system adopts the same formulation as that of \eqref{Eqns-Exp1-u-NN} and the corresponding computational results are reported in \autoref{Exp3-fig-1d-Linear-disCsolu-a1} and \autoref{table-relative-errors}.

As can be seen from (b), (e), and (h) in \autoref{Exp3-fig-1d-Linear-disCsolu-a1}, our trained neural network solution can effectively capture the discontinuous solution without notable smearing or oscillation, even as the time evolves to $t=2\pi$. For comparison, readers are referred to textbook \cite{gustafsson2013time} for the numerical solution obtained using classical finite difference methods. Furthermore, superior accuracy is achieved compared to that in \autoref{Exp1-fig-1d-Linear-disCsolu-a1}, as the augmented variable, or equivalently, the input data of our neural network, coincides with the exact solution for this specific setting.

%---------------------------------------%
%---------------------------------------%

%%%%%%%%%%%%%%%%%%%%%%%%%%%%%%%%%%%%%%%%%%%%%%%%%%%%%
\subsubsection{Discontinuous Solution with A Large Convection Coefficient}\label{sec-Exp5-1d-Linear-disCsolu-a50}

Recent empirical studies suggest that standard learning methods may struggle to resolve the exact solution when the convection coefficient becomes large \cite{cuomo2022scientific}, therefore requiring complex adjustments like sequence-to-sequence learning \cite{krishnapriyan2021characterizing} or retain-resample-release sampling \cite{daw2022mitigating}. To wrap up our study on linear problems, we explore a scenario marked by not only a large convection coefficient but also a discontinuous initial condition, that is,
\begin{equation}
\left\{
\begin{array}{ll}
\partial_t u(x,t) - 50 \partial_x u(x,t) = 0, \ &\ \ \ \textnormal{for}\ \ (x,t)\in \Omega = (0,2\pi)\times (0,\frac{\pi}{5}],\\
u_0(x) = H(x-\frac{2\pi}{3}) - H(x-\frac{4\pi}{3}), \ &\ \ \ \textnormal{for}\ \ x\in (0,2\pi), \\
u(0,t) = u(2\pi, t),\ &\ \ \ \textnormal{for}\ \ t\in (0,\frac{\pi}{5}].
\end{array}\right.
\label{Exp5-Eqns-u-exact}
\end{equation}
Furthermore, the augmented variable for \eqref{Exp5-Eqns-u-exact} is formulated as a step function to avoid potential coincidence with the exact solution. To be precise, we employ the identical configuration for solution representation as outlined in \eqref{Exp1-solution-ansatz}, with the only distinction being the presence of a large convection coefficient within the current problem.

\begin{figure}[!htbp]
\centering
\begin{subfigure}{0.32\textwidth}
\centering
\includegraphics[width=\textwidth]{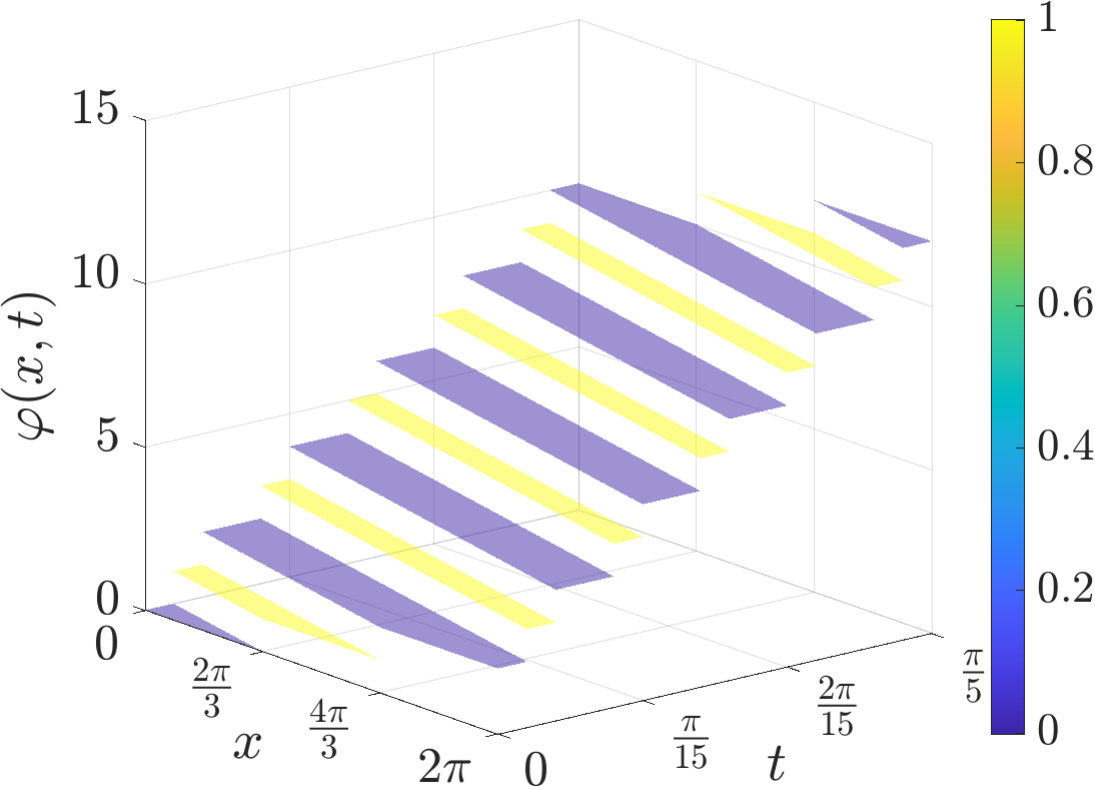}
\caption{$\hat{u}(x,t,\varphi(x, t))$}
\end{subfigure}
\hspace{0.15cm}
\begin{subfigure}{0.32\textwidth}
\centering
\includegraphics[width=0.835\textwidth]{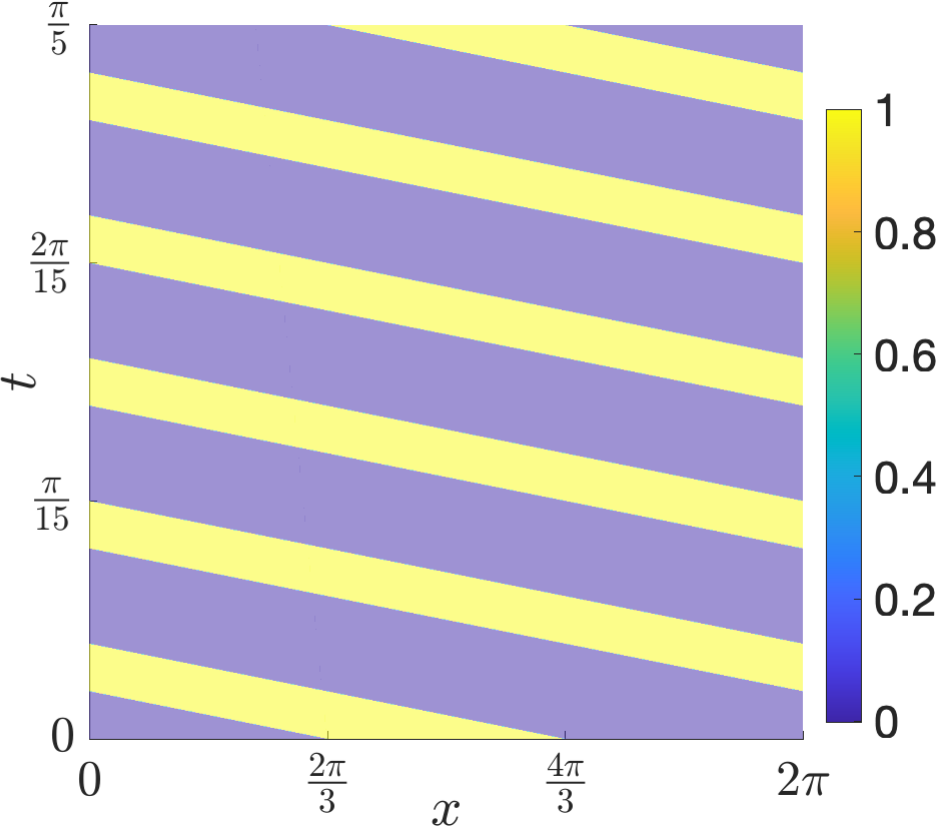}
\caption{$\check{u}(x,t)$}
\end{subfigure}
\hspace{0.15cm}
\begin{subfigure}{0.32\textwidth}
\centering
\includegraphics[width=0.89\textwidth]{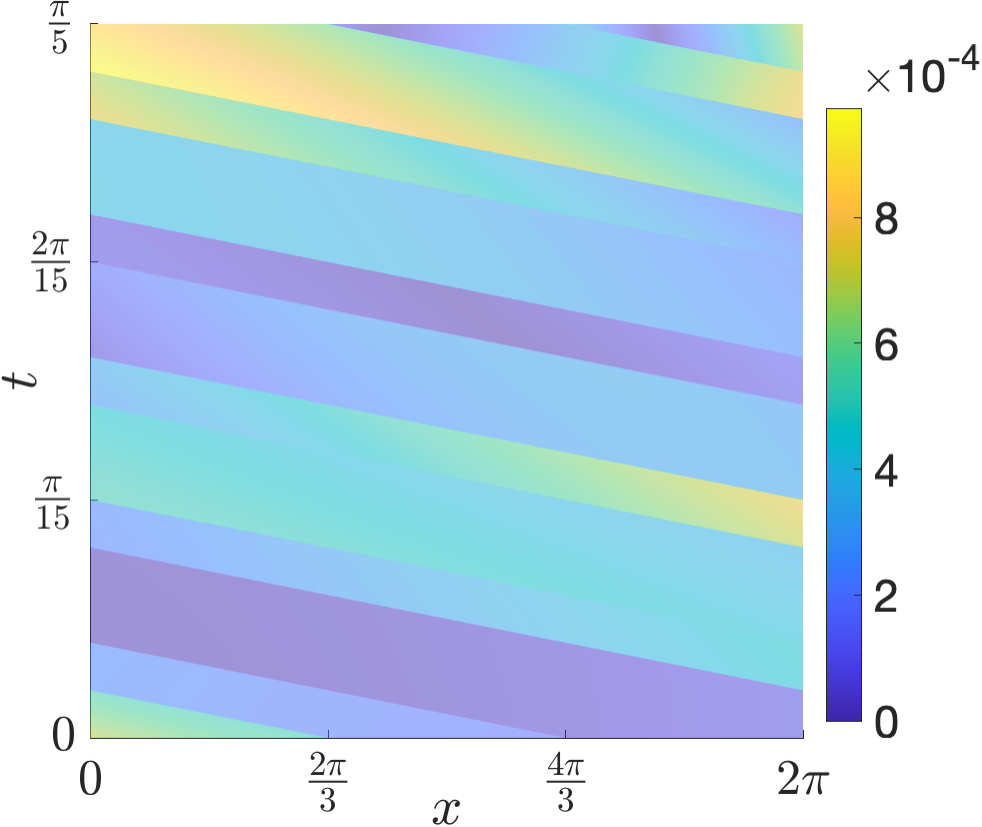}
\caption{$|\check{u}(x,t) - u(x,t)|$}
\end{subfigure}
\vspace{0.3cm}

\centering
\begin{subfigure}{0.32\textwidth}
\centering
\includegraphics[width=\textwidth]{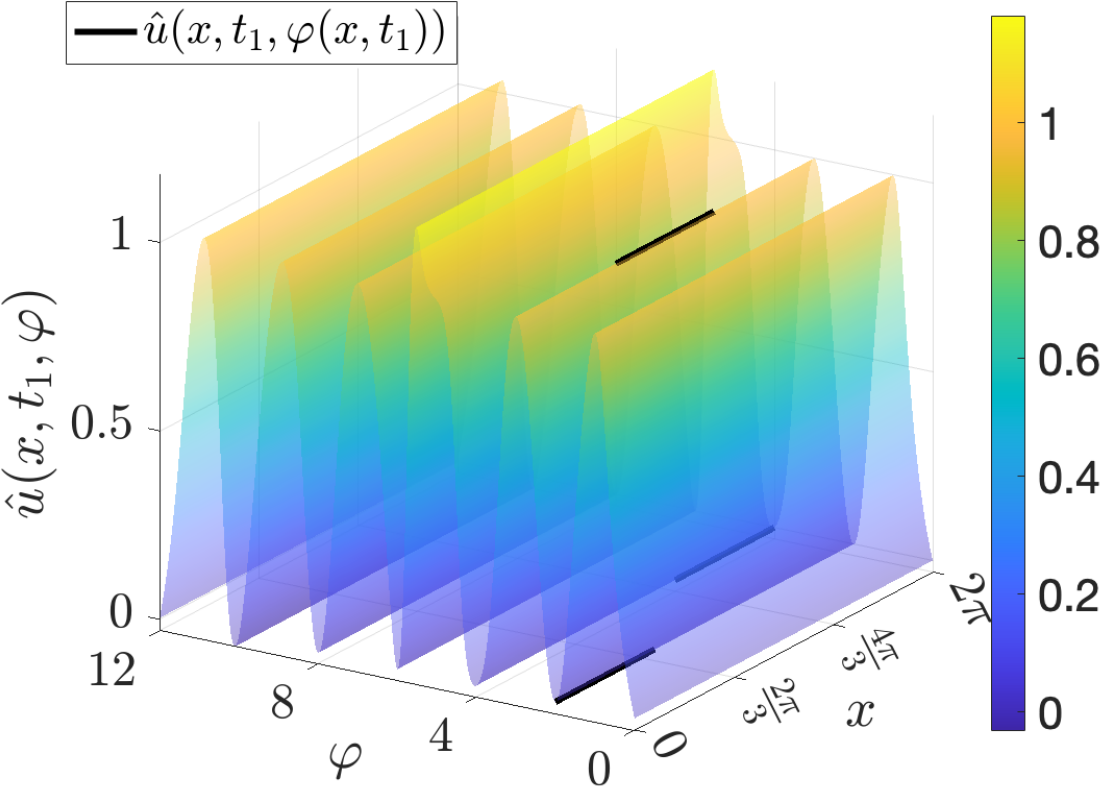}
\caption{$\hat{u}(x,t_1,\varphi)$}
\end{subfigure}
\hspace{0.15cm}
\begin{subfigure}{0.32\textwidth}
\centering
\includegraphics[width=0.9\textwidth]{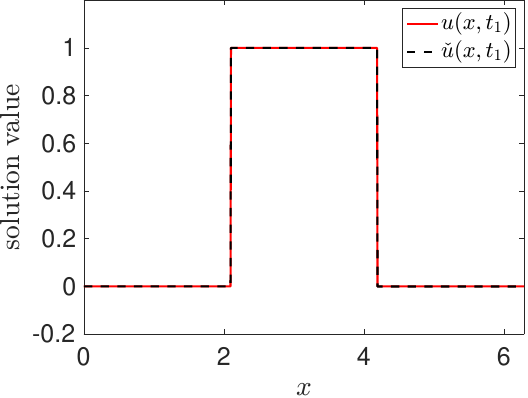}
\caption{$\check{u}(x,t_1)$ and $u(x,t_1)$}
\label{Exp5-fig-e}
\end{subfigure}
\hspace{0.15cm}
\begin{subfigure}{0.32\textwidth}
\centering
\includegraphics[width=0.9\textwidth]{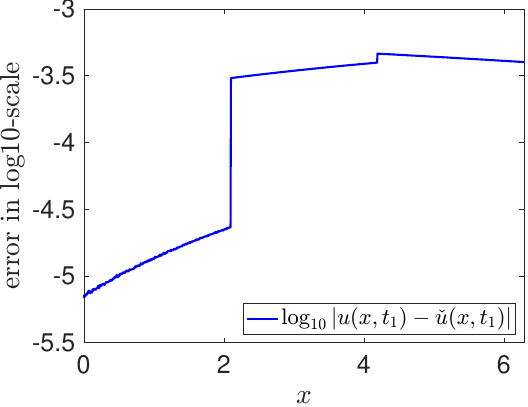}
\caption{$\log_{10}|\check{u}(x,t_1) - u(x,t_1)|$}
\end{subfigure}
\vspace{0.3cm}

\centering
\begin{subfigure}{0.32\textwidth}
\centering
\includegraphics[width=\textwidth]{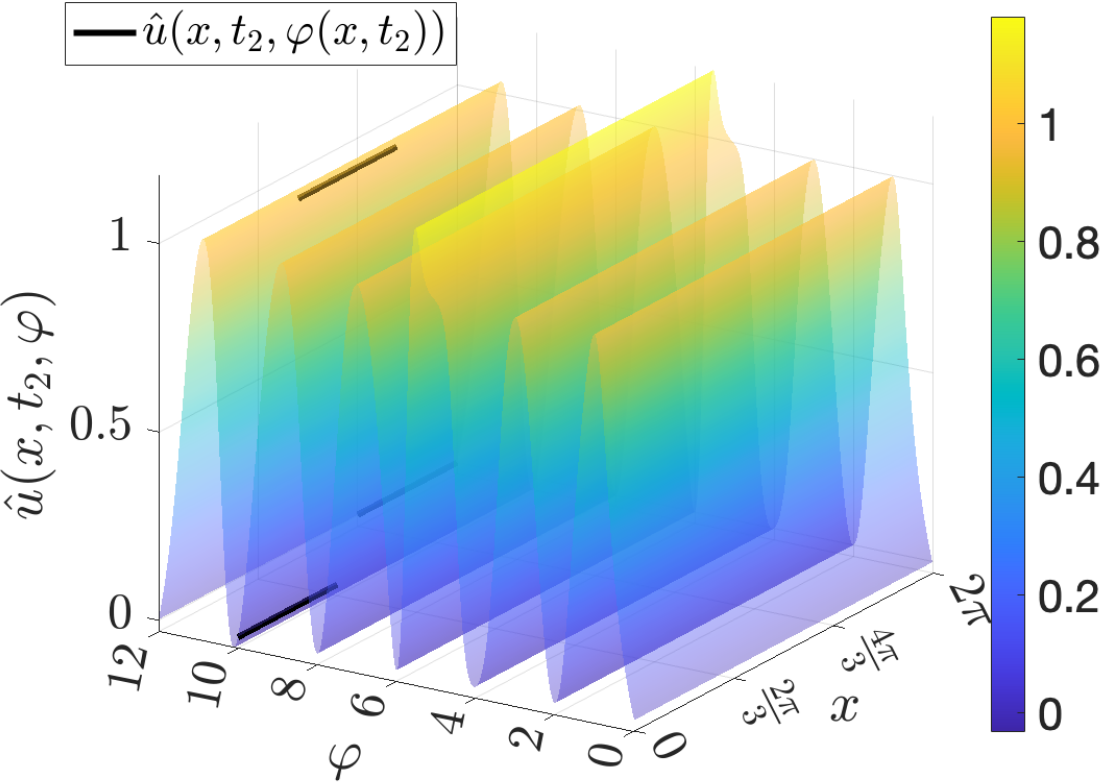}
\caption{$\hat{u}(x,t_2,\varphi)$}
\end{subfigure}
\hspace{0.15cm}
\begin{subfigure}{0.32\textwidth}
\centering
\includegraphics[width=0.9\textwidth]{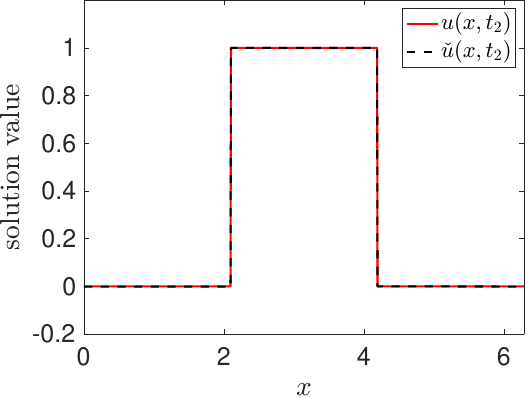}
\caption{$\check{u}(x,t_2)$ and $u(x,t_2)$ }
\label{Exp5-fig-h}
\end{subfigure}
\hspace{0.15cm}
\begin{subfigure}{0.32\textwidth}
\centering
\includegraphics[width=0.9\textwidth]{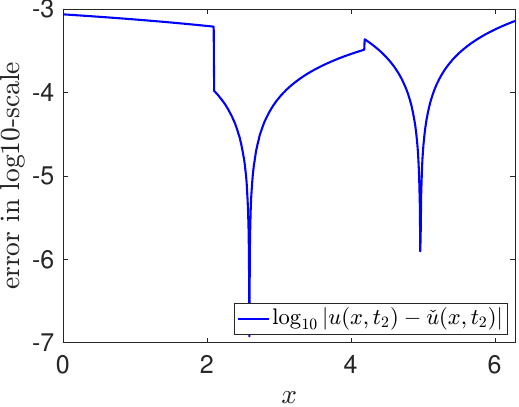}
\caption{$\log_{10}|\check{u}(x,t_2) - u(x,t_2)|$}
\end{subfigure}
\caption{Numerical results for the convection equation \eqref{Exp5-Eqns-u-exact} with a large convection coefficient ($t_1=\frac{\pi}{25}$, $t_2=\frac{\pi}{5}$).}
\label{Exp5-fig-1d-Linear-disCsolu-a50}
\vspace{-0.3cm}
\end{figure}

The neural network solution and its projection back to the lower-dimensional flat plane are shown in \autoref{Exp5-fig-1d-Linear-disCsolu-a50}, as well as error profiles over the entire domain or at specific temporal slices (see (c), (f), and (i) in \autoref{Exp5-fig-1d-Linear-disCsolu-a50}). Clearly, our numerical solution closely matches the exact solution without spurious smearing and oscillation (see \autoref{Exp5-fig-e} and \autoref{Exp5-fig-h}), regardless of the presence of numerous discontinuity interfaces incurred by the large convection coefficient. 

%---------------------------------------%
%%---------------------------------------%

%%%%%%%%%%%%%%%%%%%%%%%%%%%%%%%%%%%%%%%%%%%%%%%%%%%%%
\subsection{Inviscid Burgers' Equations with Identified Shock Locations}

Next, we present numerical studies on the one-dimensional inviscid Burgers' equation, focusing on shock-shock and rarefaction-shock interactions, to show the effectiveness of our methods. Following this, a two-dimensional problem is examined in \autoref{sec-Exp7-2d-Burgers} to showcase the neural network's capability in tackling high-dimensional functions. Here, Algorithm \ref{Algorithm-LPLM} is again used to train the neural network solution.

%%%%%%%%%%%%%%%%%%%%%%%%%%%%%%%%%%%%%%%%%%%%%%%%%%%%%
\subsubsection{Shock-Shock Interaction}\label{sec-Exp6-1d-Burgers-Merging-Shocks}

%---------------------------------------%
\begin{figure}[t!]
\centering
\begin{subfigure}{0.32\textwidth}
\centering
\includegraphics[width=\textwidth]{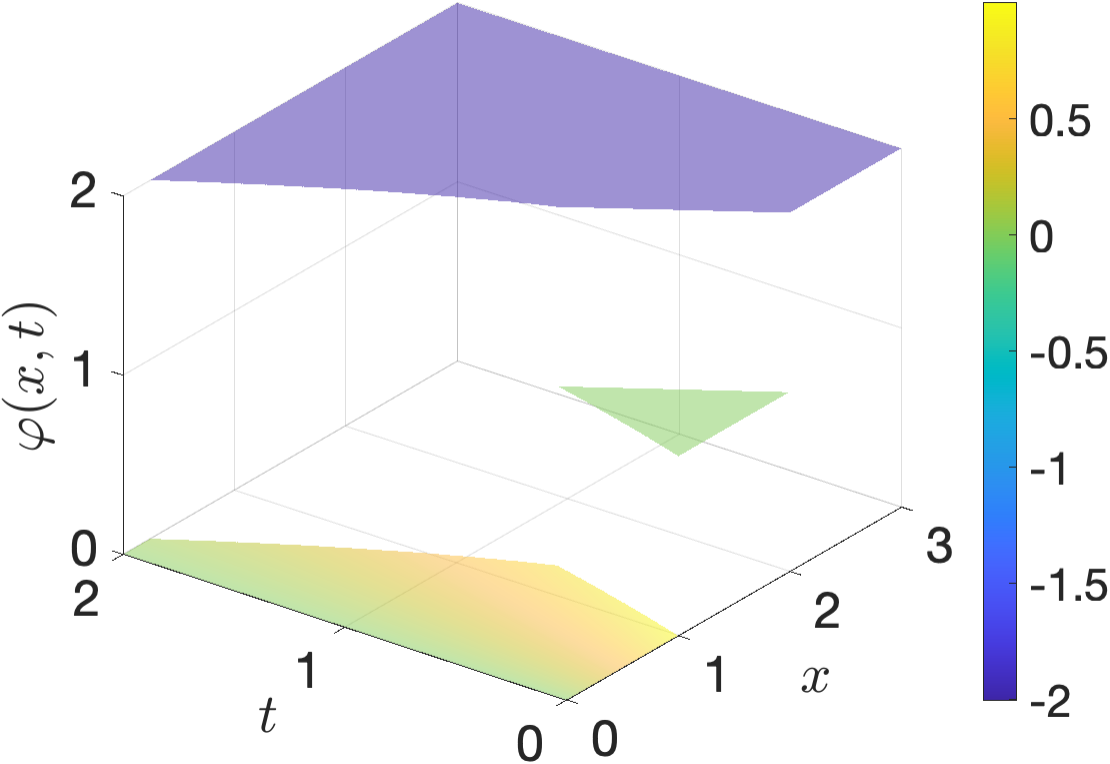}
\caption{$\hat{u}(x,t,\varphi(x,t))$}
\end{subfigure}
\hspace{0.15cm}
\begin{subfigure}{0.32\textwidth}
\centering
\includegraphics[width=0.85\textwidth]{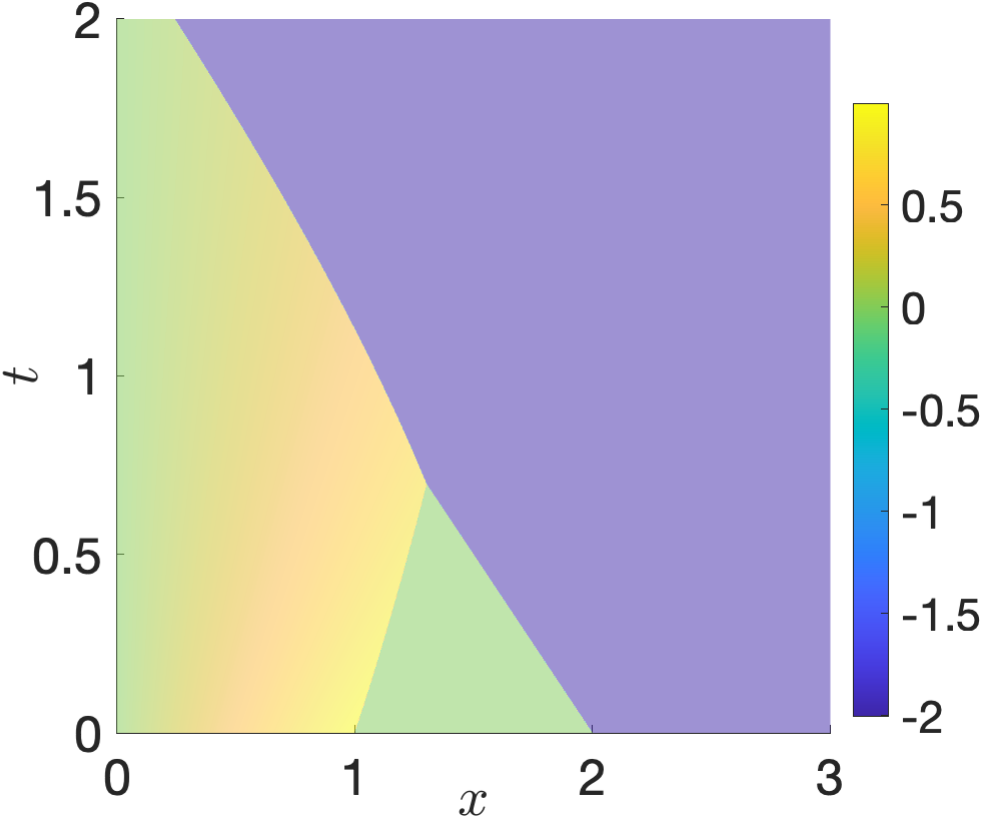}
\caption{$\check{u}(x,t)$ }
\label{Exp6-fig-b}
\end{subfigure}
\hspace{0.15cm}
\begin{subfigure}{0.32\textwidth}
\centering
\includegraphics[width=0.88\textwidth]{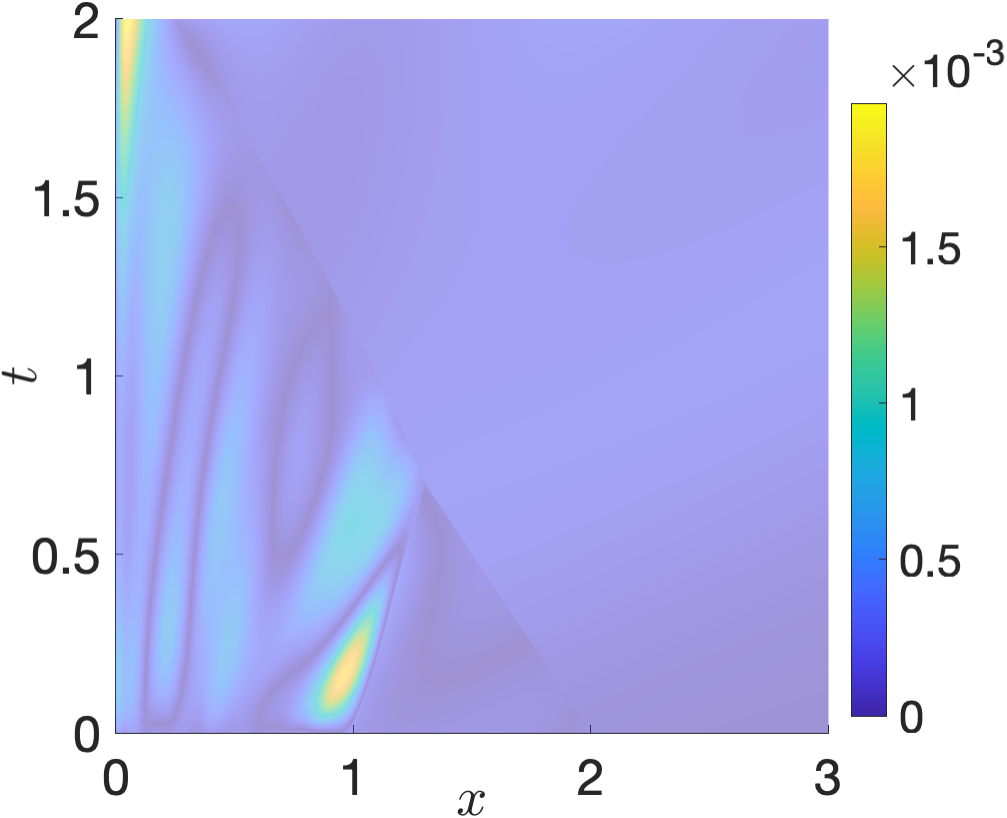}
\caption{$|\check{u}(x,t) - u(x,t)|$}
\label{Exp6-fig-c}
\end{subfigure}
\vspace{0.3cm}

\centering
\begin{subfigure}{0.32\textwidth}
\centering
\includegraphics[width=\textwidth]{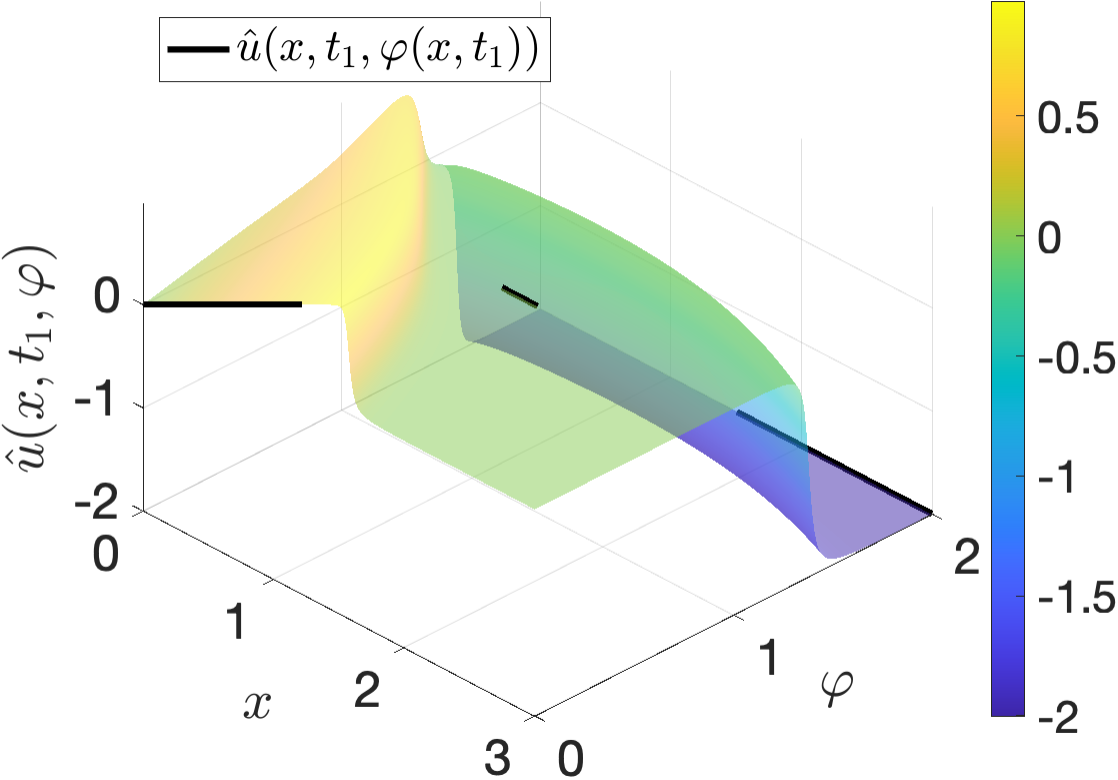}
\caption{$\hat{u}(x,t_1,\varphi)$}
\label{Exp6-fig-d}
\end{subfigure}
\hspace{0.15cm}
\begin{subfigure}{0.32\textwidth}
\centering
\includegraphics[width=0.9\textwidth]{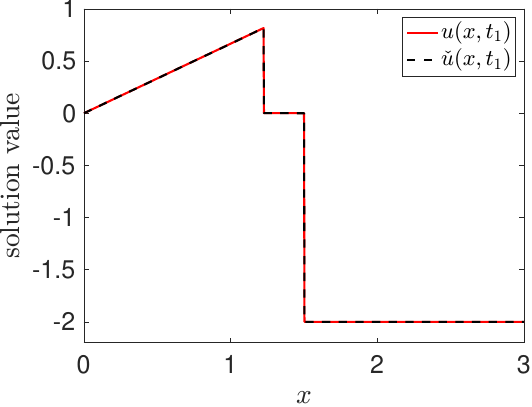}
\caption{$\check{u}(x,t_1)$ and $u(x,t_1)$ }
\end{subfigure}
\hspace{0.15cm}
\begin{subfigure}{0.32\textwidth}
\centering
\includegraphics[width=0.9\textwidth]{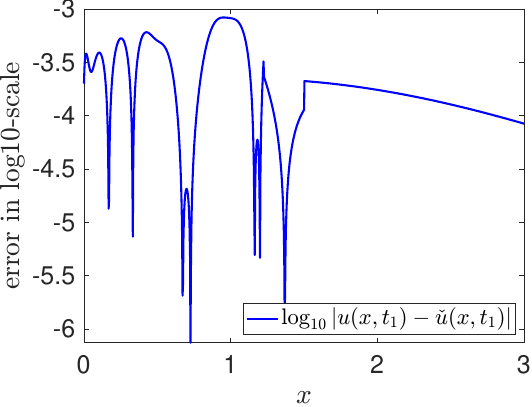}
\caption{$\log_{10}|\check{u}(x,t_1) - u(x,t_1)|$}
\label{Exp6-fig-f}
\end{subfigure}
\vspace{0.3cm}

\centering
\begin{subfigure}{0.32\textwidth}
\centering
\includegraphics[width=\textwidth]{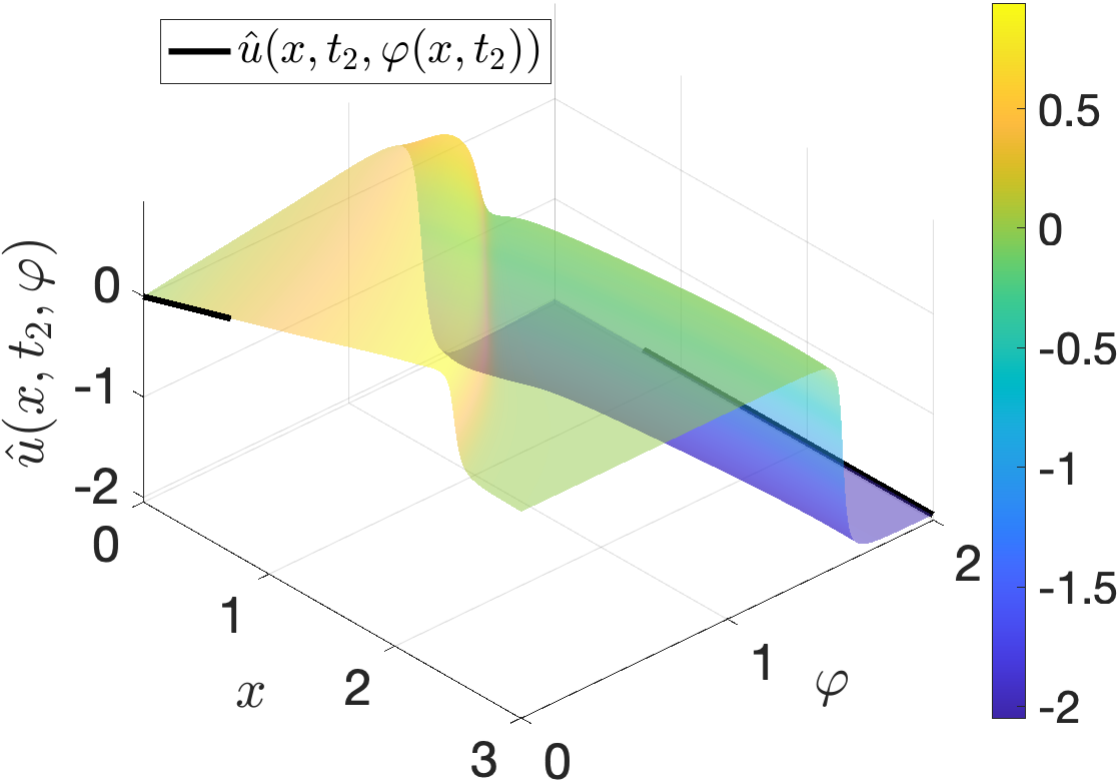}
\caption{$\hat{u}(x,t_2,\varphi)$}
\label{Exp6-fig-g}
\end{subfigure}
\hspace{0.15cm}
\begin{subfigure}{0.32\textwidth}
\centering
\includegraphics[width=0.9\textwidth]{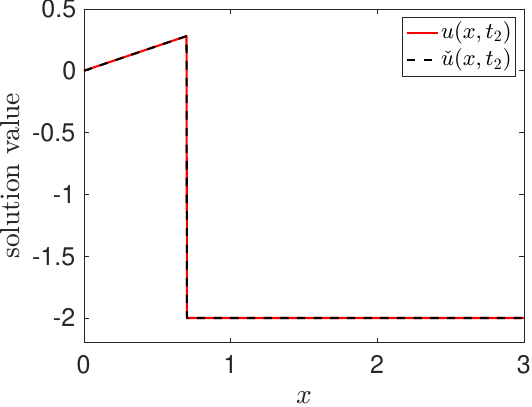}
\caption{$\check{u}(x,t_2)$ and $u(x,t_2)$ }
\end{subfigure}
\hspace{0.15cm}
\begin{subfigure}{0.32\textwidth}
\centering
\includegraphics[width=0.9\textwidth]{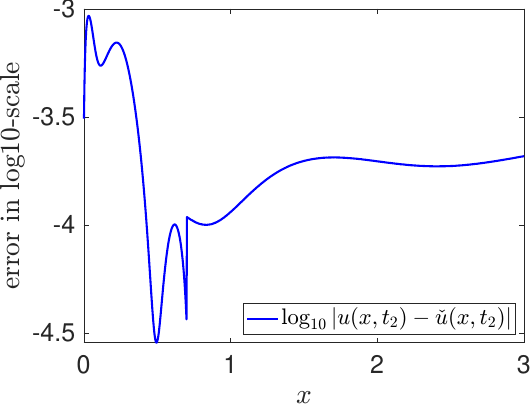}
\caption{$\log_{10}|\check{u}(x,t_2) - u(x,t_2)|$}
\label{Exp6-fig-i}
\end{subfigure}
\caption{Numerical results for the Burgers' equation \eqref{Exp6-Eqns-u-exact} with two shocks coalesce into one ($t_1=0.5$, $t_2=1$).}
\label{Exp6-fig-1d-Burger-disCsolu-2shocks}
\vspace{-0.3cm}
\end{figure}
%---------------------------------------%

To begin with, we consider the inviscid Burgers' equation in the one dimension, i.e.,
\begin{equation}
\left\{
\begin{array}{ll}
\displaystyle \partial_t u(x,t) + u(x,t) \partial_x u(x,t) = 0, \ &\ \ \ \textnormal{for}\ \ (x,t)\in \Omega = (0, 3)\times (0, 2],\\
\displaystyle u_0(x) = \left\{
    \begin{array}{l}
        x,\\
        0, \\
        -2,
    \end{array}\right.
 \ & \ \ \,
    \begin{array}{l}
       \textnormal{for} \ \ 0\leq x<1,\\
       \textnormal{for} \ \ 1 \leq x < 2, \\
       \textnormal{for} \ \ 2\leq x\leq 3,
    \end{array} \\
\displaystyle u(0,t) = 0,\  u(3,t) = -2,\ & \ \ \ \textnormal{for}\ \ t\in (0, 2], 
\end{array}\right.
\label{Exp6-Eqns-u-exact}
\end{equation}
which can be decomposed into two Riemann problems at $x_0=1$ and $x_0=2$. As depicted in \autoref{fig-characteristics-inviscid-burgers-2}, the characteristic line emanating from $x_0\in (0,1)$ takes on the form
\begin{equation*}
	x(t) = u_0(x_0)t + x_0 = x_0(1+t)\ \ \ \textnormal{with} \ \ \ u(x,t) = u_0(x_0) = \frac{x}{1+t},
\end{equation*}
while that departing from $x_0\in (1,0)$ is vertical with $u(x,t) = u_0(x_0) = 0$. Consequently, the shock curve $x = \gamma_1(t)$ emerging from $x_0=1$ must satisfies the Rankine-Hugoniot jump condition (see \autoref{Thm-Speed-Discontinuities}), i.e.,
\begin{equation*}
    s_1(t) = \frac{d\gamma_1(t)}{dt} = \frac{1}{2} \left( \frac{\gamma_1(t)}{1+t} + 0 \right) \ \ \ \textnormal{with} \ \ \ \gamma_1(0)=1,
\end{equation*}
which implies that the first shock curve reads $\gamma_1(t)=\sqrt{1+t}$. On the other hand, the second shock curve $x = \gamma_2(t)$ emerging from $x_0=2$ is given by
\begin{equation*}
     s_2(t) = \frac{d\gamma_2(t)}{dt} = \frac12 (0-2) \ \ \ \textnormal{with} \ \ \ \gamma_2(0)=2,\ \ \ \textnormal{or equivalently,}\ \ \ \gamma_2(t) = -t + 2,
\end{equation*}
which will reach the right-propagating shock $\gamma_1(t)$ at the point $(x^*,t^*) = (\frac{\sqrt{13}-1}{2}, \frac{5-\sqrt{13}}{2})$.

At the moment when two shocks collide, we obtain
\begin{equation*}
	s_3(t) = \frac{d\gamma_3(t)}{dt} = \frac{1}{2} \left( \frac{\gamma_3(t)}{1+t} - 2\right) \ \ \ \textnormal{with} \ \ \ \gamma_3(t^*) = x^*,
\end{equation*}
or equivalently, $\gamma_3(t) = \sqrt{13(1+t)} - 2(1+t)$, which indicates that two shocks merge to become a single shock at later times. Now, we are prepared to embed all these jump conditions by incorporating an augmented variable
\begin{equation*}
\varphi(x, t) = \left\{
\begin{array}{ll}
H( x - \gamma_1(t) ) + H (x - \gamma_2(t) ), \ &\ \ \ \textnormal{for}\ \ 0\leq  t\leq t^*,\\
2H( x - \gamma_3(t) ), \ &\ \ \ \textnormal{for}\ \ t^* < t \leq 2, 
\end{array}\right.
\end{equation*}
and the hyperbolic system satisfied by our solution ansatz $\hat{u}(x,t,\varphi(x,t))$ is then defined as
\begingroup
\renewcommand*{\arraystretch}{1.2}
\begin{equation*}
\left\{
\begin{array}{ll}
\partial_t \hat{u}(x,t,\varphi(x,t)) + \hat{u}(x,t,\varphi(x,t)) \partial_x \hat{u}(x,t,\varphi(x,t)) = 0, \ &\ \ \ \textnormal{for}\ \ (x,t)\in \Omega \setminus \bigcup_{i=1}^3 \Gamma_i,\\
\hat{u}(x,t,\varphi^+(x,t)) + \hat{u}(x,t,\varphi^-(x,t)) = 2 s_i (t),  \ &\ \ \ \textnormal{for}\ \ (x,t) \in \Gamma_i\ (1\leq i\leq 3),\\
\hat{u}(x,0,\varphi(x,0)) = u_0(x), \ &\ \ \ \textnormal{for}\ \ x\in (0,3), \\
\hat{u}(0,t,\varphi(0,t)) = 0,\  \hat{u}(3, t, \varphi(3, t)) = -2,\ &\ \ \ \textnormal{for}\ \ t\in (0,2].
\end{array}\right.
\end{equation*}
\endgroup

Numerical results obtained through our lift-and-embed learning approach are displayed in \autoref{Exp6-fig-1d-Burger-disCsolu-2shocks}, which are in good agreement with the exact solution. Notably, shock waves of \eqref{Exp6-Eqns-u-exact} not only intersect with each other but also appear as curved interfaces, which often incur higher costs for traditional numerical schemes requiring interface-fitted meshes. Thanks to the meshless nature of neural networks, our method can effectively adapt to irregular shock curves (see \autoref{Exp6-fig-b} and \autoref{Exp6-fig-c}), making it a compelling choice for scenarios with complex geometries of shock trajectories.

Additionally, as depicted in \autoref{Exp6-fig-d} and \autoref{Exp6-fig-g}, temporal slices of trained network solution are now situated on surfaces within a one-order higher-dimensional space. All solution jumps of \eqref{Exp6-Eqns-u-exact} are reconstructed as differences between distinct collocation points, which are effectively reconstructed using a fully-connected neural network with smooth activation functions (see \autoref{Exp6-fig-f} and \autoref{Exp6-fig-i}). 

%%%%%%%%%%%%%%%%%%%%%%%%%%%%%%%%%%%%%%%%%%%%%%%%%%%%%
\subsubsection{Rarefaction-Shock Interaction}\label{sec-ExpSM1-1d-Burgers-Rarefaction}

In this numerical experiment, we consider another inviscid Burgers' equation, i.e.,
\begingroup
\renewcommand*{\arraystretch}{1.1}
\begin{equation}
	\left\{
	\begin{array}{ll}
		\partial_t u(x,t) + \frac12 u(x, t) \partial_x u(x, t) = 0, \ &\ \ \ \textnormal{for}\ \ (x,t)\in \Omega = (-1, 6)\times (0, 10],\\
		\displaystyle u_0(x) = \left\{
		\begin{array}{l}
			0,\\
			1, \\
			0,
		\end{array}\right.
		\ & \ \ \,
		\begin{array}{l}
			\textnormal{for} \ \ -1\leq x\leq 0,\\
			\textnormal{for} \ \ 0 <  x < 1, \\
			\textnormal{for} \ \ 1\leq x\leq 6,
		\end{array} \\
		\displaystyle u(-1,t) = u(6,t) = 0,\ & \ \ \ \textnormal{for}\ \ t\in (0, 10],
	\end{array}\right.
	\label{ExpSM1-Eqns-u-exact}
\end{equation}
\endgroup
which is decomposed into two Riemann problems at $x_0 = 0$ and $x_0 = 1$. However, characteristics do not cover the whole plane as shown in \autoref{fig-characteristics-inviscid-burgers-3}, leading to the formation of a rarefaction wave between $x=0$ and $x=\frac12 t$. On the other hand, the shock wave originating form $x_0 = 1$ propagates along a straight line with the speed
\begin{equation*}
	s_1 = \frac{ \llbracket f(u) \rrbracket }{ \llbracket u \rrbracket } = \frac14 \left( u^+(x,t) + u^-(x,t) \right) = \frac{1}{4}\ \ \ \textnormal{where}\ \ \ f(u) = \frac14 u^2.
\end{equation*}

As time progresses, the shock wave encounters the rarefaction wave at $(x,t)=(2,4)$, giving rise to another shock wave with a non-constant velocity. More precisely, this shock wave satisfies the Rankine-Hugoniot condition
\begin{equation*}
	s_2(t) =\frac{d\gamma_2(t)}{dt} = \frac{1}{4} \left( \frac{2\gamma_2(t)}{t} + 0 \right) \ \ \ \textnormal{with} \ \ \ \gamma_2(4)=2,
\end{equation*}
which implies that $\gamma_2(t) = \sqrt{t}$ for $t > 4$. Consequently, our augmented variable can be defined as
\begingroup
\renewcommand*{\arraystretch}{1.3}
\begin{equation*}
	\varphi(x, t) = H(x - \gamma(t)) \ \ \textnormal{with} \ \ \gamma(t) =\left\{
	\renewcommand*{\arraystretch}{1.3}
	\begin{array}{cl}
		\frac{1}{4}t + 1, \ & \ \ \ \textnormal{for}\ \ t \in  [0, 4],\\
		\sqrt{t}, \ &\ \ \ \textnormal{for}\ t \in (4, 10],\\
	\end{array}\right.
\end{equation*}
while the hyperbolic system satisfied by our solution ansatz $\hat{u}(x, t, \varphi(x, t))$ takes on the form
\begin{equation*}
	\left\{
	\begin{array}{ll}
		\partial_t \hat{u}(x,t,\varphi(x,t)) + \frac{1}{2}\hat{u}(x,t,\varphi(x,t)) \partial_x \hat{u}(x,t,\varphi(x,t)) = 0, \ &\ \ \ \textnormal{for}\ \ (x,t)\in \Omega\setminus \bigcup_{i=1}^2 \Gamma_i,\\
		\hat{u}(x,t,\varphi^+(x,t)) + \hat{u}(x,t,\varphi^-(x,t)) = 4s_i(t), \ &\ \ \ \textnormal{for}\ \ (x,t)\in \Gamma_i\ (1\leq i\leq 2),\\
		\hat{u}(x,0,\varphi(x,0)) = u_0(x), \ &\ \ \ \textnormal{for}\ \ x\in (-1,6), \\
		\hat{u}(-1,t,\varphi(-1,t)) = \hat{u}(6, t, \varphi(6, t)) = 0,\ &\ \ \ \textnormal{for}\ t\in (0,10].
	\end{array}\right.
\end{equation*}
\endgroup

\begin{figure}[!htbp]
\centering
\begin{subfigure}{0.32\textwidth}
\centering
\includegraphics[width=\textwidth]{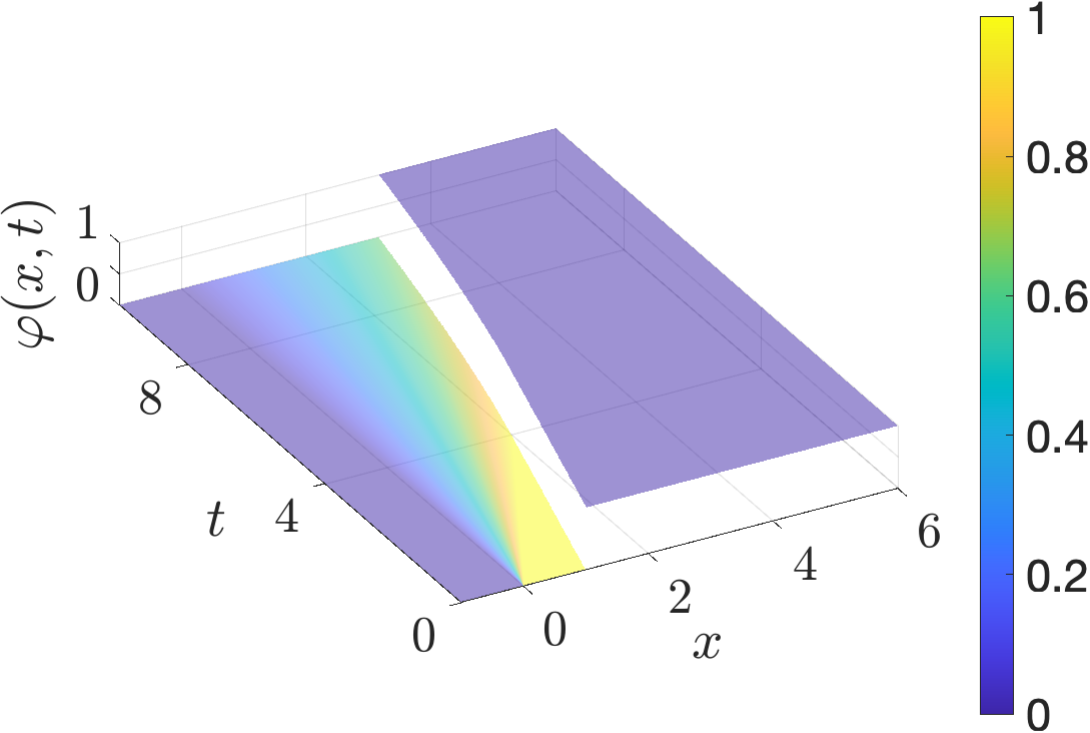}
\caption{$\hat{u}(x,t,\varphi(x,t))$}
\end{subfigure}
\hspace{0.15cm}
\begin{subfigure}{0.29\textwidth}
\centering
\includegraphics[width=0.95\textwidth]{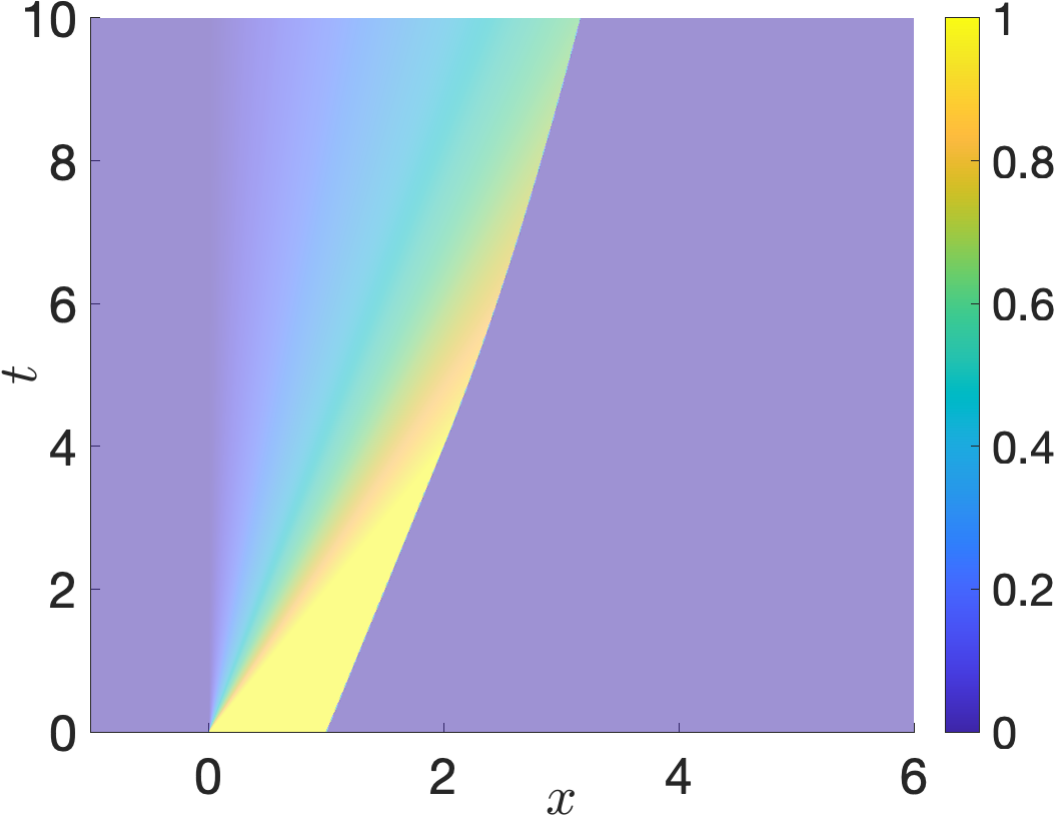}
\caption{$\check{u}(x,t)$ }
\label{ExpSM1-fig-b}
\end{subfigure}
\hspace{0.15cm}
\begin{subfigure}{0.3\textwidth}
\centering
\includegraphics[width=0.9\textwidth]{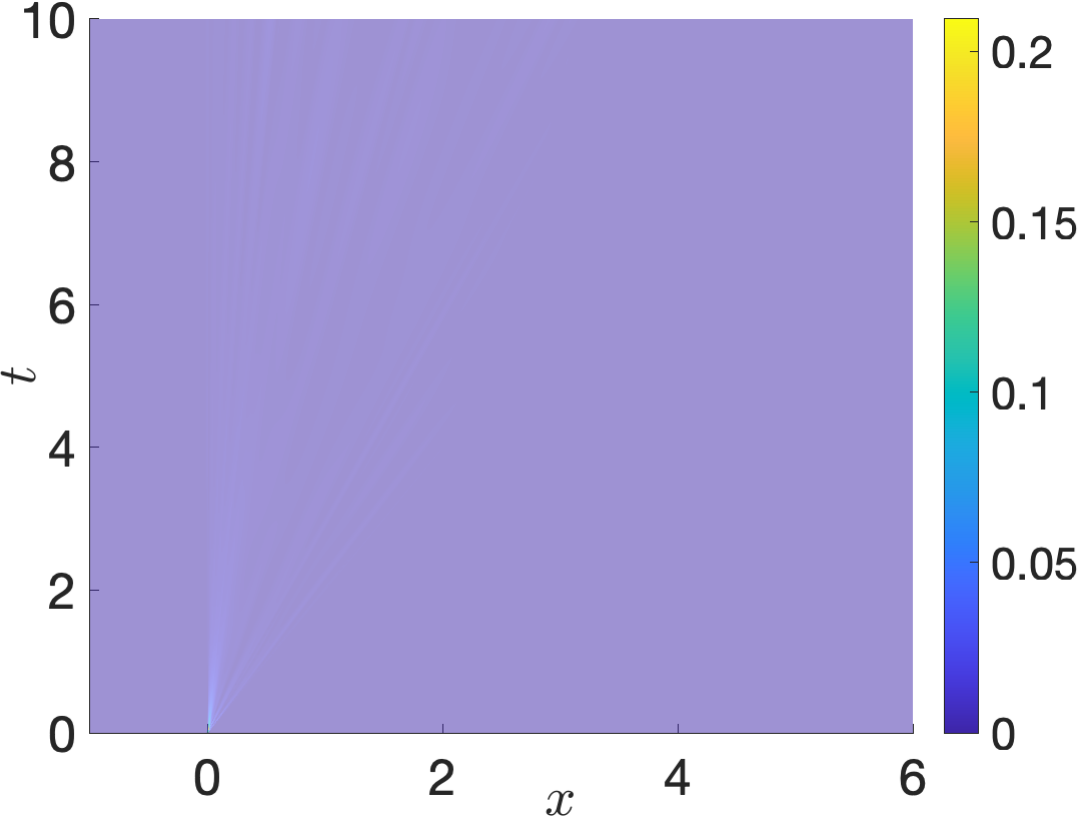}
\caption{$|\check{u}(x,t) - u(x,t)|$}
\label{ExpSM1-fig-c}
\end{subfigure}
\vspace{0.3cm}

\centering
\begin{subfigure}{0.32\textwidth}
\centering
\includegraphics[width=\textwidth]{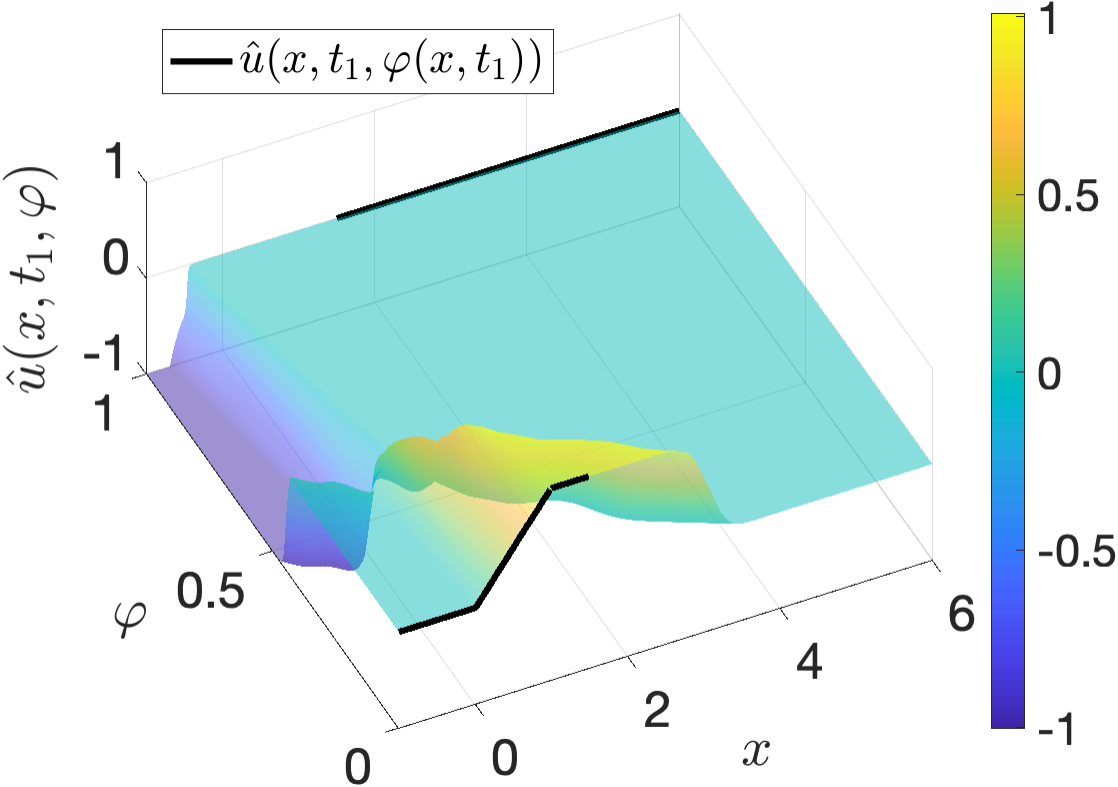}
\caption{$\hat{u}(x,t_1,\varphi)$}
\label{ExpSM1-fig-d}
\end{subfigure}
\hspace{0.15cm}
\begin{subfigure}{0.32\textwidth}
\centering
\includegraphics[width=0.9\textwidth]{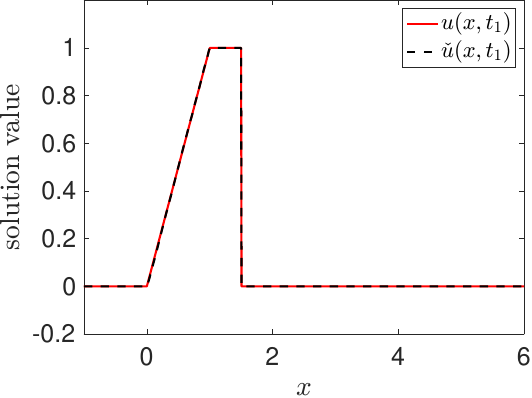}
\caption{$\check{u}(x,t_1)$ and $u(x,t_1)$ }
\end{subfigure}
\hspace{0.15cm}
\begin{subfigure}{0.32\textwidth}
\centering
\includegraphics[width=0.9\textwidth]{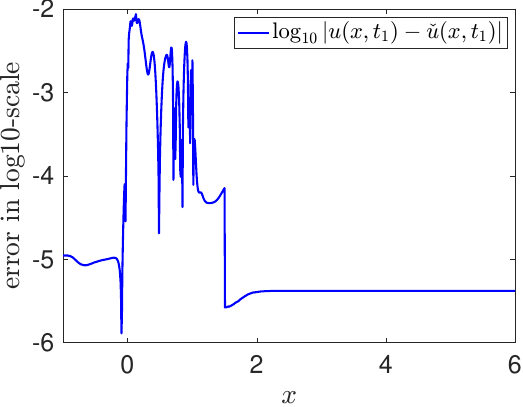}
\caption{$\log_{10}|\check{u}(x,t_1) - u(x,t_1)|$}
\label{ExpSM1-fig-f}
\end{subfigure}
\vspace{0.3cm}

\centering
\begin{subfigure}{0.32\textwidth}
\centering
\includegraphics[width=\textwidth]{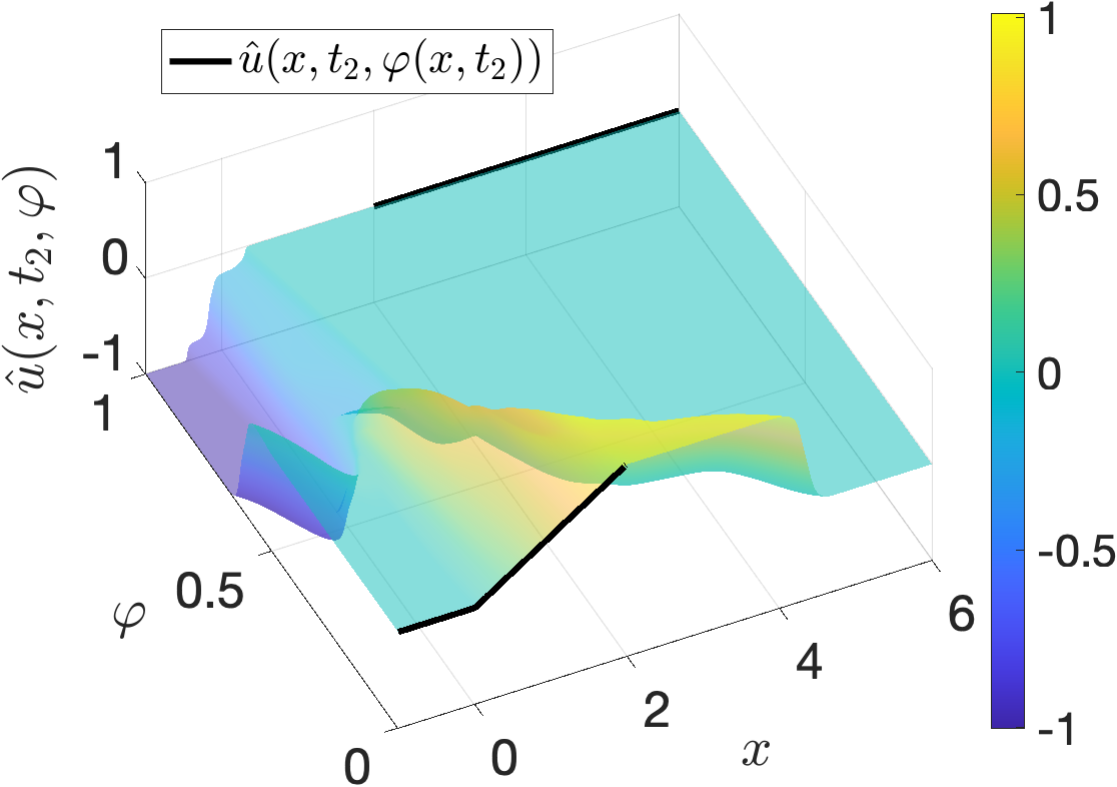}
\caption{$\hat{u}(x,t_2,\varphi)$}
\label{ExpSM1-fig-g}
\end{subfigure}
\hspace{0.15cm}
\begin{subfigure}{0.32\textwidth}
\centering
\includegraphics[width=0.9\textwidth]{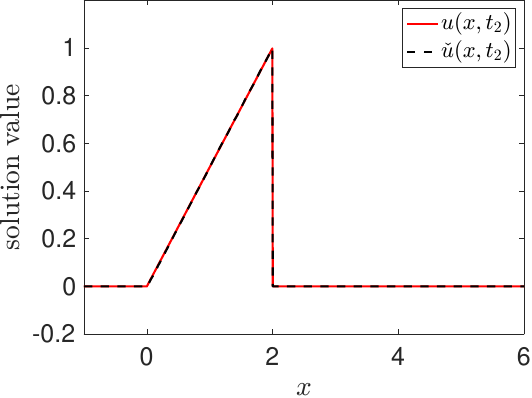}
\caption{$\check{u}(x,t_2)$ and $u(x,t_2)$ }
\end{subfigure}
\hspace{0.15cm}
\begin{subfigure}{0.32\textwidth}
\centering
\includegraphics[width=0.9\textwidth]{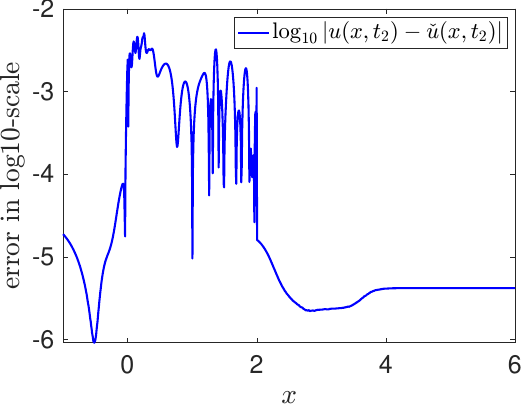}
\caption{$\log_{10}|\check{u}(x,t_2) - u(x,t_2)|$}
\label{ExpSM1-fig-i}
\end{subfigure}
\vspace{0.3cm}

\centering
\begin{subfigure}{0.32\textwidth}
\centering
\includegraphics[width=\textwidth]{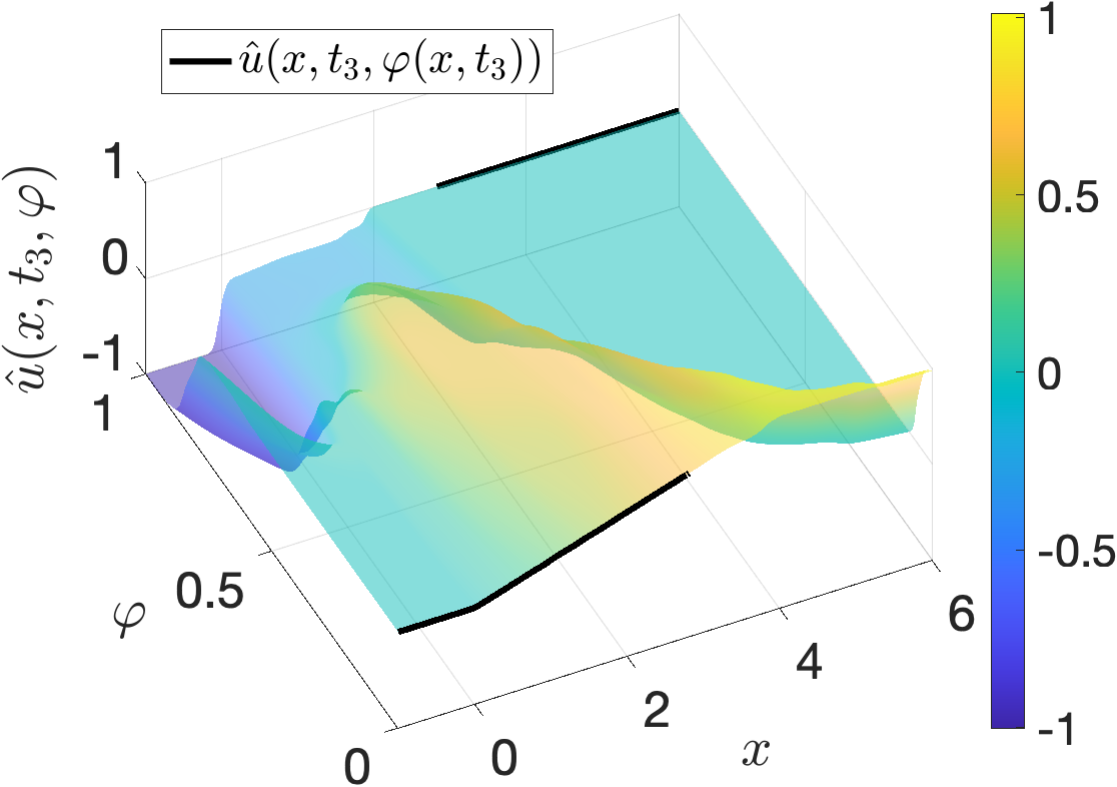}
\caption{$\hat{u}(x,t_3,\varphi)$}
\label{ExpSM1-fig-j}
\end{subfigure}
\hspace{0.15cm}
\begin{subfigure}{0.32\textwidth}
\centering
\includegraphics[width=0.9\textwidth]{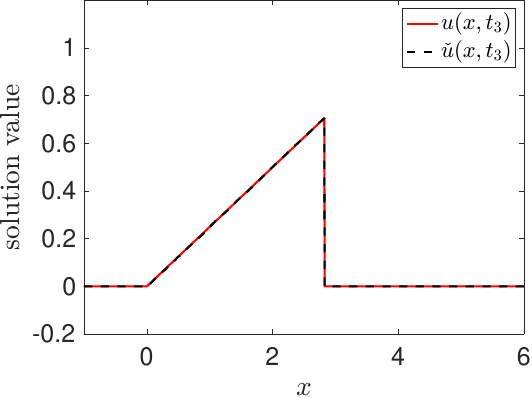}
\caption{$\check{u}(x,t_3)$ and $u(x,t_3)$ }
\end{subfigure}
\hspace{0.15cm}
\begin{subfigure}{0.32\textwidth}
\centering
\includegraphics[width=0.9\textwidth]{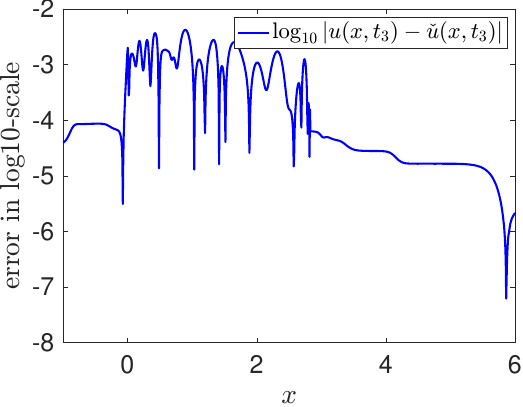}
\caption{$\log_{10}|\check{u}(x,t_3) - u(x,t_3)|$}
\label{ExpSM1-fig-l}
\end{subfigure}
\caption{Numerical results for the Burgers' equation \eqref{Exp6-Eqns-u-exact} with rarefaction-shock interaction ($t_1\!=\!2$, $t_2\!=\!4$, $t_3\!=\!8$).}
\label{ExpSM1-fig-1d-Burger-rarefaction}
\vspace{-0.3cm}
\end{figure}

Numerical results obtained through our lift-and-embed learning approach are shown in \autoref{ExpSM1-fig-1d-Burger-rarefaction}, demonstrating a satisfactory approximation of both rarefaction and shock waves. During the training process, an additional penalty coefficient $\beta= 300$ is assigned to the interior loss term \eqref{LE-Loss-Forward} in this case. Notably, the approximation error associated with the rarefaction wave is comparatively higher than that observed in other regions (see \autoref{ExpSM1-fig-f}, \autoref{ExpSM1-fig-i}, and \autoref{ExpSM1-fig-l} for instance), which is likely to be further improved by incorporating entropy conditions into the training loss function \cite{de2024wpinns} (a subject left for future investigation).

%---------------------------------------%
%%---------------------------------------%

%%%%%%%%%%%%%%%%%%%%%%%%%%%%%%%%%%%%%%%%%%%%%%%%%%%%%
\subsubsection{Two-Dimensional Inviscid Burgers' Equation}\label{sec-Exp7-2d-Burgers}

%---------------------------------------%
\begin{figure}[t!]
	\centering
	\begin{subfigure}{0.32\textwidth}
		\centering
		\includegraphics[width=\textwidth]{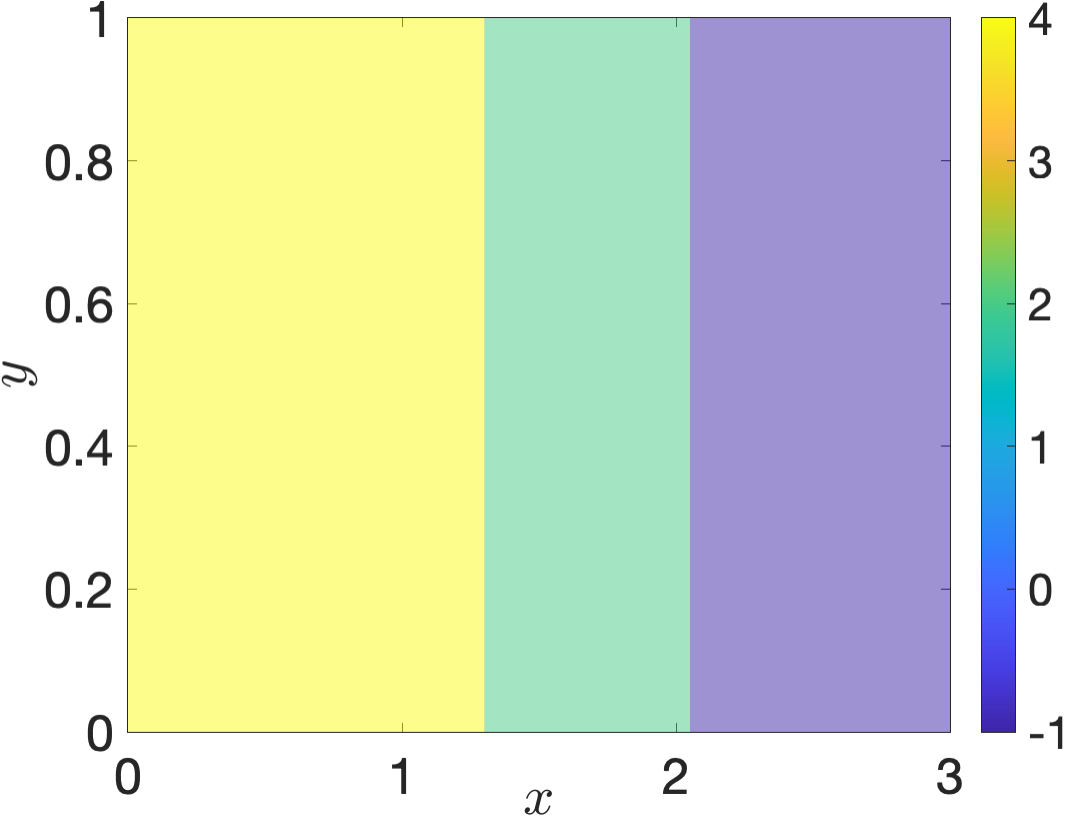}
		\caption{$u(x,y, t_1)$}
	\end{subfigure}
	\hspace{0.1cm}
	\begin{subfigure}{0.32\textwidth}
		\centering
		\includegraphics[width=\textwidth]{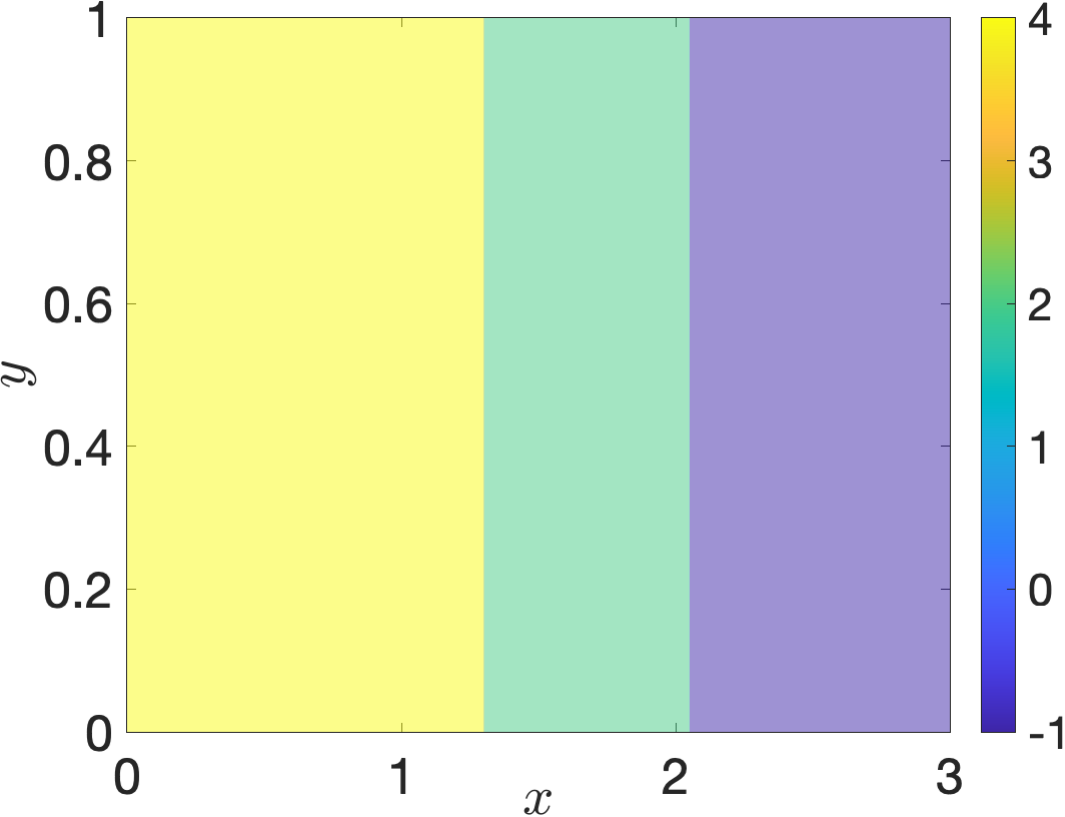}
		\caption{$\check{u}(x,y, t_1)$ }
	\end{subfigure}
	\hspace{0.1cm}
	\begin{subfigure}{0.33\textwidth}
		\centering
		\includegraphics[width=\textwidth]{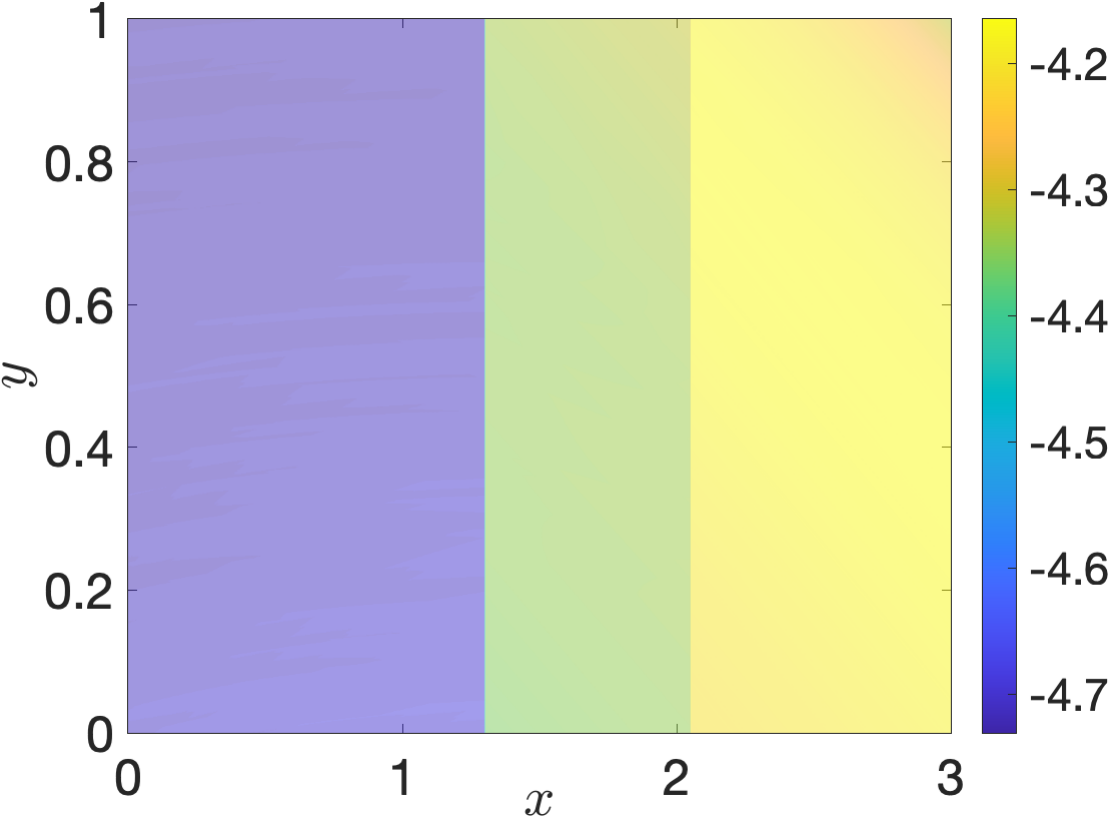}
		\caption{$\log_{10}|\check{u}(x,y,t_1) - u(x,y,t_1)|$}
	\end{subfigure}
	\vspace{0.3cm}
		
	\centering
	\begin{subfigure}{0.32\textwidth}
		\centering
		\includegraphics[width=\textwidth]{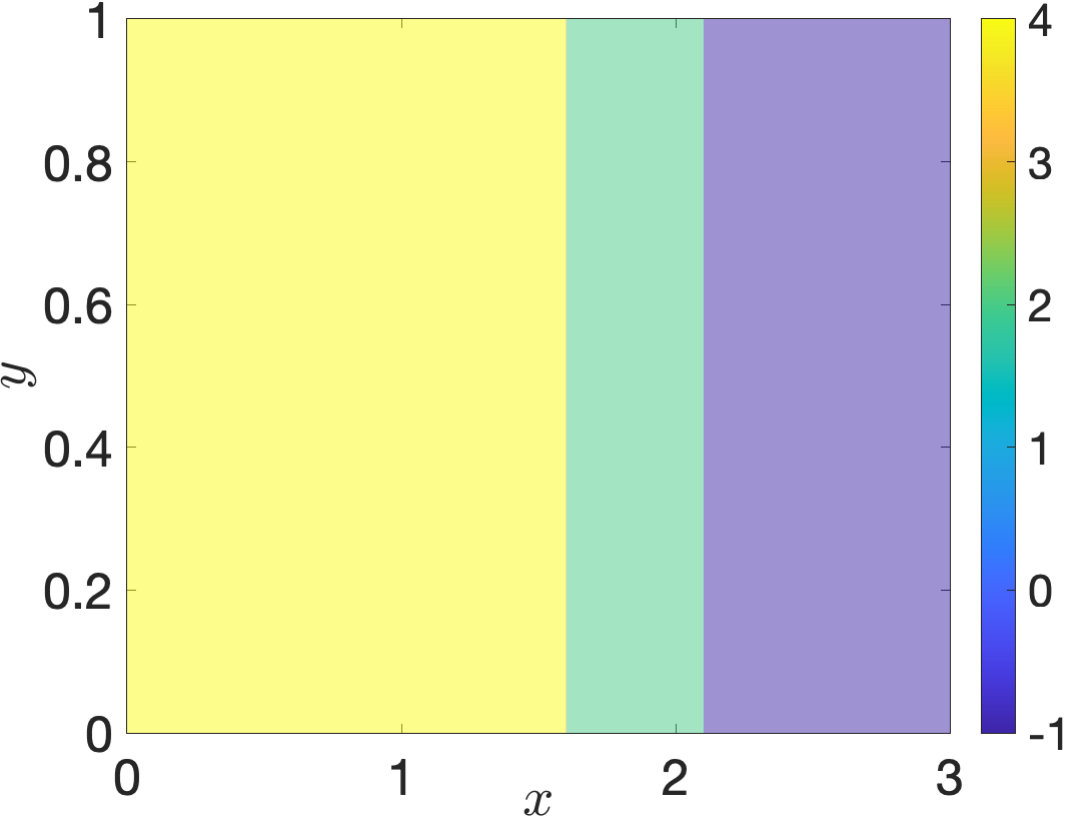}
		\caption{$u(x,y,t_2)$}
	\end{subfigure}
	\hspace{0.1cm}
	\begin{subfigure}{0.32\textwidth}
		\centering
		\includegraphics[width=\textwidth]{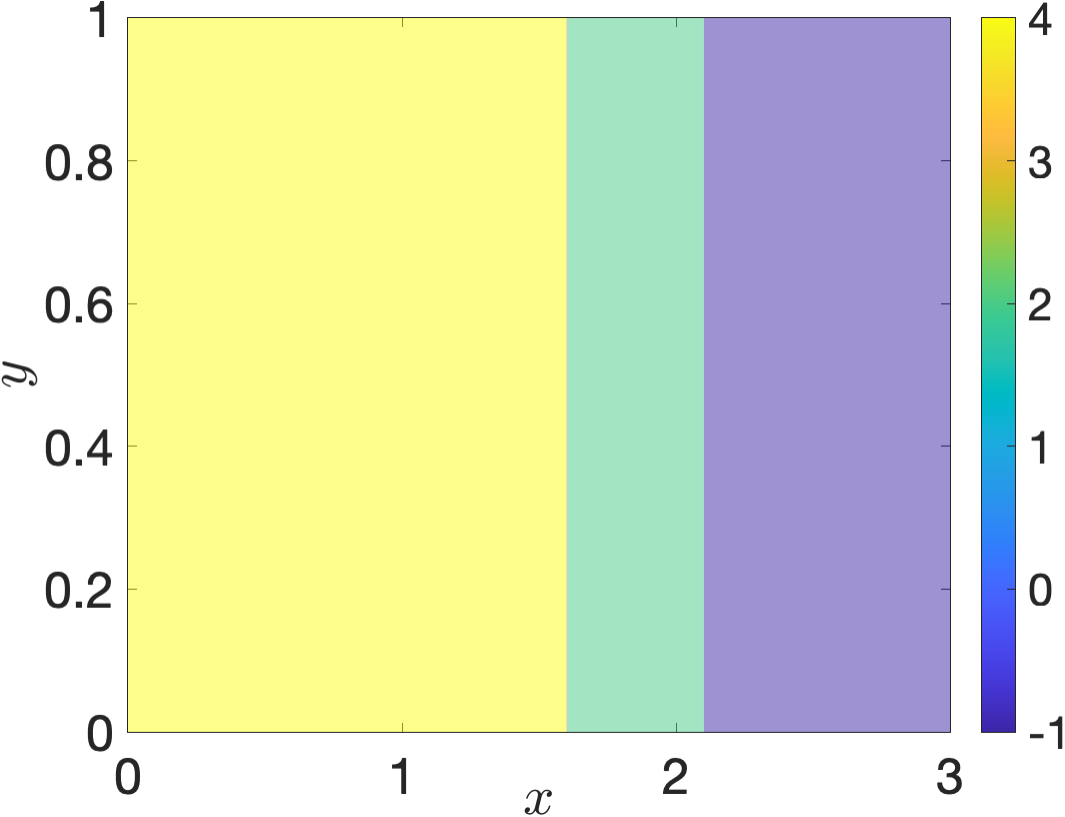}
		\caption{$\check{u}(x, y, t_2)$}
	\end{subfigure}
	\hspace{0.1cm}
	\begin{subfigure}{0.33\textwidth}
		\centering
		\includegraphics[width=\textwidth]{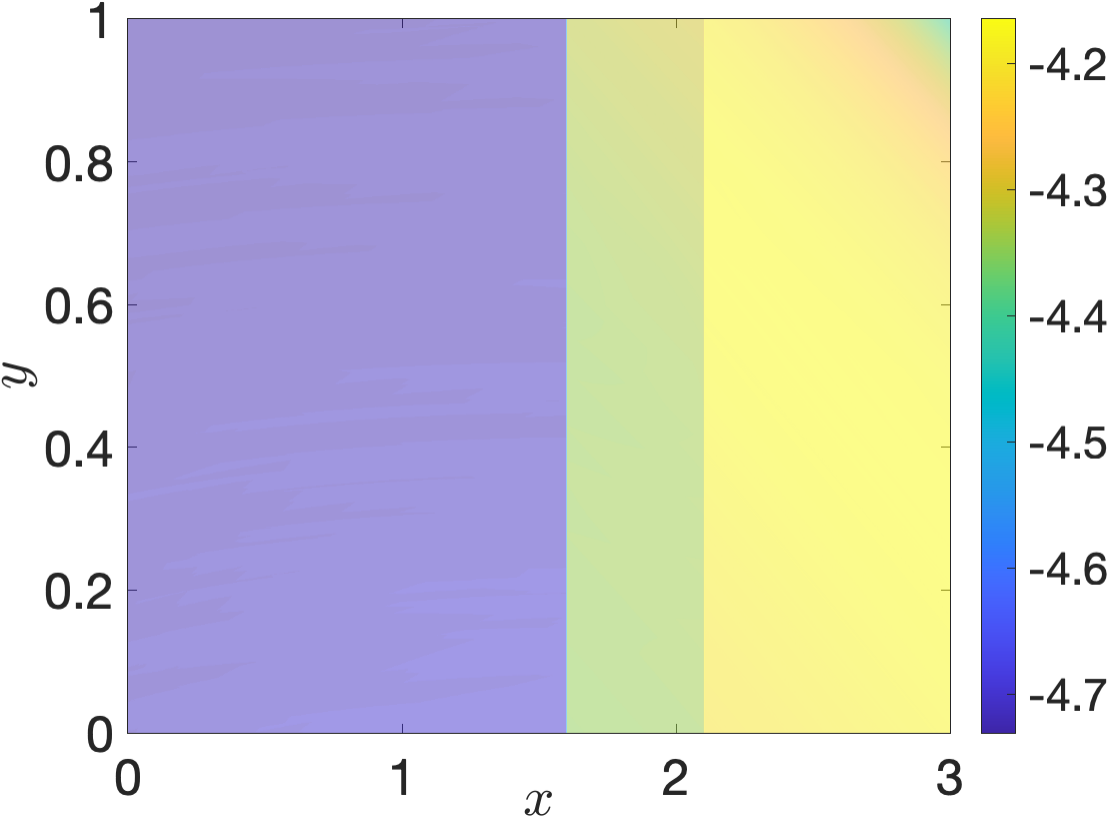}
		\caption{$\log_{10}|\check{u}(x, y, t_2) - u(x,y, t_2)|$}
	\end{subfigure}
	\vspace{0.3cm}
		
	\centering
	\begin{subfigure}{0.32\textwidth}
		\centering
		\includegraphics[width=\textwidth]{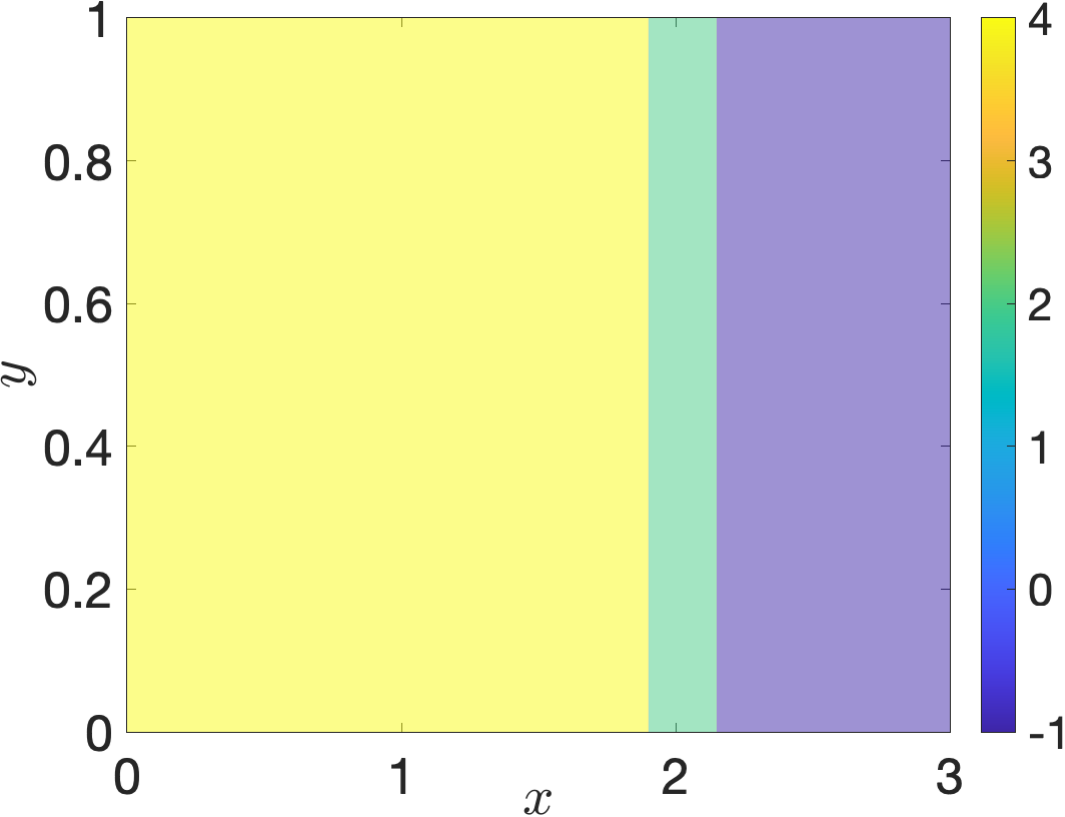}
		\caption{$u(x,y,t_3)$}
	\end{subfigure}
	\hspace{0.1cm}
	\begin{subfigure}{0.32\textwidth}
		\centering
		\includegraphics[width=\textwidth]{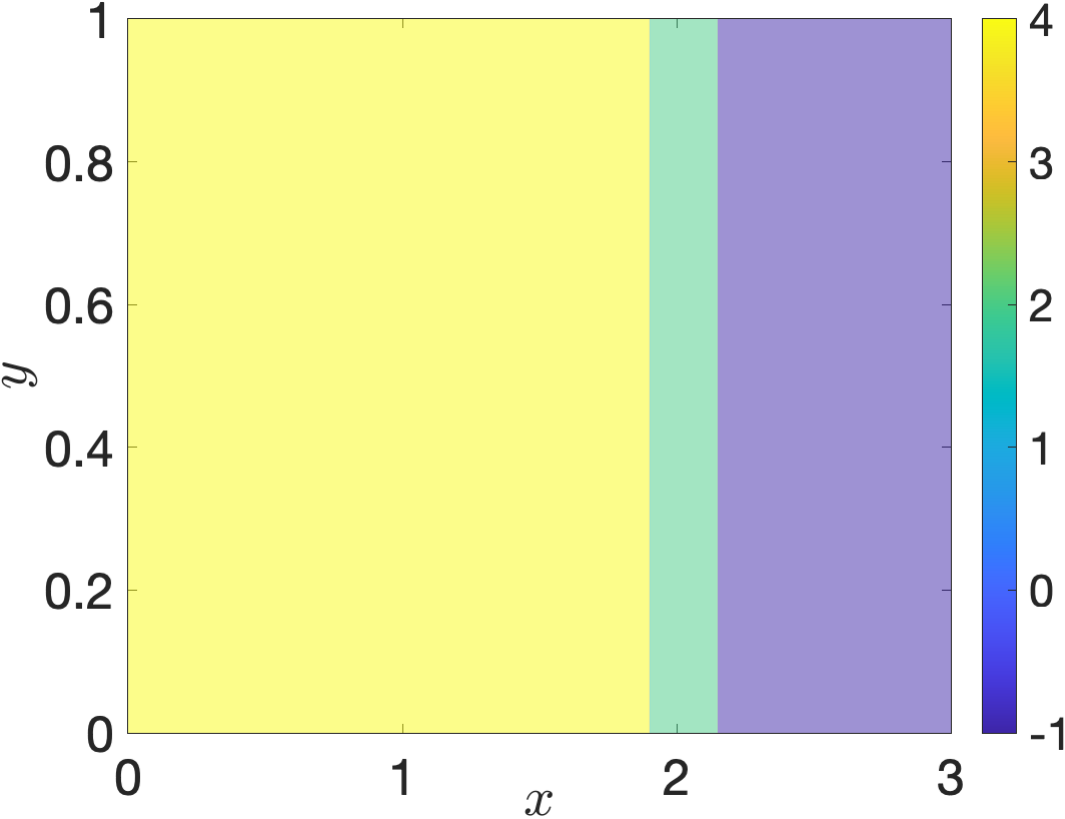}
		\caption{$\check{u}(x,y,t_3)$}
	\end{subfigure}
	\hspace{0.1cm}
	\begin{subfigure}{0.33\textwidth}
		\centering
		\includegraphics[width=\textwidth]{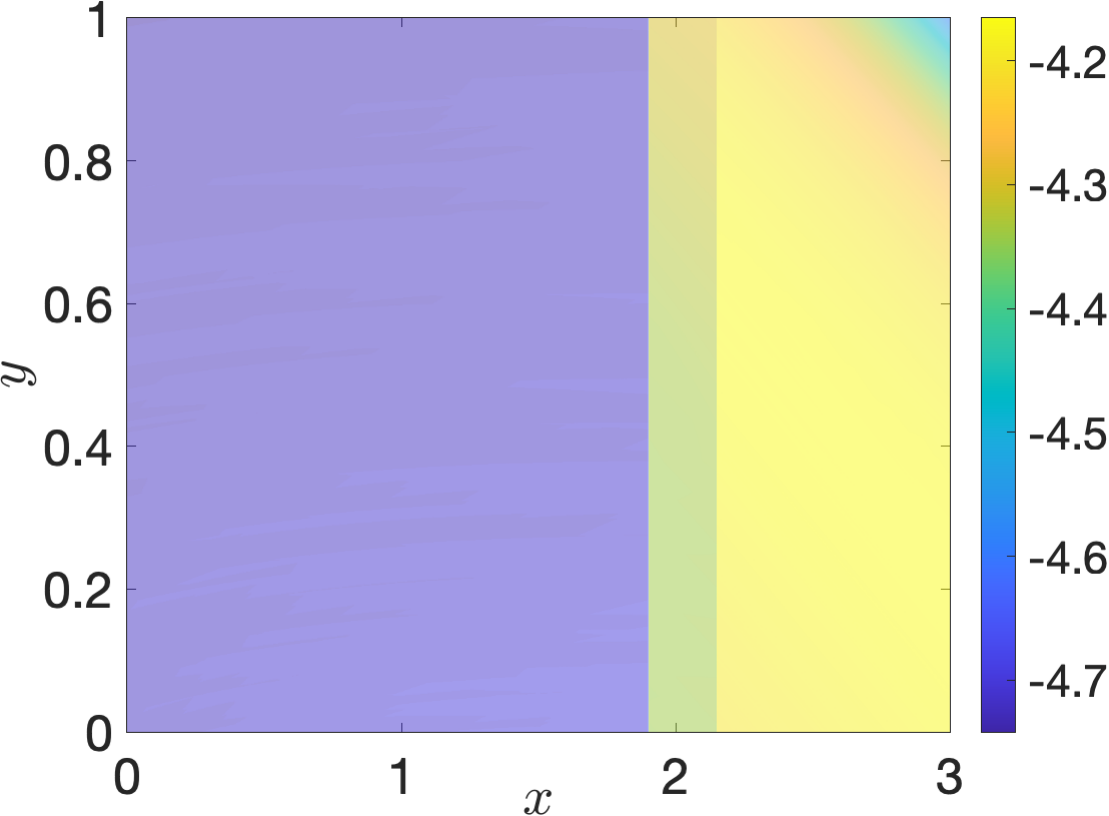}
		\caption{$\log_{10}|\check{u}(x,y,t_3) - u(x,y,t_3)|$}
	\end{subfigure}
	\caption{Numerical results for the inviscid Burgers' equation \eqref{Exp7-Eqns-u-exact} in two-dimension ($t_1= 0.1$, $t_2 =0.2$, $t_3=0.3$).}
	\label{Exp7-fig-2d-Burger-disCsolu}
	\vspace{-0.3cm}
	\end{figure}
%%---------------------------------------%

Here, we study a two-dimensional problem to show the capacity in tackling high-dimensional functions, that is, 
\begin{equation}
	\left\{
	\begin{array}{ll}
		\partial_t u(x,y,t) + u(x,y,t) \partial_x u(x,y,t) + u(x,y,t) \partial_y u(x,y,t)  = 0, \ &\ \ \ \textnormal{for}\ (x,y,t)\in \Omega = (0,3)\times(0,1)\times (0,0.4],\\
		\displaystyle u_0(x) = \left\{
		\begin{array}{l}
			4,\\
			2, \\
			-1,
		\end{array}\right.
		\ & \ \ \,
		\begin{array}{l}
			\textnormal{for} \ \ 0\leq x<1,\\
			\textnormal{for} \ \ 1 \leq x < 2, \\
			\textnormal{for} \ \ 2\leq x\leq 3,
		\end{array} \\
	\end{array}\right.
	\label{Exp7-Eqns-u-exact}
\end{equation}
with inflow boundary conditions derived from the exact solution 
\begin{equation*}
	u(x,y,t) = 4H(1+3t-x) + 2 H(x-1-3t) - 3 H(x-\frac12t-2).
\end{equation*}
We denote by $\Gamma_i = \Gamma_i(x,y,t)$ the surface of discontinuities departing from the initial position $\{ (x_i,y,0) \, | \, 0\leq y\leq 1 \}$ and $\bm{n}_i = (-s_i, \bm{\nu}_i)^T$ the normal vector to it, in which $1\leq i\leq 2$, $x_1=1$, and $x_2=2$. Due to the constancy of initial values, both surfaces are planar where $\bm{n}_1 = (-3, 0, 1)^T$ and $\bm{n}_2 = (-0.5, 0, 1)^T$, and hence can be parametrized as
\begin{equation*}
	x = 1 + 3t \ \ \ \textnormal{and} \ \ \  x = 2 + 0.5t
\end{equation*}
for $0\leq y\leq 1$, respectively. Accordingly, the augmented variable can be constructed as
\begin{equation*}
	\varphi(x, y, t) = H(x - 2 - 0.5t) +  H(x-1 -3t),
\end{equation*}
leading to a four-dimensional function $\hat{u}(x,y,t,\varphi(x,y,t))$ to be resolved from the embedded version of \eqref{Exp7-Eqns-u-exact}, i.e.,
\begingroup
\renewcommand*{\arraystretch}{1.2}
\begin{equation*}
	\left\{
	\begin{array}{ll}
		\partial_t \hat{u}(x, y, t,\varphi(x, y, t)) + 0.5 [ \partial_x +\partial_y ] \hat{u}^2(x, y, t,\varphi(x, y, t)) = 0, \ &\ \ \ \textnormal{for}\ \ (x,y,t)\in\Omega \setminus \bigcup_{i=1}^2 \Gamma_i,\\
		-s_i\llbracket\hat{u}(x,y,t,\varphi(x,y,t))\rrbracket + 0.5 \llbracket \hat{u}^2(x,y,t,\varphi(x,y,t))\rrbracket = 0,  \ &\ \ \ \textnormal{for}\ \ (x,y,t) \in \Gamma_i\ (i = 1, 2),\\
		\hat{u}(x,0,\varphi(x,0)) = u_0(x), \ &\ \ \ \textnormal{for}\ \ x\in (0,3), \\		
	\end{array}\right.
	\label{Eqns-Exp7-u-NN}
\end{equation*}
\endgroup
with inflow boundary conditions. Traditional mesh-based numerical methods are often susceptible to the curse of dimensionality, whereas neural network-based approaches offer a promising way to alleviate this issue.
	
As illustrated in \autoref{Exp7-fig-2d-Burger-disCsolu} and \autoref{table-relative-errors}, our approximate solution shows reliable alignment with the true solution (only specific temporal slices are depicted due to the high dimensionality), which demonstrates the capability of our proposed method in addressing high-dimensional discontinuous solutions.

%%%%%%%%%%%%%%%%%%%%%%%%%%%%%%%%%%%%%%%%%%%%%%%%%%%%%%
\subsection{Inviscid Burgers’ Equations with Unknown Shock Locations}

Unlike previous examples, this section addresses the inviscid Burgers' equation without specifying the location of discontinuities in advance. Instead, Algorithm \ref{Algorithm-LPLM-inverse} is utilized to dynamically infer shock curves during training.

%%%%%%%%%%%%%%%%%%%%%%%%%%%%%%%%%%%%%%%%%%%%%%%%%%%%%
\subsubsection{The Benchmark Inviscid Burgers’ Problem Revisited}\label{sec-Exp8-1d-iBurgers}

Here, we re-examine the inviscid Burgers' equation \eqref{Exp2-Eqns-u-exact} in which the shock curve takes the form of a straight line $\gamma(t) = st$. In contrast to the previous study, the shock speed $s\in\mathbb{R}$, as integrated into the formulation of our ansatz \eqref{Exp2-uNN-ansatz}, is regarded as an unknown parameter to be inferred. Specifically, the empirical loss function \eqref{Loss-Shock-Update} now reads
\begin{equation*}
L_{\textnormal{Shock}}^{\textnormal{inv}} (\theta, \hat{s}) =  L_{\textnormal{Shock}} (\theta) + \frac{1}{n} \sum_{i=1}^{n} \left| \hat{s} \llbracket \hat{u}(\gamma(t_i),t_i, \varphi(\gamma(t_i), t_i); \theta) \rrbracket - \llbracket f (\hat{u}(\gamma(t_i),t_i, \varphi(\gamma(t_i), t_i); \theta)) \rrbracket \right|^2.
\end{equation*}

\begin{figure}[t!]
\centering
\begin{subfigure}{0.32\textwidth}
\centering
\includegraphics[width=\textwidth]{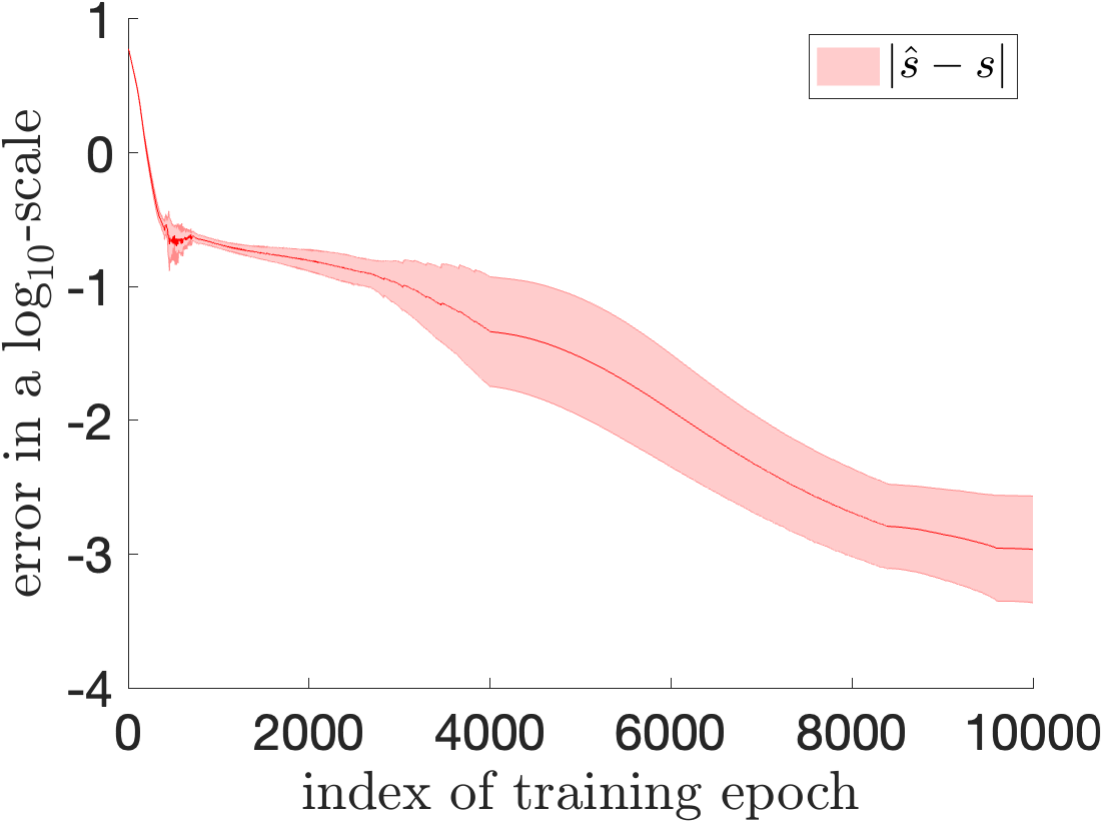}
\caption{$\log_{10}^{|\hat{s} -s |}$ (mean $\pm$ deviation)}
\end{subfigure}
\hspace{0.15cm}
\begin{subfigure}{0.32\textwidth}
\centering
\includegraphics[width=\textwidth]{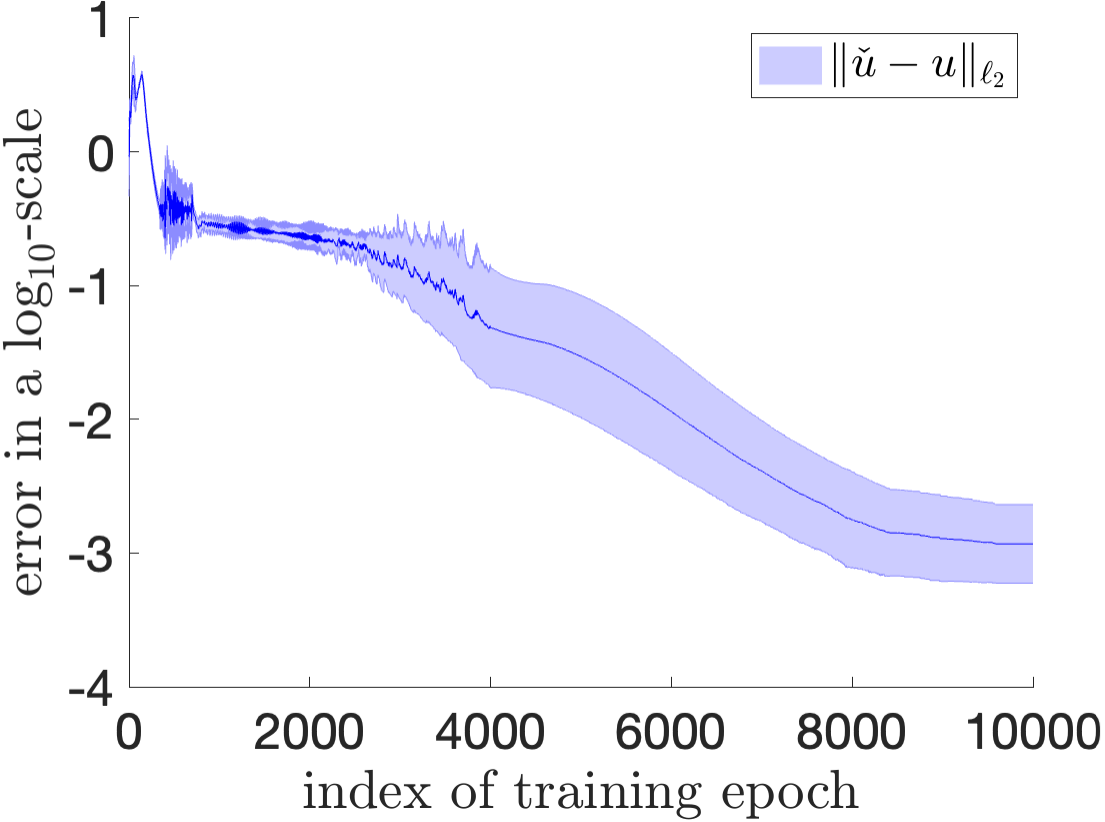}
\caption{test loss (mean $\pm$ deviation)}
\end{subfigure}
\hspace{0.15cm}
\begin{subfigure}{0.32\textwidth}
\centering
\includegraphics[width=\textwidth]{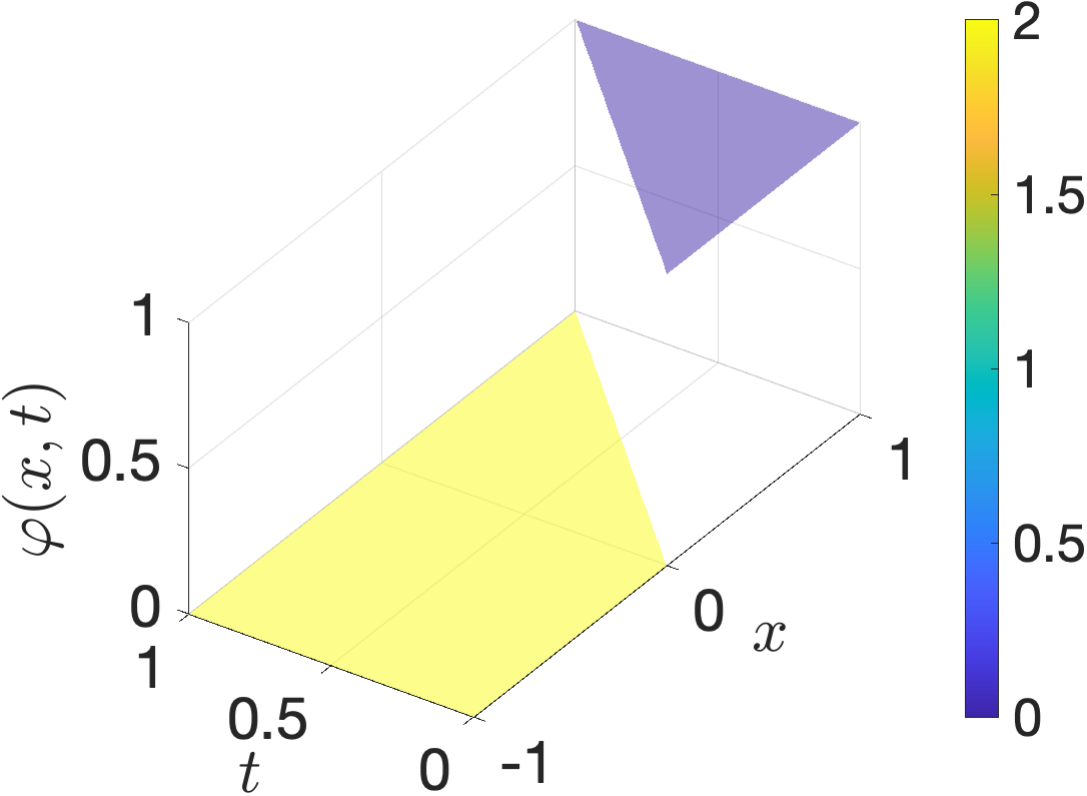}
\caption{$\hat{u}(x,t,\varphi(x,t))$ in one run}
\end{subfigure}
\caption{Numerical results for the Burgers' equation \eqref{sec-Exp2-1d-Burgers} with the initial guess of shock speed set to $\hat{s} = -5$.}
\label{fig-Burgers-Inv-PwC-i-5}
\end{figure}

\begin{figure}[t!]
\centering
\begin{subfigure}{0.32\textwidth}
\centering
\includegraphics[width=\textwidth]{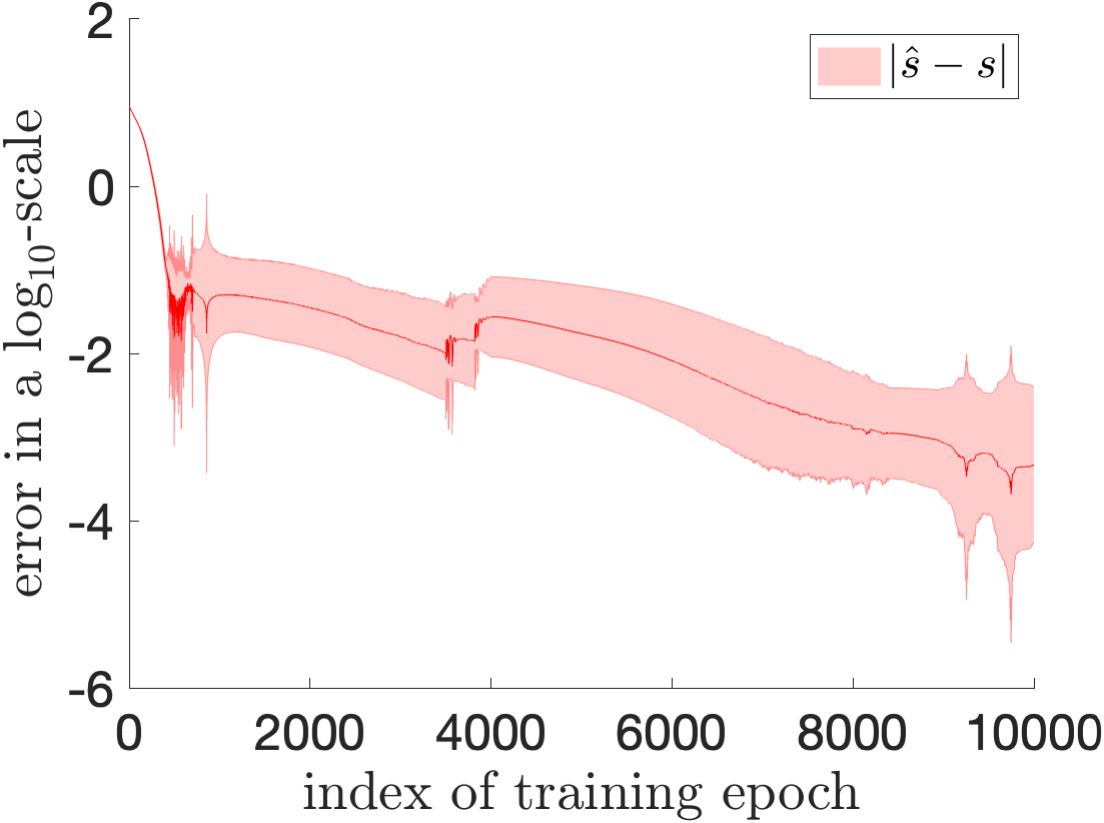}
\caption{$\log_{10}^{|\hat{s} -s |}$ (mean $\pm$ deviation)}
\end{subfigure}
\hspace{0.15cm}
\begin{subfigure}{0.32\textwidth}
\centering
\includegraphics[width=\textwidth]{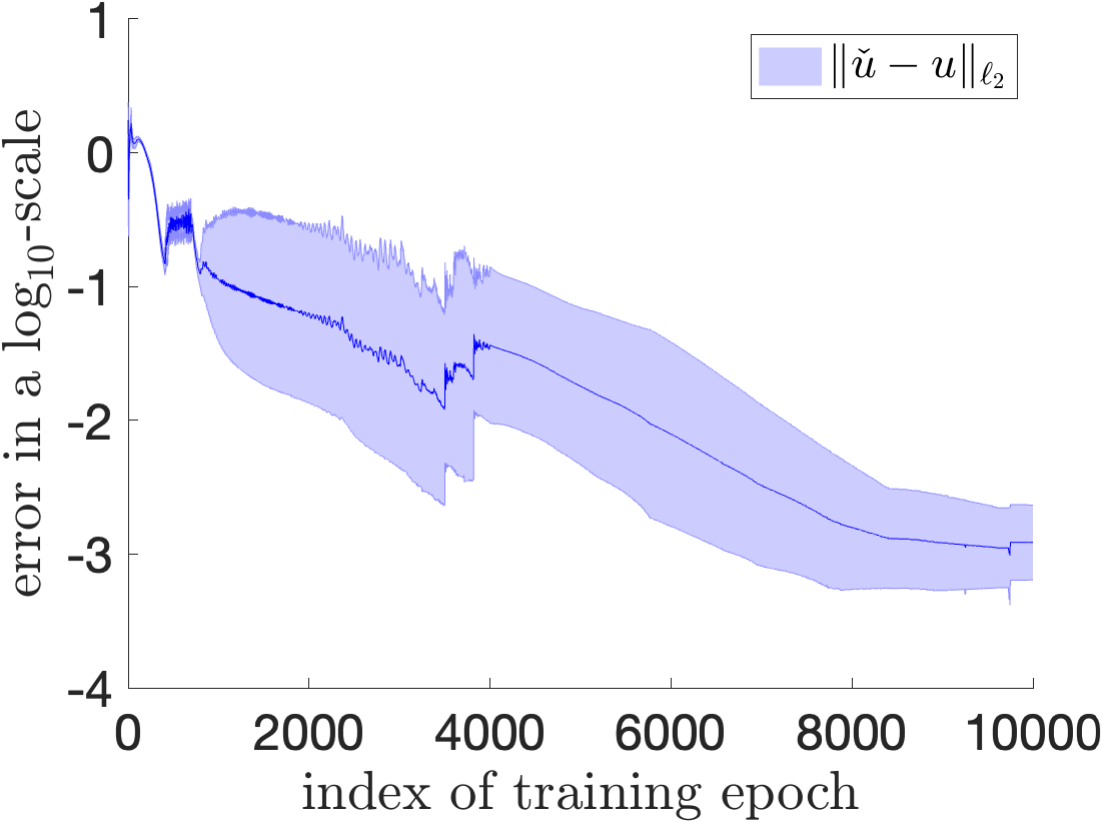}
\caption{test loss (mean $\pm$ deviation)}
\end{subfigure}
\hspace{0.15cm}
\begin{subfigure}{0.32\textwidth}
\centering
\includegraphics[width=\textwidth]{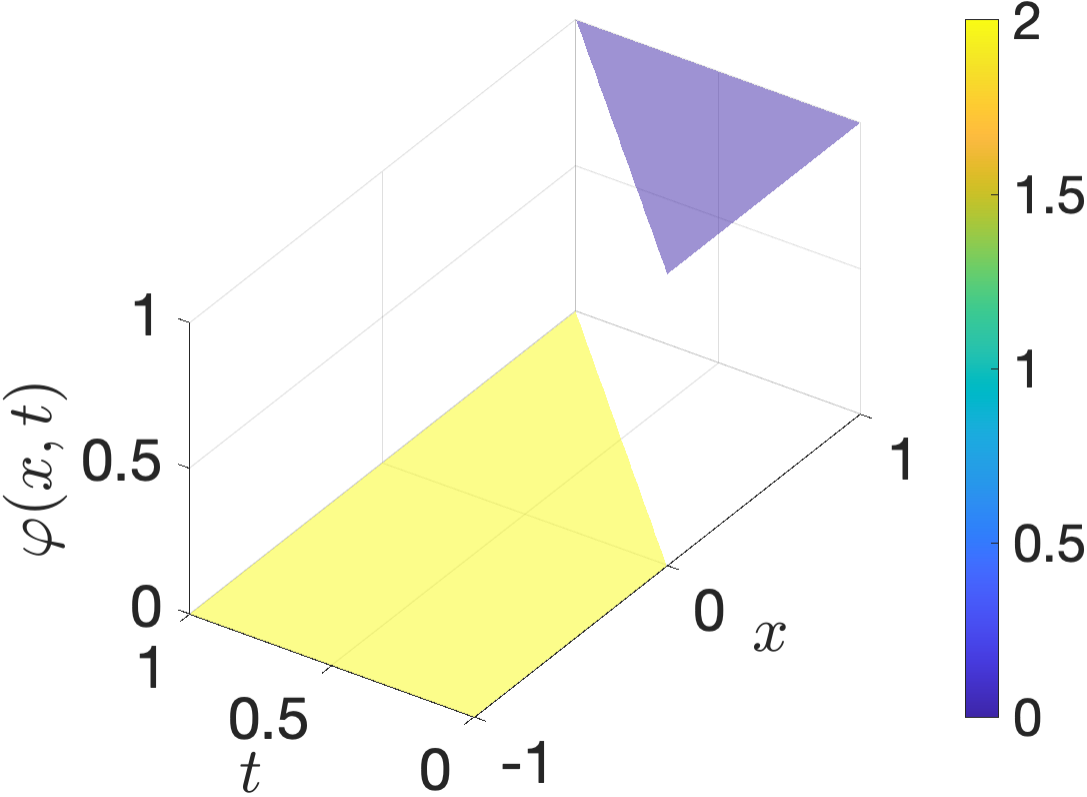}
\caption{$\hat{u}(x,t,\varphi(x,t))$ in one run}
\end{subfigure}
\caption{Numerical results for the Burgers' equation \eqref{Exp2-Eqns-u-exact} with the initial guess of shock speed set to $\hat{s} = 10$.}
\label{fig-Burgers-Inv-PwC-i10}
\end{figure}

The initial guess of shock speed $\hat{s}$ is set to $-5$ and $10$ respectively, while the true value $s = 1$ as discussed in Section \ref{sec-Exp2-1d-Burgers}.  Algorithm \ref{Algorithm-LPLM-inverse} is evaluated over 6 independent runs, with the mean and standard deviation of parameter error $|\hat{s} - s|$ and testing loss $\| \check{u} - u \|_{\ell_2}$ being reported in \autoref{fig-Burgers-Inv-PwC-i-5} and \autoref{fig-Burgers-Inv-PwC-i10}, as well as our trained neural network solution $\hat{u}(x,t,\varphi(x,t))$ in a typical simulation. Predicted shock speeds are in good agreement with the exact value in both cases, which validates the flexibility of our learning approach. It is also noteworthy that the approximation accuracy in this case is slightly lower compared to that in Section \ref{sec-Exp2-1d-Burgers} (see \autoref{table-relative-errors}), as the discrepancy between $\hat{s}$ and $s$ would compromise the training accuracy of our neural network solution.

%%---------------------------------------%

%%---------------------------------------%

%%%%%%%%%%%%%%%%%%%%%%%%%%%%%%%%%%%%%%%%%%%%%%%%%%%%%
\subsubsection{Curved Shock Trajectory}\label{sec-ExpSM2-1d-iBurgers-Curved}

Next, we consider the following inviscid Burgers' equation 
\begin{equation}
	\left\{
	\begin{array}{ll}
		\displaystyle \partial_t u(x,t) + u(x,t) \partial_x u(x,t) = 0, \ &\ \ \ \textnormal{for}\ (x,t)\in \Omega = (0, 2)\times (0, 0.5],\\
		\displaystyle u_0(x) = \left\{
		\begin{array}{l}
			4x,\\
			0, 
		\end{array}\right.
		\ & \ \ \,
		\begin{array}{c}
			\textnormal{for} \ \ 0\leq x<1,\\
			\textnormal{for} \ \ 1 \leq x \leq 2, \\
		\end{array} \\
		\displaystyle u(0,t) = u(2,t) = 0,\ &\ \ \ \textnormal{for}\ \ t\in (0, 0.5],
	\end{array}\right.
	\label{ExpSM2-Eqns-u-exact}
\end{equation}
where the shock wave, departing from the initial discontinuity at $\gamma(0)=x_0=1$, propagates along a curved trajectory $\gamma(t)=\sqrt{1+4t}$. Note that, in contrast to the numerical example reported in Section \ref{sec-Exp8-1d-iBurgers}, the shock speed $s(t)$ in \eqref{ExpSM2-Eqns-u-exact} varies over time and is thus discretized on an equidistant mesh with meshwidth $h=1/50$. 

We choose $\hat{\gamma}(t) = 0.5t+\bar{\gamma}(t)$ as our approximated shock curve, with a heuristic shock trajectory $\gamma(t)\approx 0.5t$ being used as the hypothesis model, which is determined by numerically solving \eqref{ODE-shock-speed-inverse}, i.e.,
\begingroup
\renewcommand*{\arraystretch}{1.2}
\begin{equation}
	\left\{
	\begin{array}{l}
		\displaystyle \frac{d \bar{\gamma}(t)}{dt} = \hat{s}(t)-0.5, \ \ \ \textnormal{for}\ \ t\in (0,0.5],\\
		\bar{\gamma}(0) = 1, 
	\end{array}\right.
	\label{ExpSM2-ODE}
\end{equation}
\endgroup
through a fourth-order Runge-Kutta method on the equidistant mesh $\{ t_i = i h \}_{i=0}^{25}$ with the initial guess of instantaneous values set to $\{ \hat{s}_i = \hat{s}(t_i) = 0 \}_{i=1}^{26}$. Accordingly, the augmented variable can be constructed as
\begin{equation*}
	\varphi(x,t) = H(x - \hat{\gamma}(t)),
\end{equation*} 
and Algorithm \ref{Algorithm-LPLM-inverse} is employed to train the neural network solution and infer the shock speed. Here, an additional penalty coefficient $\beta = 100$ is assigned to the interior loss term.

\begin{figure}[t!]
	\centering
	\begin{subfigure}{0.32\textwidth}
		\centering
		\includegraphics[width=\textwidth]{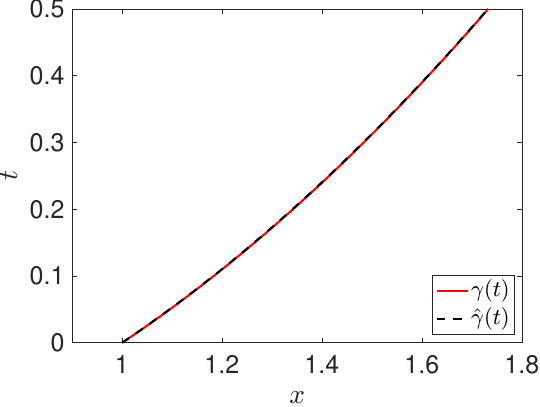}
		\caption{$\gamma(t)$ and $\hat{\gamma}(t)$}
		\label{expSM2-fig-a}
	\end{subfigure}
	\hspace{0.1cm}
	\begin{subfigure}{0.32\textwidth}
		\centering
		\includegraphics[width=0.97\textwidth]{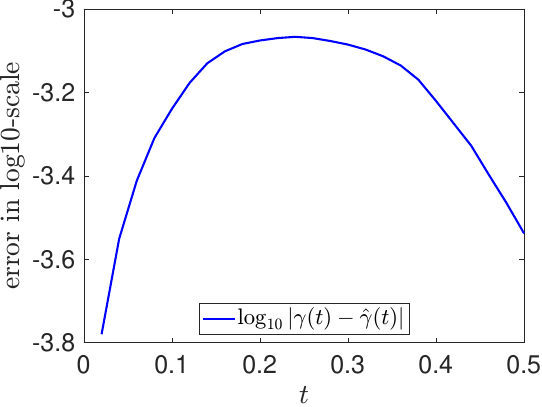}
		\caption{$\log_{10} |\gamma(t) - \hat{\gamma}(t)| $}
		\label{expSM2-fig-b}
	\end{subfigure}
	\hspace{0.1cm}
	\begin{subfigure}{0.32\textwidth}
		\centering
		\includegraphics[width=0.98\textwidth]{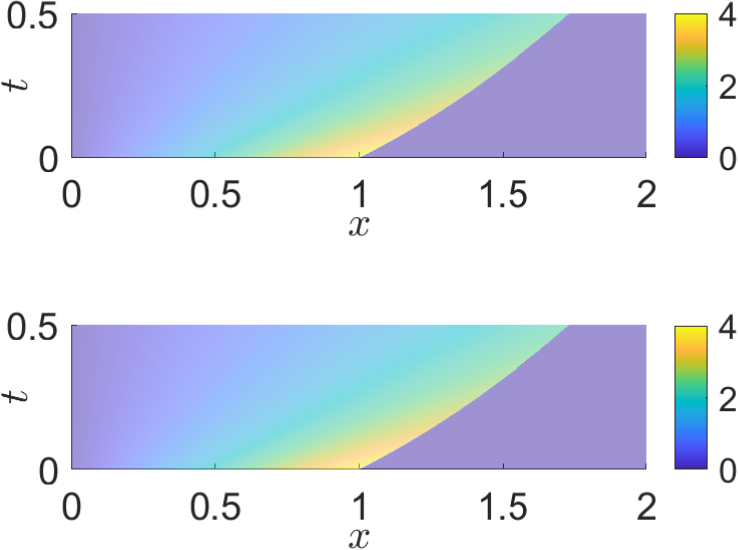}
		\caption{$u(x, t)$ and $\check{u}(x, t)$}
		\label{expSM2-fig-c}
	\end{subfigure}
	\caption{Numerical results for the Burgers' equation \eqref{ExpSM2-Eqns-u-exact} with a curved shock trajectory $\gamma(t)=\sqrt{1+4t}$.}
	\label{fig-Burgers-Inv-PwF}
\end{figure}

Numerical results obtained through our proposed learning approach are displayed in \autoref{fig-Burgers-Inv-PwF}, in which the shock curve and the discontinuous solution are recovered with commendable accuracy. Though the error arising from the numerical reconstruction of shock curve would further affect the accuracy of the trained neural network solution, the relative error reported in \autoref{table-relative-errors} maintains reasonably good performance.

%%---------------------------------------%
%%---------------------------------------%

%%%%%%%%%%%%%%%%%%%%%%%%%%%%%%%%%%%%%%%%%%%%%%%%%%%%%
%%%%%%%%%%%%%%%%%%%%%%%%%%%%%%%%%%%%%%%%%%%%%%%%%%%%%

%%%%%%%%%%%%%%%%%%%%%%%%%%%%%%%%%%%%%%%%%%%%%%%%%%%%%
%%%%%%%%%%%%%%%%%%%%%%%%%%%%%%%%%%%%%%%%%%%%%%%%%%%%%
\section{Conclusions} \label{Section-Conclusion}

In this paper, we present a lift-and-embed learning method for solving scalar hyperbolic equations with discontinuous solutions, accompanied by extensive numerical experiments that validate the effectiveness and flexibility of our method. By including an augmented variable to embed our solution ansatz into a one-order higher-dimensional space, the Rankine-Hugoniot jump condition is expressed at separate collocation points, enabling not only the effective management of multivalued functions but also the use of smooth neural networks to reconstruct discontinuous solutions. Moreover, our method facilitate the numerical solution of both linear and quasi-linear problems within a unified learning framework, accommodating scenarios with known and unknown locations of discontinuities. It is also noteworthy that our modified equations are defined on piecewise surfaces within the elevated-dimensional space, rather than fulfilling the entire domain, hence the sampling of collocation points remains unchanged despite the increased dimensionality. Experimental studies on various problems are reported, which numerically demonstrate that our proposed method is capable of resolving discontinuities without spurious smearing and oscillations. 

For future work, we anticipate that our framework could be extended to systems of hyperbolic conservation laws \cite{godlewski2013numerical}, offering a promising avenue for resolving discontinuous solutions and detecting shock locations. Additionally, conducting a comprehensive error analysis is of critical importance and will be a key focus of future investigations. It is also noteworthy that, with recent advancements in physics-informed operator learning \cite{lin2023operator,rosofsky2023applications}, applying our lift-and-embed method to infer solution operators for hyperbolic problems offers a compelling research direction.

%%%%%%%%%%%%%%%%%%%%%%%%%%%%%%%%%%%%%%%%%%%%%%%%%%%%%
%%%%%%%%%%%%%%%%%%%%%%%%%%%%%%%%%%%%%%%%%%%%%%%%%%%%% 

%%%%%%%%%%%%%%%%%%%%%%%%%%%%%%%%%%%%%%%%%%%%%%%%%%%%%
%%%%%%%%%%%%%%%%%%%%%%%%%%%%%%%%%%%%%%%%%%%%%%%%%%%%% 
\begin{acknowledgements}
Q. Sun is supported in part by National Natural Science Foundation of China (grant 12201465), Science and Technology Commission of Shanghai Municipality (grant 23JC1400502), and Shanghai Municipal Science and Technology Major Project (No. 2021SHZDZX0100). 
X. Xu is supported in part by National Natural Science Foundation of China (grants 12071350 and 12331015). This research was conducted using the computational resources and services of the HPC Center at the School of Mathematical Sciences, Tongji University.
\end{acknowledgements}
%%%%%%%%%%%%%%%%%%%%%%%%%%%%%%%%%%%%%%%%%%%%%%%%%%%%%
%%%%%%%%%%%%%%%%%%%%%%%%%%%%%%%%%%%%%%%%%%%%%%%%%%%%% 
\noindent {\bf \small Data Availability } \small Enquiries about data availability should be directed to the authors. All source codes are publicly accessible on GitHub (\url{https://github.com/q1sun/LELM_HyperbolicPDE}).


\section*{Declarations}

{\bf\small Conflict of interest } The authors have no conflict of interest.

%\section*{Conflict of interest}
%The authors declare that they have no conflict of interest.
%%%%%%%%%%%%%%%%%%%%%%%%%%%%%%%%%%%%%%%%%%%%%%%%%%%%%
%%%%%%%%%%%%%%%%%%%%%%%%%%%%%%%%%%%%%%%%%%%%%%%%%%%%% 

\bibliographystyle{spmpsci}
\bibliography{refs}
%%%%%%%%%%%%%%%%%%%%%%%%%%%%%%%%%%%%%%%%%%%%%%%%%%%%%
%%%%%%%%%%%%%%%%%%%%%%%%%%%%%%%%%%%%%%%%%%%%%%%%%%%%% 

\end{document}